\documentclass[3p, number,sort&compress,times,fleqn]{elsarticle}      

\graphicspath{{./figs/}}
\usepackage{url}
\usepackage{booktabs}
\usepackage{multirow}
\usepackage[table]{xcolor}
\usepackage{./style/csml}

\usepackage[colorlinks = true,
linkcolor = blue,
urlcolor  = blue,
citecolor = blue,
anchorcolor = blue]{hyperref}

\begin{document}
	
\begin{frontmatter}

\title{Stochastic PDE representation of random fields for large-scale Gaussian process regression and statistical finite element analysis\tnotemark[1]}

\author{Kim Jie Koh}
\author{Fehmi Cirak\corref{cor1}}
\ead{f.cirak@eng.cam.ac.uk}

\tnotetext[1]{This paper is dedicated to honouring the lifetime achievements of Thomas J.R. Hughes (Tom), whose numerous prescient contributions to computational science, engineering, and mathematics have offered a constant source of inspiration and enjoyment for us, as well as an entire generation or two (so far). The second author is deeply indebted to Tom for his invaluable mentorship and continued fellowship.}

\cortext[cor1]{Corresponding author}

\address{Department of Engineering, University of Cambridge, Cambridge, CB2 1PZ, UK }

\begin{abstract}
The efficient representation of random fields on geometrically complex domains is crucial for Bayesian modelling in engineering and machine learning, including Gaussian process regression and statistical finite element analysis. Today's prevalent random field representations are either intended for unbounded domains or are too restrictive in terms of possible field properties. Because of these limitations, new techniques leveraging the historically established link between stochastic PDEs (SPDEs) and random fields have been gaining interest in the statistics and engineering literature. The SPDE representation is especially appealing for engineering applications with complex geometries which already have a finite element discretisation for solving the physical conservation equations. In contrast to the dense covariance matrix of a random field, its inverse, the precision matrix, is usually sparse and equal to the stiffness matrix of an elliptic SPDE. In this paper, we use the SPDE representation to develop a scalable framework for large-scale statistical finite element analysis and Gaussian process (GP) regression on geometrically complex domains. The statistical finite element method (statFEM)  introduced by Girolami et al.~(2022) is a novel approach for synthesising measurement data and finite element models. In both statFEM and GP regression, we use the SPDE formulation to obtain the relevant prior probability densities with a sparse precision matrix. The properties of the priors are governed by the parameters and possibly fractional order of the SPDE so that we can model on bounded domains and manifolds anisotropic, non-stationary random fields with arbitrary smoothness. We use for assembling the sparse precision matrix the same finite element mesh used for solving the physical conservation equations. The observation models for statFEM and GP regression are such that the posterior probability densities are Gaussians with a closed-form mean and precision. The expressions for the mean vector and the precision matrix do not contain dense matrices and can be evaluated using only sparse matrix operations. We demonstrate the versatility of the proposed framework and its convergence properties with one and two-dimensional Poisson and thin-shell examples.
 \end{abstract}
	
\begin{keyword}
Bayesian modelling, Gaussian processes,  statistical finite elements, physics-informed priors, stochastic PDEs, fractional PDEs
\end{keyword}

\end{frontmatter}


%
\section{Introduction \label{sec:introduction}}
%
%
\subsection{Motivation \label{sec:introductionMotivation}}
%
Gaussian process-based Bayesian models and related techniques play a crucial role in probabilistic engineering and machine learning~\cite{williams2006gaussian, santner2003design, sobester2008engineering}. They provide the means to consistently blend random fields, or processes, with varying levels of uncertainty and to take into account any information or constraints available about those fields. In Bayesian statistics, all uncertainties (both epistemic and aleatoric) are represented using random fields or variables and their respective probability measures. The prior probability measures are updated according to an assumed data-generating process, yielding the likelihood measure, in light of the observation data~\cite{stuart2010inverse,ghattas2021learning}. The resulting posterior probability measure provides the expectation of the random variables of interest, and the spread of the posterior is indicative of the confidence we can place on those values. In GP-based models, the likelihood and the prior measures are all Gaussians, so that the posterior measure is a Gaussian with a closed-form mean and covariance, hugely simplifying its computation.  In physical systems, random variables are subject to certain constraints. For instance, the constitutive parameters must satisfy certain positivity and symmetry conditions, and the solution field must satisfy conservation equations expressed in the form of PDEs. Conceptually, in Bayesian approaches, these constraints can be enforced by modifying either the prior or the data-generating model, which makes them well-suited for combining with established models and methods from computational mechanics. However, the application of Bayesian approaches to engineering problems involving the solution of PDEs  is hampered by their vast computing requirements and the efficient description of random fields on domains with complex geometries, such as engineering structures consisting of solids, beams and shells. As we will demonstrate in this paper, both problems can be overcome by the stochastic PDE representation of random fields and their discretisation using standard finite elements.
  
%
\subsection{Related research \label{sec:introductionRelated}}
%
The discretisation of a Gaussian random field~$s(\vec x) \in \mathbb R$  with \mbox{$\vec x \in \mathbb R^d$} leads to a Gaussian random vector \mbox{$\vec s \in \mathbb R^{n}$}, where each component $s_i$ corresponds to a point $\vec x_i \in \mathbb R^d$, where $d \in \{1, \, 2, \, 3 \}$ and $i=1, \dotsc, n$. The multivariate Gaussian probability density of~$\vec s$ has the mean $\overline{\vec s} \in \mathbb R^n$ and the covariance matrix $\vec C_s \in \mathbb R^{n \times n}$. The covariance matrix is usually obtained by pointwise evaluating a covariance function such as the Mat\'ern kernel. Although Gaussian probability densities are conspicuously easy to handle analytically, their numerical treatment becomes, unfortunately, increasingly challenging for larger problems. This is because covariance matrices are dense and full-rank, so that their storage and inversion have $O(n^2)$ and $O(n^3)$ complexity, respectively, making it impossible to consider random vectors with more than a few thousand components. This poor scalability is a significant limitation when Gaussian processes are combined with finite element discretised random fields, which usually have at least several tens of thousands of finite element degrees of freedom.  As an additional difficulty, the covariance functions commonly used in statistics, including the Mat\'ern kernel,  are restricted to Euclidean domains and are inadequate for random fields on non-Euclidean  domains, like a shell structure or a truss structure consisting of several members. These and other limitations of covariance functions and Gaussian processes can be efficiently dealt with by switching to a stochastic PDE representation of random Mat\'ern fields. 

\begin{figure} 
	\centering
	\subfloat[][Isotropic and stationary \label{fig:bunnyStationarySample}]
	{
		\includegraphics[width=0.275\textwidth]{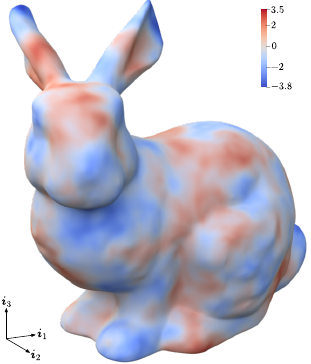}
	}
	\hfil
	\subfloat[][Anisotropic and stationary \label{fig:bunnyAnisotropicSample}]
	{
		\includegraphics[width=0.275\textwidth]{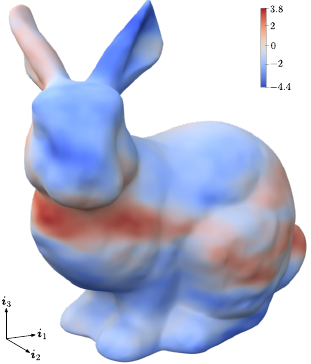}
	}
	\hfil
	\subfloat[][Isotropic and non-stationary  \label{fig:bunnyNonStationarySample}]
	{
		\includegraphics[width=0.275\textwidth]{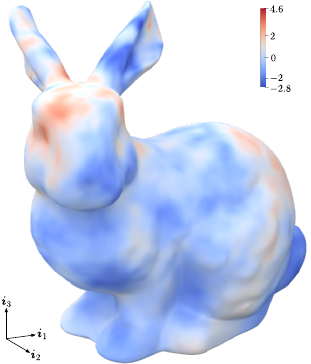}
	} 
	\caption{Samples from three different Gaussian random fields on the surface of the Stanford bunny. The covariance function for each random field has the properties mentioned in the respective caption. Notice in (b) the different length-scales along the~$\vec i_3$  axis and orthogonal to it and in (c) the decrease in  length-scale along the~$\vec i_3$ axis.  For visualisation of the respective covariance functions and further details see Figures~\ref{fig:bunnyStationary}, \ref{fig:bunnyAnisotropic} and~\ref{fig:bunnyNonStationary}. The surface is discretised with $1271808$ linear triangular elements and has the bounding box~\mbox{$(15.5, \, 12, \, 15.3)$}.
	  \label{fig:bunniesIntro}}
\end{figure}

The link between Gaussian processes, specifically Mat\'ern fields, and stochastic PDEs has been known since Whittle~\cite[Sect. 9]{whittle1954stationary} and was more recently reintroduced to statistics by Lindgren et al.~\cite{lindgren2011explicit}.  It can be shown that the solution of a second-order elliptic PDE with a fractional-order exponent and a Gaussian white noise on the right-hand side is equal to a Gaussian process with a Mat\'ern covariance function and a zero mean.  In~\cite{lindgren2011explicit}, it is shown how the stiffness matrix of the discretised stochastic PDE corresponds to the inverse covariance matrix~$\vec C_s^{-1}$, i.e.\ the precision matrix~$\vec Q_s=\vec C_s^{-1}$, of the multivariate Gaussian density. The precision matrix~$\vec Q_s$ is sparse, especially when the fractional order of the SPDE is within a specific range~\cite{bolin2020rational}. Indeed, the precision matrix is always sparse when the PDE is non-fractional and has a finite power. Intuitively, a zero entry of the precision matrix expresses the conditional independence between two components of a Markov random vector when all the other components are known, which is, in fact, a plausible assumption for most physical fields, see also~\cite{rue2005gaussian,povala2022variational}.  As an additional benefit, the SPDE representation makes it obvious how to generalise the Mat\'ern covariance function to non-Euclidean domains and to non-stationary (i.e., non-homogeneous) and anisotropic random fields. See Figure~\ref{fig:bunniesIntro} for different random fields on the surface of the Stanford bunny obtained by solving an SPDE with a Gaussian white noise as the right-hand side.

The seminal paper by Lindgren et al.~\cite{lindgren2011explicit}, combined with the ubiquity of Gaussian process-based models in probabilistic engineering and machine learning, has lately generated immense interest in the SPDE representation of Mat\'ern fields. The SPDE formulation is especially appealing for random fields $s(\vec x) \in \mathbb R$ with $\vec x \in \Omega \subset \mathbb R^d$ on low-dimensional domains with $d \in \{1, 2, 3\}$, such as the diffusivity field on a domain of interest  $\Omega$, which may be a manifold, like a surface embedded in $\mathbb R^3$. In engineering, the domain~$\Omega$ is usually given as a computer-aided design (CAD) model and is discretised with a finite element mesh for solving the governing equations. Hence, given that there is already a mesh, the SPDE representation of random fields is especially appealing for engineering applications. Recently, Zhang et al.~\cite{zhang2021stochastic} and Wang et al.~\cite{wang2022stochastic} considered the non-fractional, second-order SPDE with a linear exponent for representing geometric uncertainties on two-manifolds in finite element discretised probabilistic forward problems. As discussed in  Chen et al.~\cite{chen2022spatially},  the SPDE formulation is not restricted to Gaussian random fields as it can represent, after a suitable transformation, for instance, random constitutive parameters subject to constraints. In data assimilation, a non-fractional SPDE formulation has been utilised by Rouse et al.~\cite{paul2022probabilistic} and Poot et al.~\cite{poot2023bayesian}. Before the mentioned works,  Bui-Thanh et al.~\cite{bui2013computational} used the inverse of the finite element discretised Laplace operator as the prior covariance in  Bayesian inverse problems. The mathematical analysis of  inverse problems with the inverse fractional Laplace operator as the prior covariance was discussed by Stuart~\cite{stuart2010inverse}. However, most papers on the SPDE formulation are currently from geostatistics and primarily concerned with spatial random fields over two and three-dimensional domains; see the recent review~\cite{lindgren2022spde}. Most of the applications mentioned so far consider SPDEs with integer exponents. Yet, to obtain Mat\'ern fields with arbitrary smoothness, it is necessary to allow for fractional exponents. Among various definitions of fractional PDEs~\cite{kwasnicki2017ten,lischke2020fractional}, the most relevant for the present paper is the rational series approximation of the fractional PDE operator introduced by Bolin and Kirchner~\cite{bolin2020rational} and Haziranov et al.~\cite{harizanov2018optimal}. The principal idea in~\cite{bolin2020rational} is to determine first a series approximation of the function~$x^{-\beta}$ with respect to $x$, where $x, \, \beta \in \mathbb R$, and to subsequently replace $x$ with the differential operator of the fractional PDE. This approach is rooted in the well-known spectral mapping theorem. Although any series expansion could be used in principle, barycentric rational interpolation is the most robust and efficient~\cite{berrut2004barycentric,hofreither2021algorithm}. Instead of the rational series expansion, it is possible to use the integral representation of fractional PDEs and to evaluate the integral numerically~\cite{bolin2020numerical}. We refer to Higham~\cite{higham2008functions} for related approximation techniques for fractional matrix functions. 

The SPDE representation of random fields is one of the many techniques to improve the scalability of Gaussian process-based models. Especially in machine learning, primarily concerned with inferring a random field~$s(\vec x) \in \mathbb R$ from datapoints~$\{( \vec x_i, \, s_i ) \}_{i=1}^n $,  a wide range of approximation techniques have been proposed, see Rasmussen and Williams~\cite[Ch.\ 8]{williams2006gaussian} and reviews~\cite{hensman2013gaussian, liu2020gaussian}. These techniques usually aim to either improve the sparsity or reduce the size or rank of the covariance matrix $\vec C_s \in \mathbb R^{n \times n}$.  For instance,  in so-called sparse Gaussian processes, a limited number of inducing points are introduced to obtain a low-rank approximation of the original covariance matrix~\cite{quinonero2005unifying}. Other approximation techniques from geostatistics are more suitable for spatial Gaussian processes over~$\Omega \subseteq \mathbb R^d$ with~$d \in \{ 1, \, 2, \, 3\}$, see the review~\cite{heaton2019case}. The conventional covariance approximation techniques from machine learning and geostatistics do not  require a discretisation of the problem domain~$\Omega$ and are usually only suitable for problems defined over unbounded domains, i.e.~$\Omega \equiv \mathbb R^d$. More recently, with the continued convergence of machine learning and computational mechanics, truncated-series expansions of the covariance function using orthogonal global basis functions are gaining interest. For instance, Solin and S{\"a}rkk{\"a}~\cite{solin2020hilbert} use  as basis functions the eigenfunctions of the  finite element discretised Laplace operator on~$\Omega$. This approach is similar to the truncated-series expansion methods from computational mechanics, including the Karhunen--Lo{\`e}ve method~\cite{betz2014numerical,marzouk2009dimensionality} and the spectral representation method~\cite{shinozuka1972digital, shinozuka1991simulation}, which are based on some other orthogonal basis functions. All the covariance approximation techniques mentioned in this paragraph are restricted to stationary (homogeneous) and isotropic covariance functions, i.e.\ the covariance function depends only on the distance between two points, and exclusively focus on approximating the dense covariance matrix. In addition, the techniques from machine learning are primarily intended for unbounded Euclidean domains, while truncated-series expansions can, in principle, be applied to bounded domains, see e.g.~\cite{uribe2020bayesian,gulian2022gaussian}. 

As noted, GP-based Bayesian models are particularly attractive for consistently blending data and information from different physical models, empirical knowledge and observational data. In contrast to machine learning, which relies mainly on observational data, physical models and empirical knowledge play a crucial role in engineering. Physical models are expressed as PDEs with boundary and initial conditions, which can all be taken into account by choosing a suitable prior probability density or observation model, i.e.\ likelihood, see the comprehensive review~\cite{swiler2020survey}. In the case of GP priors, their mean and covariance can be tailored so that any constraints are precisely or approximately satisfied and empirical knowledge is considered. Subsequent conditioning of the GP-based model on the available observational data yields a posterior density which implicitly takes into account the data and the constraints imposed by the prior. In engineering, this approach has been first rigorously explored by Raissi et al.~\cite{raissi2017machine,raissi2018hidden} for generating new prior covariance functions by applying the continuous PDE operators to standard covariance functions, like the squared-exponential kernel; related earlier and recent works include~\cite{graepel2003solving,sarkka2011linear,owhadi2015,chen2021solving,chen2022apik}. These methods are intended for supervised learning on domains with simple geometries and are framed as a replacement for current engineering analysis techniques, i.e.\ finite element analysis. Moreover, due to their reliance on standard covariance functions from machine learning, they have problems with scalability and consideration of non-Euclidean domains. 

Recently, Girolami et al.~\cite{girolami2021statistical} introduced the statistical finite element method (statFEM), which uses a conventional probabilistic forward problem, see e.g.~\cite{ghanem1991stochastic,sudret2000stochastic,matthies2005galerkin}, to obtain the prior density. Consequently, statFEM is inherently well-suited for domains with complex geometries. The statFEM observation model is inspired by the seminal Kennedy and O'Hagan~\cite{kennedy2001bayesian}  paper on Bayesian calibration, which introduces a data-driven GP-based framework for blending data from a black-box simulator and observations. Since its inception, numerous extensions and applications of the Kennedy and O'Hagan framework have been proposed, far too many to discuss here, see e.g.~\cite{higdon2004combining,xiong2009better,ling2014selection,yang2019physics,maupin2020model,jiang2020sequential}. A crucial component of statFEM is a random discrepancy, or inadequacy, term taking into account the mismatch between the finite element model and the actual system. In practice, such a mismatch is inevitable because of the very assumptions and simplifications necessary in creating a numerical model. The conditioning of the statFEM model on the observation data yields a posterior that consistently blends the finite element prior with the observation data. Any hyperparameters of the statistical model are determined by maximising the marginal likelihood as in standard GP regression.  The finite element prior and posterior probability densities in statFEM are Gaussians and are inevitably high-dimensional, given that finite element models can have several hundred thousand of unknowns.

%
\subsection{Contributions \label{sec:introductionContrib}}
%
We introduce a scalable framework for large-scale GP regression and statistical finite element analysis on domains with complex geometries and problems with generalised random fields. As will be discussed, statFEM reduces to standard GP regression after certain simplifying assumptions are introduced. According to the observation model underpinning statFEM, the observed data vector is decomposed into a finite element, a model inadequacy and a noise component.  The finite element and inadequacy components are independent Gaussian random fields and the noise component is independent and identically distributed. The mentioned finite element component is the solution of the physical governing equations, i.e. conservation equations, and should not be confused with the solution of the SPDE for representing random fields. In statFEM, the priors for the finite element and inadequacy components are chosen as Gaussian processes so that the posterior obtained is a Gaussian process with a closed-form mean and covariance. We determine the finite element prior by solving a conventional stochastic forward problem. In the present paper, we consider for the sake of illustration only linear governing equations and assume that only the source, or forcing term, is a Gaussian process  leading to a  finite element solution which is a Gaussian process. In the case of nonlinear governing equations or a PDE operator with random parameters, like random diffusivity, we can use a first-order perturbation to approximate the random solution as a Gaussian process~\cite{girolami2021statistical}. 

The key novelty of our approach is to express the two statFEM random fields, namely the source and model inadequacy, using their SPDE representations, assuming that both are generalised Mat\'ern random fields. The two independent SPDEs are discretised with the same finite element approach and mesh as the physical governing equations. In the presented examples, we use for finite element discretisation either conventional Lagrange or isogeometric subdivision basis functions~\cite{Hughes:1987aa,hughes2005isogeometric,Cirak:2000aa}. As mentioned, the stiffness matrix of the SPDE operator is equal to the precision matrix of the respective generalised Mat\'ern field. To discretise SPDE operators with arbitrary power, we first decompose the exponent into an integer and a fractional part yielding two SPDEs which are discretised in turn.  The solution of the SPDE with the integer exponent is obtained using the recursion technique suggested by Lindgren et al.~\cite{lindgren2011explicit}. This approach can be interpreted as a mixed finite element discretisation of the SPDE and circumvents the need for smooth basis functions. The solution of the SPDE with the integer exponent serves as the source term for the SPDE with the fractional exponent, which we solve by following the rational approximation technique proposed by Bolin and Kirchner~\cite{bolin2020rational}.  Departing from~\cite{ bolin2020rational}, we use the barycentric rational interpolation algorithm of Hofreither~\cite{hofreither2021algorithm}  for stably computing the coefficients in the rational series expansion. For implementational convenience, we use the same polynomial degree in the numerator and denominator of the rational series approximation. In the sketched approach, the sparse precision matrix of the fractional SPDE operator is obtained by repeated multiplication of standard stiffness matrices. Consequently, it has a slightly larger memory footprint than standard stiffness matrices but is much smaller than the corresponding usually dense covariance matrix.

We introduce anisotropic and non-stationary random fields by slightly altering the SPDE operator without impacting the efficiency of the solution process.  In the case of random fields on shells, we replace the Laplace operator in the SPDE with the Laplace-Beltrami operator. The statFEM posterior mean vector and  precision matrix are  determined using only sparse matrix operations involving the prior precision matrices of the solution and inadequacy fields. With the obtained posterior precision matrix, we can evaluate selected components of the respective covariance matrix by factorising the precision matrix once and solving for different right-hand sides. 

%
\subsection{Overview \label{sec:introductionOverview}}
%
The rest of this paper is organised as follows. In Section~\ref{sec:matern}, we begin by first reviewing the strong and weak forms of the SPDE formulation of Mat\'ern fields and introduce their finite element discretisation. Subsequently, we briefly outline the generalisation of the SPDE representation to anisotropic and non-stationary random fields and manifolds. We then introduce in Section~\ref{sec:gpRegression} the use of the obtained precision matrix for large-scale GP regression on the finite element mesh corresponding to the SPDE.  We provide closed-form expressions for the multivariate Gaussian posterior and the marginal likelihood needed for learning the hyperparameters of the SPDE. In Section~\ref{sec:statFEM}, this is followed by the generalisation of  large-scale GP regression to statistical finite element analysis. We first discuss the computation of a PDE-informed prior with a sparse precision matrix by solving a conventional probabilistic forward problem. After that, we review the observation model for statFEM and discuss the treatment of the random discrepancy field. We again provide closed-form expressions for the multivariate Gaussian posterior and the marginal likelihood. Finally, in Section~\ref{sec:examples}  we introduce five examples of increasing complexity demonstrating the convergence of the obtained prior and posterior random fields using the SPDE representation of generalised Mat\'ern fields. In particular, we study convergence in terms of mesh refinement, the order of the fractional series expansion, the number of observation points and repeated readings per observation point. Finally, we provide three appendices reviewing rational interpolation, summarising properties of multivariate Gaussian densities and generalising our results to the case of repeated readings.

%
\section{Mat\'ern random fields\label{sec:matern}}
%
%
\subsection{Fractional PDE representation\label{sec:fractionalPDE}}
\begin{figure}
	\centering
	\subfloat[][Mat\'ern covariance function] {
		\includegraphics[height=0.315 \textwidth]{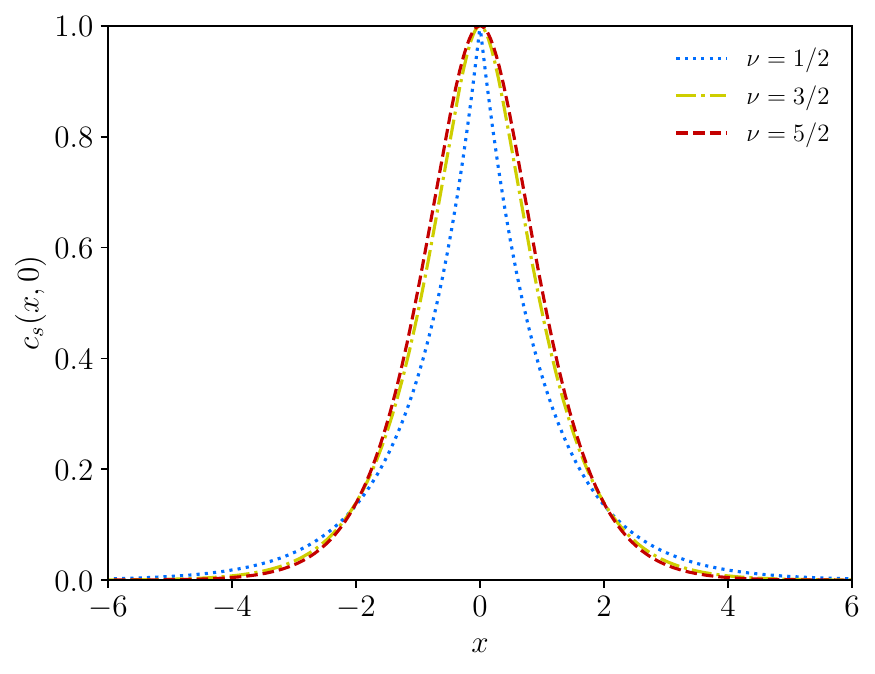} } 
	\hfil
	\subfloat[][Samples drawn from Mat\'ern Gaussian processes] {
		\includegraphics[height=0.315 \textwidth]{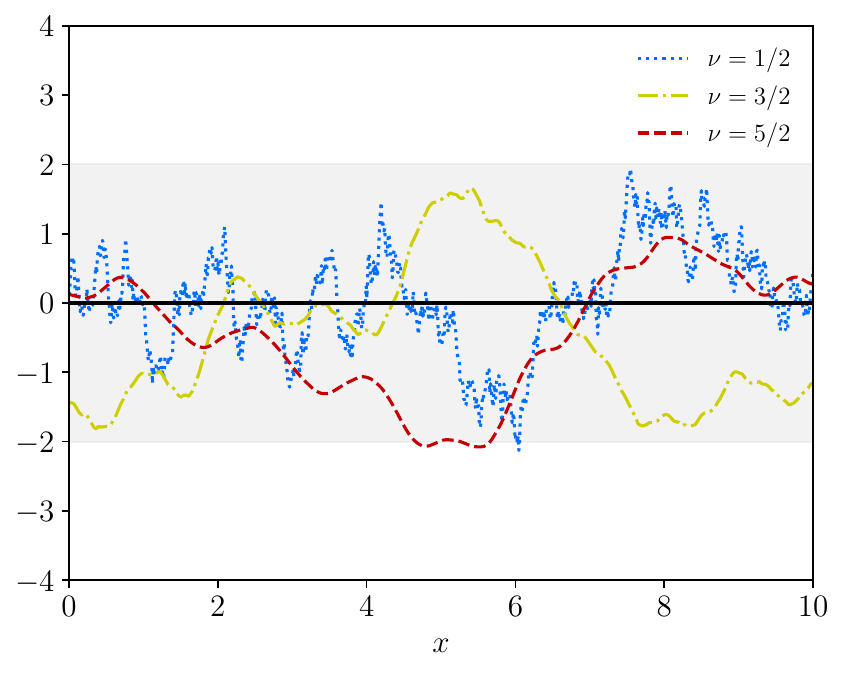} }
	\caption{Mat\'ern covariance function with~$\sigma = 1$,~$\ell = 1$ and~$\nu \in \{ 1/2, \, 3/2, \, 5/2 \}$,  and three samples drawn from respective Mat\'ern Gaussian processes with zero mean. In (b) the black line represents the zero mean and the shaded region an offset of $\pm 2 \sigma$ from the mean. \label{fig:maternPrior}}
\end{figure}

%
We consider in~$\mathbb R^d$, with $d \in \{ 1, \, 2, \, 3\}$, a zero mean Gaussian process, also referred to as a Gaussian random field,
\begin{equation} \label{eq:maternField}
	s(\vec x) \sim \set{GP} \left ( 0, \, c_s(\vec x, \, \vec x')  \right ) \, ,
\end{equation}
with the Mat\'ern covariance function
\begin{equation} 
	c_s(\vec x, \, \vec x') = \cov \left ( s(\vec x), \, s(\vec x') \right ) =   \expect \left [ s(\vec x) s(\vec x') \right ] = \frac{\sigma^2}{2^{\nu -1} \Gamma(\nu)} \left (  \frac{\sqrt{2 \nu}}{\ell}  \| \vec x - \vec x' \| \right )^{\nu}  K_{\nu} \left (  \frac{\sqrt{2 \nu}}{\ell}  \| \vec x - \vec x' \| \right ) \, ,
	\label{eq:maternKernel}
\end{equation}
where~\mbox{$\vec x, \, \vec x' \in \mathbb R^d$} are point coordinates,~$\expect $ is the expectation operator,~\mbox{$\sigma \in \mathbb R^+ $} the standard deviation, \mbox{$\nu \in \mathbb R^+ $} a smoothness parameter,  \mbox{$\ell \in \mathbb R^+ $} a length-scale parameter, $\Gamma$ the Gamma function, and~$K_{\nu}$ the modified Bessel function of the second kind of order~$\nu$. The smoothness parameter~$\nu$ governs the smoothness of the Gaussian process in the mean square sense. In the limit $\nu \rightarrow \infty $ the Mat\'ern covariance function converges to the infinitely smooth squared exponential kernel~\cite{stein1999interpolation}. For~$d = 1$, Figure~\ref{fig:maternPrior} shows the covariance function for smoothness parameters $\nu \in \{ 1/2, \, 3/2, \, 5/2 \}$  and three random functions drawn from the respective Mat\'ern Gaussian processes with zero mean. As can be seen, the smoothness of the covariance function and the samples increase with increasing~$\nu$.

The Mat\'ern  random field~$s (\vec x) $ is the solution of the (stochastic) partial differential equation 
\begin{equation} \label{eq:fracPDE}
	\left (  \kappa^2   - \Delta \right )^\beta s(\vec x) = \frac{1}{\tau} g (\vec x) \,,
\end{equation}
where~$\Delta$ is the Laplace operator,~$g(\vec x)$ the Gaussian white noise process
\begin{equation} \label{eq:whiteNoise}
	g (\vec x) \sim  \set{GP} \left ( 0, \, \delta (\vec x -  \vec x')  \right ) \, ,
\end{equation}
 and the remaining parameters are defined as 
\begin{equation}\label{eq:fracPDEvars}
	\kappa = \frac{\sqrt{2 \nu}}{ \ell} \, , \quad  \beta = \frac{\nu}{2} + \frac{d}{4} \, , \, \quad \tau^2 = \frac{\Gamma (\nu)}{ \sigma^2 \Gamma (\nu+d/2) (4 \pi )^{d/2} \kappa^{2\nu}}   \, .
\end{equation}
The practically most relevant exponents are~$\beta \in \{(1+d)/4, \, (3+d)/4, \,  (5+d)/4 \}$ ~\cite{williams2006gaussian}. In the following we shorten~\eqref{eq:fracPDE} to
\begin{equation} \label{eq:fracPDE2}
	\set L^\beta  s(\vec x) = \frac{1}{\tau}  g (\vec x) \, .
\end{equation}
The order of this fractional partial differential equation depends on the exponent~$\beta$, which, in turn, depends on the choice of the smoothness parameter~$\nu$ and the dimension~$d$. The exponent can take any value~\mbox{$\beta > d/4$}. Furthermore, note that the domain of the partial differential equation is all of~$\mathbb R^d$.

It is expedient to decompose~\eqref{eq:fracPDE2}  into an integer part and a fractional part in the form 
\begin{subequations} \label{eq:fracPDEdec}
\begin{align}
	\set L^\alpha \tilde{s} (\vec x) &=  \frac{1}{\tau} g (\vec x)  \, ,  \label{eq:fracPDEdecA} \\
	\set L^{\beta -\alpha} s (\vec x) &= \tilde{s} (\vec x)  \, , \label{eq:fracPDEdecB}
\end{align} 
\end{subequations}
with the exponent~\mbox{$\alpha \in \mathbb N_{\ge 1} $} chosen as 
\begin{equation}
	\alpha = \max \{ 1, \, \lfloor \beta \rfloor  \} 
\end{equation}
so that the fractional exponent takes the values~\mbox{$\beta -\alpha \in  (d/4 - 1, \, 1) \subset \mathbb R $.} The solution of the integer part~\eqref{eq:fracPDEdecA} is fairly standard and is discussed in Section~\ref{sec:fractionalFE}. The solution of the fractional part~\eqref{eq:fracPDEdecB} requires additional techniques. There are several equivalent definitions of fractional operators~\cite{kwasnicki2017ten}. The most apparent  definition relies on the spectral decomposition of the operator~$\set{L}$ and is according to the spectral mapping theorem given by the expansion
\begin{equation} \label{eq:fracCont}
	\set L^{\beta-\alpha}  s  (\vec x) = \sum_{i=1}^\infty \lambda_i^{\beta - \alpha} \left  ( s  (\vec x), \,  \zeta_i(\vec x) \right  ) \zeta_i (\vec x) \, , 
\end{equation}
where~$ ( \cdot,  \cdot  ) $ denotes the~$L^2$ inner product, and $\zeta_i (\vec x)$ and~$\lambda_i$ are respectively the orthonormal eigenfunctions and eigenvalues satisfying
\begin{equation}
	\set L  \zeta_i(\vec x)  = \lambda_i \zeta_i (\vec x)  \, , \quad   (\zeta_i, \, \zeta_j ) = \delta_{ij} \,.
\end{equation}
Hence, the solution of~\eqref{eq:fracPDEdecB} reads
\begin{equation}
	s  (\vec x) = \sum_{i=1}^\infty \lambda_i^{ \alpha -\beta}  \left  ( \tilde{s}  (\vec x), \,  \zeta_i(\vec x) \right  ) \zeta_i (\vec x) \,	.
\end{equation}
In the above expansion, we assume that the eigenvalues are sorted in ascending order. Although it is possible to approximately consider a truncated spectral expansion, we do not pursue it further because of the difficulties in solving large eigenvalue problems. We use instead a rational series expansion to approximate the fractional operator as detailed in Section~\ref{sec:arbiraryExponent}. For later, note that the eigenvalues of the operator~$\set L$ lie within the interval~$[ \lambda_1, \, \infty )$ and the eigenvalues of its inverse~$\set L^{-1}$ within the interval~$[ 0, \, 1/ \lambda_1 ]$.

%
\subsection{Finite element approximation for integer exponents \label{sec:fractionalFE}}
%
In this section, we focus on the solution of the non-fractional partial differential equation~\eqref{eq:fracPDEdecA} with~\mbox{$\alpha \in \mathbb N_{\ge 1}$}. Its finite element discretisation for~\mbox{$\alpha =1$} is straightforward. However,~\mbox{$\alpha \ge 2$} requires in a standard finite element discretisation basis functions with higher-order smoothness. To avoid the use of smooth basis functions  we resort to a mixed variational formulation. To this end, the solution of the partial differential equation~\eqref{eq:fracPDEdecA} for a fixed~$\alpha \ge 2$, i.e.
\begin{equation}
	 \tilde {s}_{\alpha} (\vec x) \equiv  \tilde s (\vec x)  \,,
\end{equation}
can be obtained by recursively solving 
\begin{subequations} \label{eq:PDErecursive}
\begin{align} 
		\set  L \tilde{s}_{1} (\vec x) &= \frac{1}{\tau} g (\vec x) \, , \label{eq:PDErecursiveA} \\
		\set L \tilde s_{k} (\vec x) &= \tilde {s}_{k - 1} (\vec x) \, , \quad k = 2, \dotsc, \alpha \,. \label{eq:PDErecursiveB} 
\end{align}
\end{subequations}
Hence, the solution~$\tilde{s}(\vec x)$ is obtained by repeatedly solving the same partial differential equation with different right-hand sides.

To discretise the system of equations~\eqref{eq:PDErecursive}, it is sufficient to detail the discretisation of the first equation~\eqref{eq:PDErecursiveA} in the recursion. As mentioned its domain is all of $\mathbb R^d$  which we approximate with a sufficiently large domain~\mbox{$\widehat \Omega \subset \mathbb R^d$} with homogeneous Neumann boundary conditions. The  weak form reads: find~$\tilde{s}_1 (\vec x) \in \set H^1 (\widehat \Omega)$ such that 
\begin{equation}
	\int_{\widehat \Omega} \left ( \kappa^2  \tilde s_1 (\vec x ) v (\vec x)  +   \nabla \tilde s_1 (\vec x) \cdot  \nabla v  (\vec x) \right)  \D \vec x =  \frac{1}{\tau}  \int_{\widehat \Omega}  g(\vec x)v  (\vec x) \D \vec x \, \quad \forall \ v(\vec x)  \in \set H^1 (\widehat \Omega) \, ,
	\label{eq:spdeWeakForm}
\end{equation}
where~$ \nabla$ is the gradient operator, $v(\vec x)$ is a test function and~$\set H^1 (\widehat \Omega)$ is a standard Sobolev space.  To obtain the finite element approximation,~$\tilde s_1 (\vec x)$ and~$v(\vec x)$ are approximated by
\begin{equation} \label{eq:discIntExp}
	\tilde s_1 (\vec x) \approx \tilde s_1^h (\vec x)  = \sum_{i=1}^{n_u } \phi_i (\vec x) \tilde s_{1, i} \, ,  \quad   v  (\vec x) \approx v^h  (\vec x)  = \sum_{i=1}^{n_u } \phi_i (\vec x) v_i   \, , 
\end{equation}
where~$\phi_i (\vec x)$ are the finite element basis functions, and~$\tilde s_{1, i}$ and $v_i$ are the respective~$n_u$ nodal coefficients. The so discretised weak form yields the discrete system of equations
\begin{equation} \label{eq:PDErecursive1}
	\left ( \kappa^2 \vec M  + \vec K   \right )  \tilde {\vec s}_1 = \frac{1}{\tau} {\vec g} \, ,
\end{equation}
where~$\vec M$ is the mass matrix,~$\vec K$ is the discretised Laplacian matrix and~$\vec g$ is a random vector with the components
\begin{equation}
	g_i = \int_{\widehat \Omega} g(\vec x) \phi_i(\vec x)  \D \vec x  \, .
\end{equation}
After defining the stiffness matrix~$\vec L=  \kappa^2 \vec M  + \vec K  $, we can write  
\begin{equation} \label{eq:PDErecursive2}
	\vec L  {\tilde{\vec s}}_{1} = \frac{1}{\tau} { \vec g} \, .
\end{equation}
Considering that~$g (\vec x)$ is a Gaussian white noise, the components of the covariance matrix of~$\vec g$ are given by
\begin{align}
	\begin{split}
		\cov \left (  g_i , g_j \right ) &= \expect \left [  \int_{\widehat \Omega} \int_{\widehat \Omega} g(\vec x) \phi_i(\vec x) g(\vec x') \phi_j(\vec x') \D \vec x \D \vec x'  \right ] 
		=  \int_{\widehat \Omega} \int_{\widehat \Omega}  \phi_i(\vec x)  \expect \left [   g(\vec x)   g(\vec x') \right ] \phi_j(\vec x') \D \vec x \D \vec x'  \\
		&=  \int_{\widehat \Omega} \int_{\widehat \Omega}  \phi_i(\vec x)  \delta (\vec x - \vec x') \phi_j(\vec x') \D \vec x \D \vec x' =  \int_{\widehat \Omega} \phi_i(\vec x) \phi_j(\vec x) \D \vec x \, . 
	\end{split}	
\end{align}
Consequently,~$\vec g$ has the multivariate Gaussian density 
\begin{equation}
	p( \vec g ) = \set N (\vec 0, \vec M) \, .
\end{equation}
Hence, by recourse to~\eqref{eq:PDErecursive2} and linear transformation property of Gaussian densities, see~\ref{sec:multivariateNormalCov},   the solution vector~$\tilde{\vec{s}}_1$ has the density
\begin{equation}
	p ( \tilde {\vec s}_1 ) =  \set N \left (\vec 0, \, \vec C_{\tilde{s}, 1}  \right ) =  \set N \left (\vec 0, \,  \frac{1}{\tau^2} \vec L^{-1} \vec M \vec L^{-\trans} \right )   \, , 
\end{equation}
which can also be expressed as
\begin{equation}
	p ( \tilde {\vec s}_1 ) =  \set N \left (\vec 0, \, \vec Q_{\tilde{s}, 1}^{-1} \right )=  \set N \left (\vec 0, \,  \left (  \tau^2 \vec L^\trans \vec M^{-1} \vec L \right )^{-1}  \right )  \, .
\end{equation}
The covariance matrix~$\vec C_{\tilde{s}, 1}$ is dense. However, the precision matrix~$\vec Q_{\tilde{s}, 1}$ becomes sparse when~$\vec M$ is replaced by a diagonal lumped mass matrix. With a slight abuse of notation we denote in the following the lumped mass matrix~$\vec M$ with the same symbol like the consistent mass matrix. We can deduce after the finite element discretisation of~\eqref{eq:PDErecursiveB} that each of the solution vectors must have a zero-mean multivariate Gaussian density with the respective precision matrices
\begin{align}~\label{eq:PDEprecission} 
	\begin{split}
		\vec Q_{\tilde{s}, 1} & = \tau^2 \vec L^\trans \vec M^{-1} \vec L  \,, \\
		\vec Q_{\tilde{s}, {k}}&= \vec L^\trans \vec M^{-1} \vec Q_{\tilde{s},  k-1} \vec M^{-1} \vec L \,, \quad k = 2, \dotsc, \alpha \,, 
	\end{split}
\end{align}
so that 
\begin{equation} \label{eq:nonFracDensity}
	p ( \tilde{\vec s} )  = \set N \left ( \vec 0, \, \vec Q_{\tilde{s}}^{-1} \right ) \quad \text{with} \quad \vec Q_{\tilde{s}} \equiv \vec Q_{\tilde{s}, \alpha} \,.
\end{equation}
The precision matrix~$\vec Q_{\tilde s}$ is sparse when a lumped mass matrix is used throughout and becomes denser with increasing~$\alpha$.

%
%
\subsection{Generalisation to arbitrary exponents\label{sec:arbiraryExponent}}
\begin{figure}
	\centering
	\subfloat[][{Rational approximation}] {
		\includegraphics[height=0.315\textwidth]{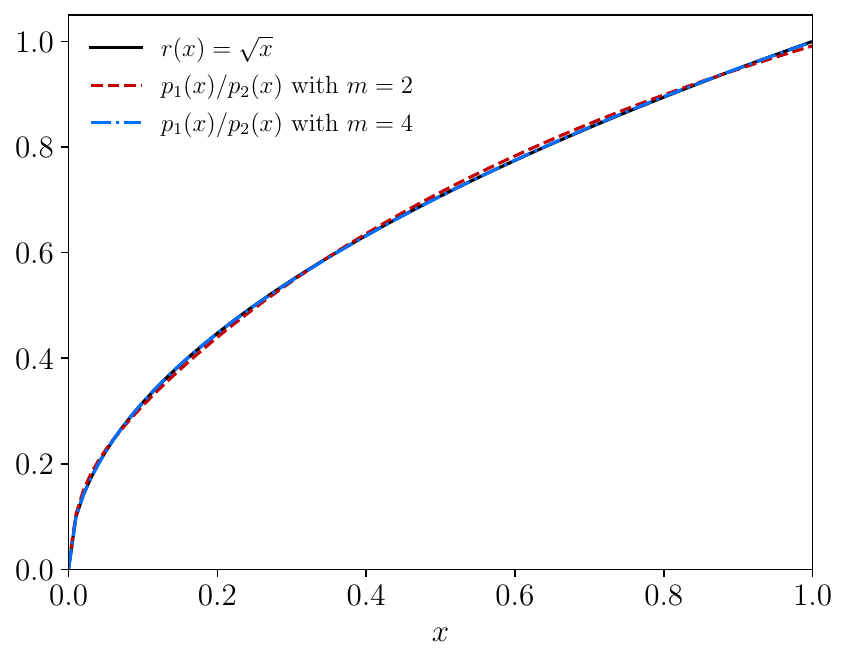} } 
	\hfil
	\subfloat[][{Approximation error for~$m = 2$}] {
		\includegraphics[height=0.315\textwidth]{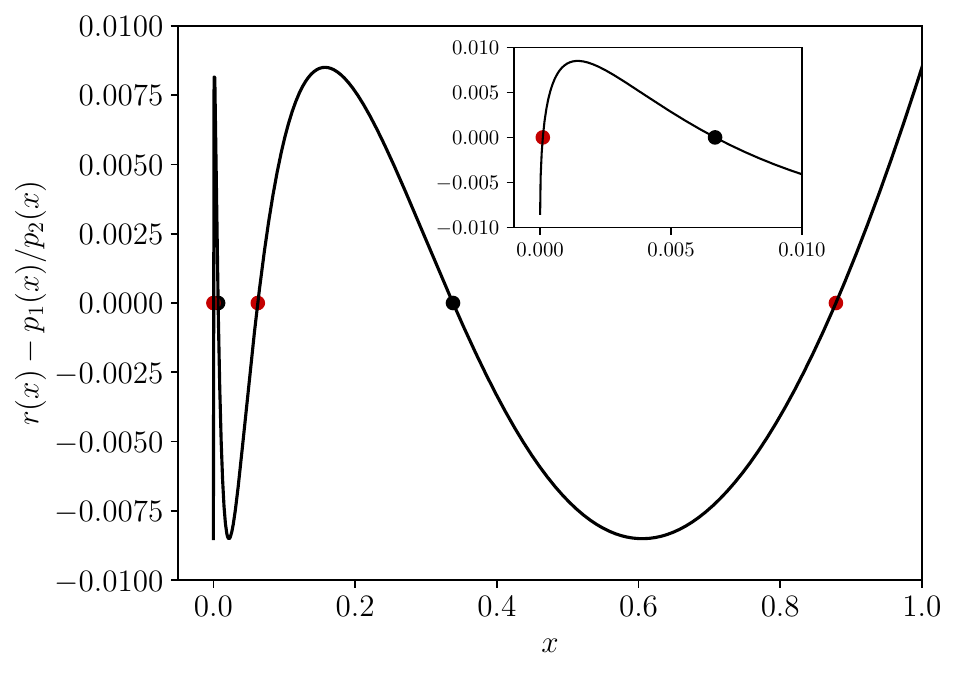} }
	\caption{Rational approximation of~$r(x) = \sqrt{x}$ using the BRASIL algorithm~\cite{hofreither2021algorithm}. As expected the accuracy of the approximation~\mbox{$r(x) \approx p_1(x) / p_2(x)$} improves with an increase of the polynomial degree~$m$ as shown in (a). The approximation error exhibits the equioscillation behaviour shown in~(b). The BRASIL algorithm interpolates~$r(x)$ at~$2m + 1$ points as indicated by the red and black dots in~(b). Specifically, the red dots represent the~$m+1$  chosen interpolation points $x_i$ and the black dots follow from the equioscillation property, see~\ref{sec:rationalInterp}.
	 \label{fig:rationalInterpolation}}
\end{figure}

%
We discuss next the solution of the fractional partial differential equation~\eqref{eq:fracPDEdecB}. To this end, we first consider the barycentric rational interpolation of a function
\begin{equation}
	r(x) = x^{\beta - \alpha} \quad \text{with} \quad \beta -\alpha \in  (d/4 -1, \, 1) \subset \mathbb R  \quad \text{and} \quad x \in (\varepsilon, \, 1) \subset \mathbb R \,,
\end{equation}
where $\varepsilon > 0$ is a small number. At the yet to be specified interpolation points~\mbox{$x_1, \, x_2, \dotsc, \, x_{m+1}$} the function values~$r(x_i)$ are assumed to be given. The rational interpolation of~$r(x)$ reads
\begin{equation} \label{eq:ratInterp}
	 r( x)  \approx  \frac{\sum_{i=1}^{m+1} a_i x^{i - 1} }{\sum_{j=1}^{m+1} b_j  x^{j-1}} =  \frac{a_{m+1} \prod_{i=1}^{m}  (x- c_i   ) }{b_{m+1} \prod_{j=1}^{m} ( x- d_j  ) } \, .
\end{equation}
We use the BRASIL (best rational approximation by successive interval length adjustment) library by Hofreither~\cite{hofreither2021algorithm} to obtain an optimal set of monomial coefficients~$a_i$ and~$b_j$ that minimises the interpolation error. Figure~\ref{fig:rationalInterpolation} illustrates the interpolation of~\mbox{$r(x) = \sqrt{x}$} using the BRASIL algorithm. See also~\ref{sec:rationalInterp} for a brief discussion on rational interpolation.  The roots of the numerator and denominator~$c_i$ and~$d_j$ are real and distinct and are computed after determining the coefficients~$a_i$ and~$b_j$. 

Formally, replacing~$x$ by the inverse operator~$\set L^{-1}$ in the rational interpolation of~$r(x)$ yields an approximation for the inverse fractional operator~$\set L^{\alpha - \beta}$. As pointed out by Bolin and Kirchner~\cite{bolin2020rational}, approximating the inverse operator~$\set L^{-1}$ rather than the operator~$\set L$ itself gives a more efficient algorithm. Recall that~$\lambda_1$ denotes the smallest eigenvalue of~$\set{L}$ so that the eigenvalues of~$\set L^{-1}$ lie within the bounded interval~$[0, \, 1/\lambda_1]$ while the eigenvalues of~$\set L$ lie within the half-bounded interval~$[\lambda_1, \, \infty)$, see Section~\ref{sec:fractionalPDE}. The domain of the variable~$x$ in~\eqref{eq:ratInterp} must be chosen in dependence of the interval for the eigenvalues. A bounded interval is evidently a more convenient choice of a domain. Furthermore, a simple rescaling of the variables can lead to an inverse operator with eigenvalues within the interval~$[0, \, 1]$. 

Prior to introducing the inverse operator we re-express \eqref{eq:ratInterp} as  
\begin{equation} \label{eq:ratInterpInv}
	 r \left ( x^{-1} \right )  \approx   \frac{a_{m+1} \prod_{i=1}^{m}  \left (  x^{-1} - c_i \right ) }{b_{m+1} \prod_{j=1}^{m} \left ( x^{-1} - d_j \right ) } =  \frac{a_{m+1} \prod_{i=1}^{m}  \left (1- c_i x \right ) }{b_{m+1} \prod_{j=1}^{m} \left ( 1- d_j  x   \right ) }  \, , 
\end{equation}
and obtain for the inverse operator
\begin{equation}
	\set L^{\alpha - \beta} \approx a_{m+1} \prod_{i=1}^{m}  \left (\set I - c_i \set L  \right ) \left ( b_{m+1}  \prod_{j=1}^{m}  \left (\set I - d_j \set L   \right )  \right )^{-1}  \, , 
\end{equation}
or, more compactly,
\begin{equation}
	\set L^{\alpha - \beta} \approx \set F_r(\set L) \, \set F_l^{-1} (\set L) \, . 
\end{equation}
Introducing this approximation into the fractional partial differential equation~\eqref{eq:fracPDEdecB} and noting that~\mbox{$\set F_r(\set L)$} and~\mbox{$\set F_l(\set L)$}  commute yields
\begin{equation} \label{eq:fractionalCont}
	\set F_l  (\set L)  \,   s (\vec x) = \set F_r(\set L) \,   \tilde{s} (\vec x)  \, .
\end{equation}

The operators~$\set F_l  (\set L)$ and~$\set F_r(\set L) $ are discretised similar to the discretisation of the non-fractional partial differential equation~\eqref{eq:fracPDEdecA} in the preceding section, see also Bolin and Kirchner~\cite[Appendix A]{bolin2020rational}. It is easy to show that the discretisation yields the two matrices
\begin{equation}
	\vec F_l =   b_{m+1}  \prod_{j=1}^{m}  \left (\vec I - d_j \vec M^{-1} \vec L   \right )  \, , \quad \vec F_r = a_{m+1} \prod_{i=1}^{m}  \left (\vec I - c_i \vec M^{-1}\vec L  \right )  \, .
\end{equation}
After discretisation,~\eqref{eq:fractionalCont} takes the form
\begin{equation}
	\vec F_l \vec s = \vec F_r \tilde {\vec s} \, .
\end{equation}
The right-hand side vector~$\tilde {\vec s}$ is the solution of the integer component and has, as derived in Section~\ref{sec:fractionalFE}, the multivariate Gaussian density~\eqref{eq:nonFracDensity}.  Therefore, we can write for the density of the overall solution
\begin{equation}
	p (\vec s) = \set N \left ( \vec 0, \,  \vec Q_{s}^{-1}  \right )   = \set N \left ( \vec 0, \,  \left (  \vec F_r ^{-\trans}  \vec F_l^{\trans} \vec Q_{\tilde{s}}  \vec F_l  \vec F_r^{-1} \right )^{-1} \right )  \, , 
	\label{eq:overallSPDESolution}
\end{equation}
where we used the fact that~$ \vec F_r $ and~$ \vec F_l$ commute. The precision matrix~$ \vec Q_{s} $ is not sparse because~$\vec F_r^{-1}$ is a dense matrix.  However, the expression~$  \vec F_l^{ \trans}  \vec Q_{\tilde{s}}  \vec F_l$ is sparse. These observations are critical for Gaussian process regression and statistical finite element analysis for large-scale problems as detailed in Sections~\ref{sec:gpRegression} and~\ref{sec:statFEM}. 
%
\subsection{Generalisation to non-standard random fields and manifolds \label{sec:maternExtensions}}
%
The Mat\'ern random fields introduced so far are stationary, isotropic and restricted to Euclidean domains.  These restrictions can be relaxed by modifying the fractional stochastic PDE representation of Mat\'ern fields slightly. To illustrate this, we consider a two-manifold, i.e.\ surface,~$\set{M} \subset \mathbb R^3$ described by the map
\begin{equation}
	\vec \varphi(\vec \xi) : \,   \vec \xi = (\xi_1, \, \xi_2)^\trans \mapsto \vec x = (x_1, \, x_2, \, x_3 )^\trans \, ,
\end{equation}
where~\mbox{$\vec \xi  \in \mathbb R^2$} are the parametric coordinates of the surface points~$\vec x \in \mathbb R^3$. Without loss of generality, the surface is assumed to be closed to sidestep the discussion of boundaries. The covariant and contravariant surface basis vectors~$\vec a_i $ and~$\vec a^i$ are defined as 
\begin{equation}
	 \vec a_i = \frac{\partial \vec \varphi(\vec \xi) }{\partial \xi_i} \, , \quad   \vec a_i \cdot \vec a^j = \delta_i^j \, ;  
\end{equation}
see, e.g.~\cite{Ciarlet:2005aa}. The respective covariant and contravariant metric tensors~$\vec G$ and $\vec G^{-1}$ satisfy the relation~$\vec G \vec G^{-1} = \vec I$ and are given by
\begin{equation}
	 \vec G = \begin{pmatrix}   \vec a_1  \cdot \vec a_1 & \vec a_1 \cdot \vec a_2  \\[1.em]
	 \vec a_2 \cdot \vec a_1  & \vec a_2 \cdot \vec a_2  \end{pmatrix} \, , \quad   	 
	 \vec G^{-1} = \begin{pmatrix}   \vec a^1  \cdot \vec a^1 & \vec a^1 \cdot \vec a^2  \\[1.em]
	 \vec a^2 \cdot \vec a^1  & \vec a^2 \cdot \vec a^2  \end{pmatrix}  \, .
\end{equation}

\begin{figure}
	\centering
	\includegraphics[width=0.25\textwidth]{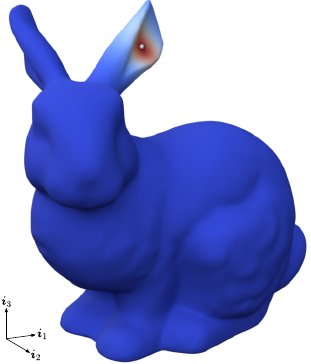}
	\hfil
	\includegraphics[width=0.25\textwidth]{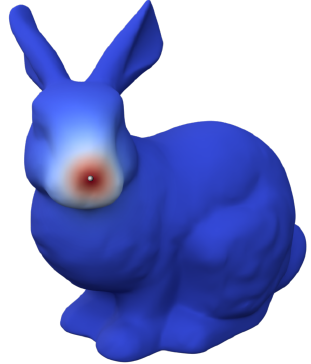}
	\hfil
	\includegraphics[width=0.25\textwidth]{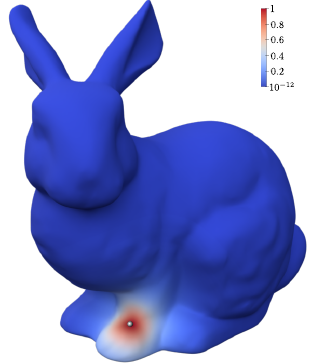}
	\caption{Three isocontours of an isotropic and stationary covariance function. In each of the plots the covariance function is centred at the point indicated by the white sphere. A sample from the respective Gaussian random field is depicted in Figure~\ref{fig:bunnyStationarySample}. The parameters of the stochastic partial differential equation are~$\kappa= \sqrt{32}/5$, $\beta=1$, and $\tau=5/\sqrt{128\pi}$ and correspond  to the Mat\'ern parameters~$\sigma = 1$, $\ell = 1.25$, and $\nu = 1$.    \label{fig:bunnyStationary}}
%
\vspace{1.5em}
%
	\centering
	\includegraphics[width=0.25\textwidth]{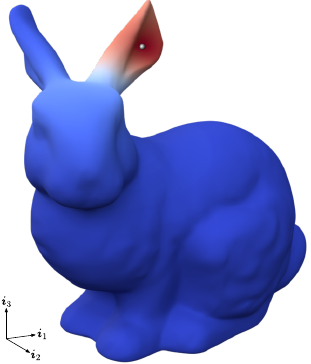}
	\hfil
	\includegraphics[width=0.25\textwidth]{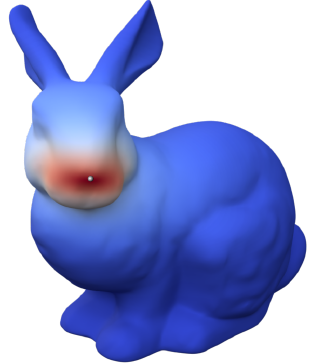}
	\hfil
	\includegraphics[width=0.275\textwidth]{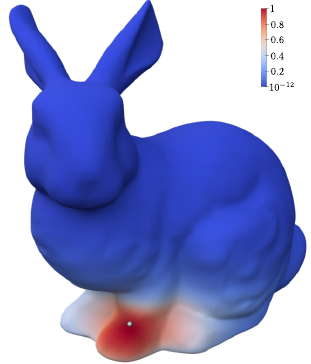}
	\caption{Three isocontours of an anisotropic and stationary covariance function. In each of the plots the covariance function is centred at the point indicated by the white sphere. A sample from the corresponding Gaussian random field is depicted in Figure~\ref{fig:bunnyAnisotropicSample}.  To introduce anisotropy the diffusion matrix is chosen as~\mbox{$\vec H = \diag (10, \,10, \, 1) $}. The other two parameters~$\sigma=1$ and~$\nu=1$ are the same as for the covariance function in Figure~\ref{fig:bunnyStationary}. Notice that the covariance function is  more spread along the~$\vec i_1$--$\vec i_2$ plane. \label{fig:bunnyAnisotropic}}
%
\vspace{1.5em}
%
	\centering
	\includegraphics[width=0.25\textwidth]{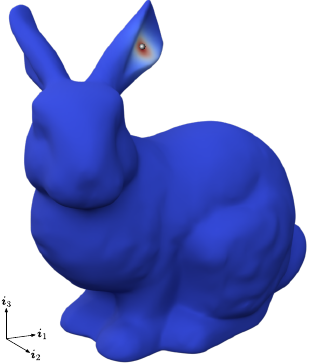}
	\hfil
	\includegraphics[width=0.25\textwidth]{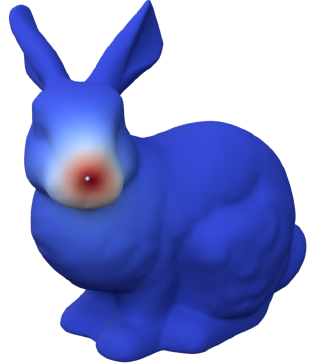}
	\hfil
	\includegraphics[width=0.25\textwidth]{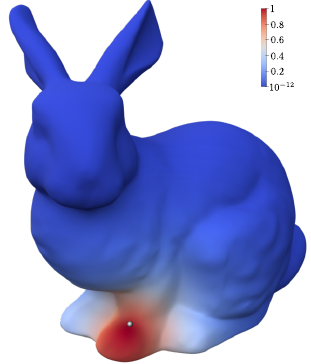}
	\caption{
	Three isocontours of a non-stationary and isotropic covariance function. In each of the plots the covariance function is centred at the point indicated by the white sphere. A sample from the respective Gaussian random field is depicted in Figure~\ref{fig:bunnyNonStationarySample}. The  Mat\'ern length-scale~$\ell(x_3) = 0.2(20 - x_3)$ is chosen to vary linearly along the height of the bunny, i.e.\ the length-scale increases linearly from $0.25$ (top) to $3.3$ (bottom). The other two parameters~$\sigma=1$ and~$\nu=1$ are the same as for the covariance function in Figure~\ref{fig:bunnyStationary}. Notice that  the convariance function is more spread out at the feet than at the head of the bunny.
	 \label{fig:bunnyNonStationary}}
\end{figure}

We focus, first, on non-fractional, isotropic, stationary random fields and choose the surface~$\set M$ as the domain of the stochastic PDE~\eqref{eq:fracPDE}. The  weak form of the respective first equation in the recursion~\eqref{eq:PDErecursive} is given by  
\begin{equation} \label{eq:maniWeak1}
	\int_{\set M} \left ( \kappa^2  \tilde{s}_1 (\vec \xi ) v (\vec \xi)  +    \nabla   \tilde{s}_1 (\vec \xi)  \cdot  \nabla v  (\vec \xi) \right)  \D \set M =   \int_{\set M} \frac{1}{\tau}  g(\vec \xi)v  (\vec \xi) \D \set M   \, ,  
\end{equation}
where the surface gradient, e.g.\ of~$\tilde{s}_1$ is defined as
\begin{equation}
	\nabla  \tilde{s}_1 = \sum_{i=1}^2 \frac{\partial \tilde{s}_1 (\vec \xi)}{\partial \xi_i} \vec a^i \, .
\end{equation}
Hence, the weak form of~\eqref{eq:maniWeak1} can be rewritten as   
\begin{equation} \label{eq:maniWeak2}
	\int_{\set M} \left ( \kappa^2  \tilde{s}_1 (\vec \xi ) v (\vec \xi)  +   \frac{  \partial \tilde {s}_1 (\vec \xi ) }{\partial \vec \xi } \cdot \vec G^{-1}  \frac{  \partial v ( \vec \xi) }{\partial \vec \xi }     \right) \sqrt{\det \vec G } \D \vec \xi =   \int_{\set M} \frac{1}{\tau}  g(\vec \xi)v  (\vec \xi) \sqrt{\det \vec G } \D \vec \xi    \, .
\end{equation}
Following the isogeometric approach, the geometry~$\vec x(\vec \xi)$, the solution field~$ \tilde{s}_1 (\vec \xi)$ and the test field~$ v(\vec \xi)$ are approximated with same basis functions
\begin{equation}
	 \vec x (\vec \xi) \approx \vec x^h (\vec \xi) = \sum_{i=1}^{n_u } \phi_i (\vec \xi) \vec  x_ i  \, ,  \quad  \tilde s_1 (\vec \xi) \approx  \tilde s_1^h (\vec \xi)  = \sum_{i=1}^{n_u } \phi_i (\vec \xi) \tilde s_{1, i} \, ,  \quad   v  (\vec \xi) \approx  v^h  (\vec \xi)  = \sum_{i=1}^{n_u } \phi_i (\vec \xi) v_i   \, , 
\end{equation}
where~$ \vec  x_ i$, $\tilde s_{1, i}$ and $v_i$ are the respective nodal values. With a slight abuse of notation, we denote the basis functions~$\phi_i(\vec \xi)$ on the parametric domain with the same symbol like the basis functions~$\phi_i(\vec x)$ on the actual domain defined in~\eqref{eq:discIntExp}.  See Figure~\ref{fig:bunnyStationarySample} for an isotropic and stationary random field on a manifold surface and Figure~\ref{fig:bunnyStationary} for the respective covariance function centred at three selected locations on the surface. 

Next, we consider an anisotropic random field by introducing a diffusion matrix~$\vec H$ given in terms of the global orthonormal basis vectors~$\left( \vec i_1, \, \vec i_2, \, \vec i_3\right)$. To transform~\mbox{$\vec H \in \mathbb R^{3 \times 3}$} into the contravariant surface basis we introduce the transformation matrix 
\begin{equation}
	\vec T = \begin{pmatrix} \vec i_1 \cdot \vec a^1 & \vec i_1 \cdot \vec a^2 \\ \vec i_2 \cdot \vec a^1 & \vec i_2 \cdot \vec a^2 \\ \vec i_3 \cdot \vec a^1 & \vec i_3 \cdot \vec a^2 \end{pmatrix} \,.
\end{equation}
To consider anisotropy it is sufficient to replace in~\eqref{eq:maniWeak2} the contravariant metric~$\vec G^{-1}$ by
\begin{equation}
	\widetilde {\vec H} = \vec T^\trans \vec H \vec T \, .
\end{equation}
As expected~$\widetilde{\vec H} = \vec G^{-1}$ when~$\vec H = \vec I$. In passing, we note that non-stationary random fields can be modelled by choosing the SPDE parameters as 
\begin{equation}\label{eq:fracPDEvarsGeneral}
	\kappa(\vec x) = \frac{ \sqrt{2 \nu} }{ \ell(\vec x) } \, ,  \quad \tau^2(\vec x) = \frac{\Gamma (\nu)}{ \sigma^2 \Gamma (\nu+d/2) (4 \pi )^{d/2} \left(\kappa(\vec x)\right)^{2\nu} \sqrt{\det{\vec H(\vec x)}}}   \, , 
\end{equation} 
where the length-scale~$\ell (\vec x)$ and the diffusion matrix~$\vec H(\vec x)$ depend on the location on the surface. 
The subsequent steps in solving the fractional SPDE with an arbitrary exponent follow Sections~\ref{sec:fractionalFE} and~\ref{sec:arbiraryExponent}.

 Figures~\ref{fig:bunnyStationarySample}, \ref{fig:bunnyAnisotropicSample} and~\ref{fig:bunnyNonStationarySample} depict the isocontours of different kinds of random fields determined with the proposed approach. In Figures~\ref{fig:bunnyStationary}, \ref{fig:bunnyAnisotropic}, and~\ref{fig:bunnyNonStationary} the corresponding covariance functions centred at three different locations on the surface are shown. See~\ref{sec:multivariateNormalSampling} for sampling from Gaussian densities. 

%
\section{Gaussian process regression \label{sec:gpRegression}}
%
\subsection{Statistical observation model}
%
We are given the observation set~$\{ (\vec x_i, \,  y_i ) \}_{i=1}^{n_y}$ with~$n_y$ observations~$ y_i \in \mathbb R$ at the locations~$\vec x_i \in \mathbb R^d$. The observation set is assumed  to correspond to a random function~$s(\vec x)$ sampled from an unknown Gaussian process. We posit for regression a Mat\'ern field, in other words a process, which represents according to the Bayesian viewpoint the prior.  As discussed, the finite element discretisation of the fractional PDE representation of the Mat\'ern field leads to the multivariate Gaussian density~\eqref{eq:overallSPDESolution}, restated here for convenience,   
\begin{equation}  	
	p (\vec s) =  \set N \left ( \vec 0, \,  \vec C_{s}  \right ) =  \set N \left ( \vec 0, \,  \vec Q_{s}^{-1}  \right )   = \set N \left ( \vec 0, \,  \left (  \vec F_r ^{-\trans}  \vec F_l^{\trans} \vec Q_{\tilde{s}}  \vec F_l  \vec F_r^{-1} \right )^{-1} \right )  \, ,  
\end{equation}
where~$\vec s \in \mathbb R^{n_u}$ is the vector of the field values at the finite element nodes, and \mbox{$\vec C_s, \, \vec Q_s \in  \mathbb R^{n_u \times n_u}$} are the covariance and precision matrices as detailed before. The covariance and precision matrices depend on the hyperparameters~\mbox{$\{ \sigma, \, \ell, \, \nu \}$} of the Mat\'ern covariance function~\eqref{eq:maternKernel} or the respective parameters~\mbox{$\{\tau, \, \kappa, \, \beta \}$} of the SPDE~\eqref{eq:fracPDEvars}. Anisotropic and non-stationary random fields introduced in Section~\ref{sec:arbiraryExponent} may have additional hyperparameters. In Gaussian process regression the hyperparameters are first assumed to be fixed and later learned from the observation data. The related approaches are in statistics and machine learning known as empirical Bayes, evidence approximation or type II maximum likelihood. 

The standard statistical observation model for Gaussian process regression takes the form
\begin{equation} \label{eq:statMgp}
	\vec y  = \vec P \vec s + \vec e  \, , 
\end{equation}
where $\vec y \in \mathbb R^{n_y}$ is the vector containing all observations, $\vec P \in \mathbb R^{n_ y \times n_u}$  is the observation matrix and~$\vec e \in \mathbb R^{n_y}$  is an additive observation noise vector. The entries of the observation matrix depend on the type and location of the observations. For instance, if~$y_i$ corresponds to a nodal value~$s_j$, the observation matrix has in the $i$-th row only the one non-zero entry~$P_{ij} =1$. For other observation locations and types it is straightforward to determine the respective observation matrix coefficients using the basis functions of the finite element discretisation. Furthermore, the vector of the field values~$\vec s$ and the noise vector~$\vec e$ are assumed to be statistically independent and the probability density of the noise vector is given by 
\begin{equation}
	p (\vec e) = \set N \left(\vec 0 , \, \sigma_e^2 \vec I\right) \, . 
\end{equation}
Hence, the noise vector components are independent and identically distributed and have the standard deviation $\sigma_e$. 

The probability density of the nodal field value vector~$\vec s$ for a given observation vector~$\vec y \in \mathbb R^{n_y}$  is according to the Bayes formula given by
\begin{equation} \label{eq:bayesGP}
	p (\vec s \vert \vec y) = \frac{p(\vec y \vert \vec  s) p(\vec s)}{p(\vec y)}    \, ,
\end{equation}
where~$p(\vec s \vert \vec  y) $ is the posterior density,  $p(\vec y \vert \vec  s)$ the likelihood, $p(\vec s)$ the prior and $p(\vec y)$ the marginal likelihood. For the given observation model~\eqref{eq:statMgp} the likelihood takes the form 
 \begin{equation}
 	p(\vec y \vert \vec s) = \set N \left ( \vec P \vec s, \, \sigma_e^2 \vec I  \right ) \, .  
 \end{equation}
The marginal likelihood, or the evidence,~$p(\vec y)$, which plays a key role in determining the hyperparameters of the prior~$p(\vec s)$, i.e.~\mbox{$\{ \sigma, \, \ell, \, \nu \}$} in~\eqref{eq:maternKernel}, is  obtained by marginalising out~$\vec s$, i.e.
 \begin{equation}
 	p(\vec y) = \int p(\vec y  \vert  \vec s) p(\vec s) \D \vec s \, .
 \end{equation}
In the presented examples, we determine the hyperparameters by maximising the marginal likelihood~$p(\vec y)$ for an observation vector~$\vec y$. Indeed, it is numerically more convenient to maximise the log marginal likelihood~$\log p (\vec y)$. Given that the~$\log$ is a monotonically increasing function, the maximisation of~$p(\vec y)$ and~$\log p (\vec y)$ are equivalent.  

The Gaussian densities appearing in the Bayes formula~\eqref{eq:bayesGP}  can be either expressed via their covariance or precision matrices.  In machine learning usually the covariance formulation is preferred.

%
\subsection{Standard covariance formulation  \label{sec:gpReview}}
%
When all the probability densities appearing in the Bayes formula~\eqref{eq:bayesGP} are Gaussians, as is the case in the presented approach, it is convenient to determine the posterior~$p(\vec s \vert \vec y)$ directly from the joint density 
\begin{equation} \label{eq:jointDistribCov}
	 p(\vec s, \, \vec y) = \set{N} \left(
	\begin{pmatrix}
		\vec 0 \\ \vec 0 
	\end{pmatrix}  ,
	\begin{pmatrix}
		\vec C_s  & \vec C_s \vec P^\trans 
		\\[0.15em]
		\vec P \vec C_s &  \vec P \vec C_s  \vec P^\trans + \sigma_e^2 \vec I
	\end{pmatrix}
	\right) \,.
\end{equation}
See~\ref{sec:multivariateNormal} for the relevant properties of Gaussians used in deriving this expression.  According to the conditioning property~\eqref{eq:gaussianDensityConditional}, the posterior is given by
\begin{equation} 
	p (\vec s \vert \vec  y )  = \set{N} \left(  \overline{\vec s}_{\vert y}, \,  \vec C_{s \vert y}  \right)  \,,
\end{equation}
where 
\begin{subequations}
	\begin{align}
		\overline{\vec s}_{\vert y} &= \vec C_s \vec P^\trans  \left(  \vec P \vec C_s  \vec P^\trans + \sigma_e^2 \vec I \right)^{-1} \vec y  \,,
		\label{eq:gpQConditionalMeanC}
		\\
		\vec C_{s\vert y} & = \vec C_s - \vec C_s \vec P^\trans  \left(  \vec P \vec C_s  \vec P^\trans + \sigma_e^2 \vec I \right)^{-1} \vec P \vec C_s   \, .
		\label{eq:gpQConditionalCovarianceC} 
	\end{align}
\end{subequations}
The expression to be inverted within the brackets is a dense matrix of size~$n_y \times n_y$. Hence, the number of observations~$n_y$ has to be relatively small for the posterior mean~$\overline{\vec s}_{\vert y} $ and covariance~$\vec C_{s\vert y} $ to be numerically computable.

Furthermore, according to the marginalisation property of Gaussians, see~\ref{sec:multivariateNormalConditional}, the marginal likelihood can be read from~\eqref{eq:jointDistribCov} as 
\begin{equation}  
	p(\vec y) =   \set N  \left( \vec 0, \, \vec P \vec C_s  \vec P^\trans + \sigma_e^2 \vec I \right) \, , 
\end{equation}
and its logarithm is given by 
\begin{equation} 
	 \log p(\vec y)  =  - \frac{1}{2} \vec y^\trans  \left  (\vec P \vec C_s  \vec P^\trans + \sigma_e^2 \vec I \right )^{-1} \vec y  - \frac{1}{2} \log \det \left (\vec P \vec C_s  \vec P^\trans  + \sigma_e^2 \vec I \right )  - \frac{n_y}{2}  \log 2\pi     \, .
\end{equation}
%
\subsection{Sparse precision formulation \label{sec:gpPrecision}}
%
The limitation of Gaussian process regression to a relatively small number of observations can be ameliorated by expressing the probability densities in the Bayes formula~\eqref{eq:bayesGP} using their precision matrices.   However, as noted in Section~\ref{sec:arbiraryExponent}, the precision matrix~$\vec Q_s$ of the prior~$p(\vec s)$ is not always sparse. It is dense when the exponent~$\beta$ in the stochastic PDE~\eqref{eq:fracPDE} is fractional. To arrive at a posterior~$p(\vec s \vert \vec y)$  and a marginal likelihood~$p(\vec y)$ with sparse precision matrices, we introduce the auxiliary (latent) variable~$\vec t \in \mathbb R^{n_u}$ such that 
\begin{equation} \label{eq:gpAuxiliaryVariable}
	\vec t = \vec F_r^{-1} \vec s  \, ,
\end{equation}
and infer first~$\vec t$ and then~$\vec s$  from the given observations~$\vec y$. We emphasise that~$\vec F_l \equiv \vec F_r \equiv \vec I$ when~$\beta \in \mathbb N$ resulting in an inherently sparse precision matrix. The probability density of the auxiliary variable~$\vec t$ is in view of~\eqref{eq:overallSPDESolution} and the linear transformation property of Gaussian vectors, see~\ref{sec:multivariateNormalCov}, given by 
\begin{equation}
	p(\vec t) = \set N \left(\vec 0, \, \vec Q_t^{-1}\right) =  \set N \left ( \vec 0, \,  \left ( \vec F_r^\trans \vec Q_s \vec F_r \right )^{-1} \right ) =   \set N \left ( \vec 0, \,  \left (   \vec F_l^{\trans} \vec Q_{\tilde{s}}  \vec F_l   \right )^{-1} \right )  \, .
\end{equation}
The precision matrix~$\vec Q_t $ is indeed sparse considering that~$\vec F_l$ and~$ \vec Q_{\tilde{s}}$ are sparse.

By introducing the auxiliary variable~$\vec t$ into the statistical observation model~\eqref{eq:statMgp}, we obtain the modified model
 \begin{equation}
	\vec y  = \vec P \vec F_r \vec t + \vec e  \, . 
\end{equation}
The respective joint density of the auxiliary and observation vectors is given by
\begin{equation}
	\begin{split}
		p (\vec t, \, \vec y) &=  \set{N} \left(
		\begin{pmatrix}
			\vec 0 \\  \vec 0
		\end{pmatrix}  ,
		\begin{pmatrix}
			 \vec Q_t^{-1}  &  \vec Q_t^{-1} \vec F_r^\trans \vec P^\trans
			\\[0.15em]
			\vec P   \vec F_r   \vec Q_t^{-1}&  \vec P   \vec F_r  \vec Q_t^{-1}  \vec F_r^\trans \vec P^\trans + \sigma_e^2 \vec I 
		\end{pmatrix}
		\right)  =
		\set{N} \left(
		\begin{pmatrix}
			\vec 0 \\  \vec 0
		\end{pmatrix}  ,
		\frac{1}{\sigma_e^2}
		\begin{pmatrix}
			 \sigma_e^2  \vec Q_t +  \vec F_r ^\trans \vec P^\trans \vec P \vec F_r  & - \vec F_r ^\trans \vec P^\trans
			\\
			-  \vec P \vec F_r   & \vec I
		\end{pmatrix}^{-1} \right) \,.
	\end{split}
	\label{eq:auxiliaryJointDensity}
\end{equation}
The last~$2\times2$ block matrix is determined by straightforward inversion of the first $2\times2$ block matrix. Thus, the posterior density of the auxiliary vector~$\vec t$ conditioned on the observations~$\vec y$ reads
\begin{equation}
	p(\vec t \vert \vec y ) = \set{N} \left( \overline{\vec t}_{\vert y},  \, \vec Q_{t\vert y}^{-1} \right)  \,,
\end{equation}
where
\begin{subequations}
	\begin{align}
		\overline{\vec t}_{\vert y} &=   \frac{1}{\sigma_e^2}  \vec Q_{t\vert y} ^{-1}   \vec F_r^\trans \vec P^\trans \vec y \, ,
		\label{eq:gpConditionalMeanT}
		\\
		\vec Q_{t\vert y} & =   \vec Q_t +  \frac{1}{\sigma_e^2} \vec F_r^\trans \vec P^\trans \vec P \vec F_r  \, ;
		\label{eq:gpConditionalCovarianceT}
	\end{align}
	\label{eq:gpConditionalMeanAndCovarianceT}
\end{subequations}
cf.~\eqref{eq:gaussianDensityConditionalQ}. The sought density of the nodal values~$\vec s$ is given by
\begin{equation}  \label{eq:sparsePrecPosterior}
	p(\vec s \vert \vec y) =  \set{N} \left(  \overline{\vec s}_{\vert y},  \, \vec Q_{s \vert y}^{-1} \right)  =   \set{N} \left( \vec F_r  \overline{\vec t}_{\vert y},  \, \left ( \vec F_r ^{-\trans} \vec Q_{t\vert y}  \vec F_r^{-1} \right )^{-1} \right)  \,. 
\end{equation}
As is clear from~\eqref{eq:gpConditionalCovarianceT} the precision matrix~$\vec Q_{t\vert y}$ is sparse so that it is sufficient to factorise it only once. Subsequently, the mean~$\overline{\vec s}_{\vert y}$ is determined by first solving for~$\overline{\vec t}_{\vert  y}$ using~$ \vec F_r^\trans \vec P^\trans \vec y$ as the right-hand side and then multiplying the result by~$\vec F_r$. The covariance matrix~$\vec C_{s \vert y} = \vec Q_{s \vert y}^{-1}  $ is dense so that it is impossible to determine all its components for large~$n_u$. However, we can determine selected components  $(i, \, j)$ by noting that 
\begin{equation}
	\vec i_i  \cdot  \vec C_{s \vert y}  \vec i_j  = \sigma_e^2  \vec i_i  \cdot \vec F_r  \vec Q_{t\vert y}^{-1}  \vec F_r^\trans  \vec i_j   \, , 
\end{equation}
where~$\vec i_i$ and~$\vec i_j$ are two Cartesian basis vectors, and making use of the already factorised~$\vec Q_{t\vert y}$.

According to~\eqref{eq:auxiliaryJointDensity} and the marginalisation property of Gaussians the marginal likelihood is given by
\begin{equation}
	p(\vec y)  = \set{N}\left(  \vec 0, \,   \vec P   \vec F_r  \vec Q_t^{-1}  \vec F_r^\trans \vec P^\trans + \sigma_e^2 \vec I   \right)   \, . 
\end{equation}
Taking its logarithm yields
\begin{equation}
	 \log p(\vec y) = -  \frac{1}{2} \vec y^\trans \left(  \vec P   \vec F_r  \vec Q_t^{-1}  \vec F_r^\trans \vec P^\trans + \sigma_e^2 \vec I   \right)^{-1}  \vec y  -  \frac{1}{2} \log \det  \left(  \vec P   \vec F_r  \vec Q_t^{-1}  \vec F_r^\trans \vec P^\trans + \sigma_e^2 \vec I  \right) - \frac{n_y}{2} \log 2 \pi  \, .
\end{equation}
The first two terms on the right-hand side are costly to evaluate given that the matrix in the bracket is dense. Using the Sherman--Morrison--Woodbury formula and the matrix determinant lemma~\cite{sherman1950adjustment,woodbury1950inverting}, the two terms can be simplified as  
\begin{subequations}
\begin{align}
	\vec y^\trans \left(    \vec P   \vec F_r  \vec Q_t^{-1}  \vec F_r^\trans \vec P^\trans + \sigma_e^2 \vec I    \right)^{-1}  \vec y&= \frac{1}{ \sigma_e^2} \vec y^\trans \vec y -\frac{1}{ \sigma_e^2}  \vec y^\trans \vec P \vec F_r \left( \sigma_e^2 \vec Q_t +  \vec F_r^\trans \vec P^\trans   \vec P \vec F_r  \right)^{-1} \vec F_r^\trans \vec P^\trans \vec y  \, ,
\\ 
	\det \left( \vec P   \vec F_r  \vec Q_t^{-1}  \vec F_r^\trans \vec P^\trans + \sigma_e^2 \vec I  \right)  &= \left(\sigma_e^2 \right)^{n_y - n_u} \frac{ \det\left( \sigma_e^2 \vec Q_t +  \vec F_r^\trans \vec P^\trans   \vec P \vec F_r  \right) }{\det \left( \vec Q_t \right) }  \, .
\end{align}
\end{subequations}
Both expressions can be efficiently evaluated for large~$n_y$ using only sparse matrix operations. As an implementation note, it is numerically more robust to determine the log determinant of large matrices directly from their sparse factorisation.

\begin{figure}
	\centering
	\hfil
	\subfloat[][Samples drawn from SPDE prior \label{fig:gpDemoPrior}] {
		\includegraphics[width=0.4 \textwidth]{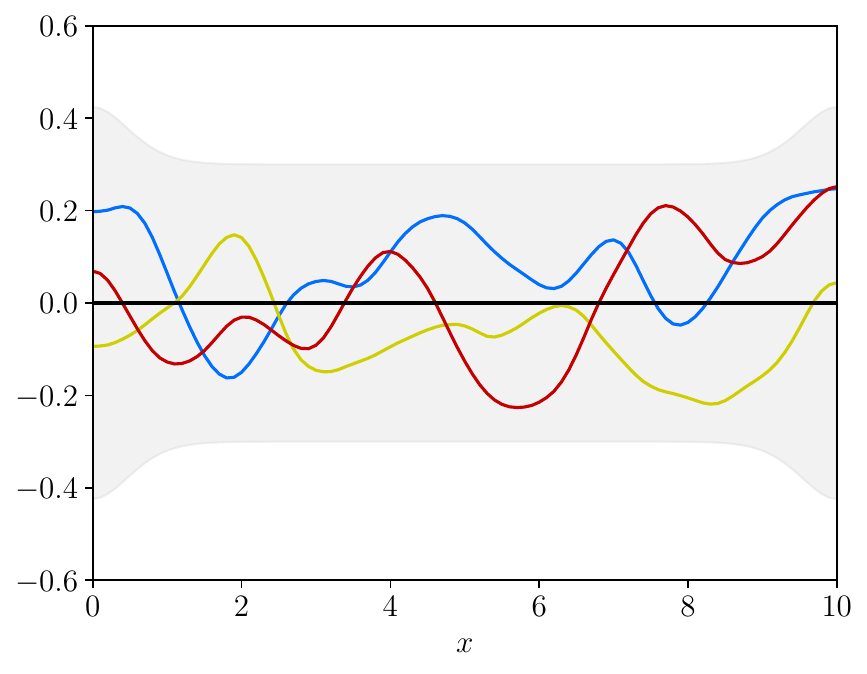} }
	\hfil
	\subfloat[][Posterior \label{fig:gpDemoPosterior}] {
		\includegraphics[width=0.4 \textwidth]{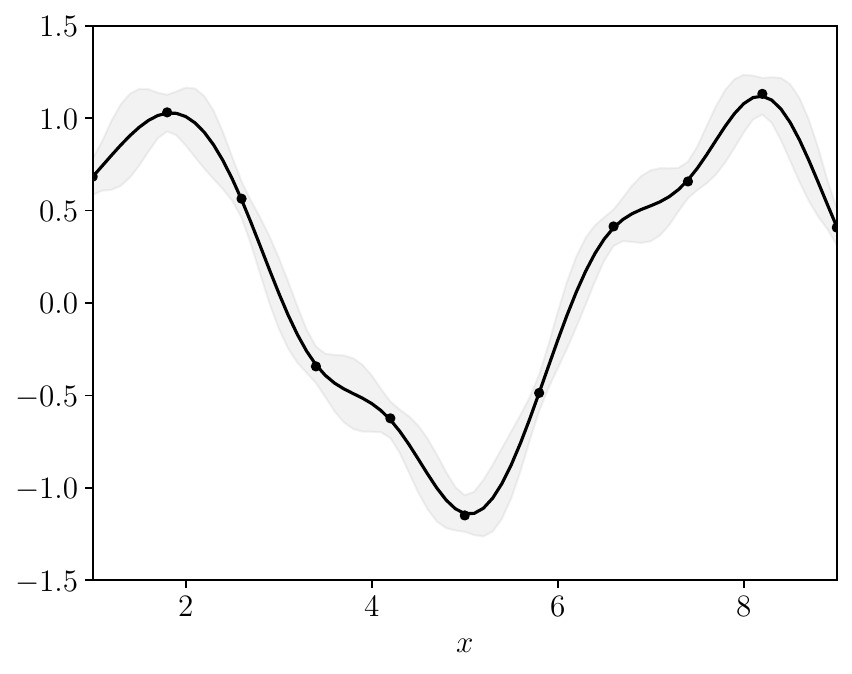} }
	\caption{Gaussian process prior and posterior densities for a one-dimensional regression problem discretised with linear Lagrange basis functions. In (a) we depict three samples from the prior density~$p(\vec s)= \set{N}\left( \vec 0, \vec Q_s^{-1} \right)$ with the Mat\'ern parameters~$\sigma = 0.15$,~$\ell = 1$, and~$\nu = 3.95$. The black line represents the zero mean and the shaded region an offset of two standard deviations~$\pm 2 \sqrt{\diag \left(\vec Q_s^{-1}\right)}$ from the mean. In (b) we show the posterior mean~$\overline{\vec{s}}_{\vert y}$ and the~$95\%$ confidence region obtained by Gaussian process regression from the~$n_y=11$ data points depicted as dots.\label{fig:gpDemo}}
\end{figure}

 As an illustrative example, Figure~\ref{fig:gpDemo} demonstrates the Gaussian process regression on a one-dimensional domain~\mbox{$\widehat{\Omega} = (0,10)$}. The prior density in Figure~\ref{fig:gpDemoPrior} is defined using a zero mean vector and a sparse precision matrix with a non-integer~$\beta$.  In the domain interior the prescribed standard deviation~$\sigma=0.15$ is well approximated by its finite element approximation~$\sqrt{\diag \left(\vec Q_s^{-1} \right ) }$. There are some discrepancies close to the boundaries due to the finite domain effect.  For more details on this example, we refer to  Section~\ref{sec:regressionLine} where we study in depth the learning of the hyperparameters that yield the posterior density shown in Figure~\ref{fig:gpDemoPosterior}. 

Lastly, in engineering applications, there are often repeated readings available at the~$n_y$ observation locations. The respective~$n_o$ observation vectors are collected in the matrix~$\vec Y \in \mathbb R^{n_y \times n_o}$.  The corresponding marginal likelihood~$p(\vec Y)$ for determining the hyperparameters is derived in~\ref{sec:repeatedReadingsAppendix2} . 

%
\section{Statistical finite elements \label{sec:statFEM}}
%
\subsection{PDE-informed prior \label{sec:statFEMprior}}
%
In this section, we choose a Poisson-Dirichlet problem as our forward problem. It is straightforward to consider other kinds of forward problems as demonstrated in the included examples. On a domain~$\Omega \subset \mathbb{R}^d$ with $d \in \{ 1, \, 2, \, 3\}$ and the boundary~$\partial \Omega$, the Poisson-Dirichlet problem is given by
\begin{subequations}	 \label{eq:poisson}
\begin{alignat}{2}
	-  \nabla^2  u(\vec x)  &=  \overline f(\vec s) + s (\vec x)    \qquad &&\text{in }  \Omega  \,,  	 \label{eq:poissonA} \\ 
	u (\vec x ) &= 0 &&\text{on } \partial \Omega  \, , 	 \label{eq:poissonB}
\end{alignat}
\end{subequations}
where~$u(\vec x) \in \mathbb R$ is the unknown solution and the source consists of a deterministic component~$\overline f(\vec x)$ and a random component~$s(\vec x)$.  The random component is a zero-mean Gaussian field with a Mat\'ern covariance function as given by~\eqref{eq:maternField}. The weak form of the forward problem reads: find~$u(\vec x) \in \set H^1(\Omega)$ such that 
\begin{equation}  \label{eq:poissonWeak}
	\int_\Omega \nabla u(\vec x) \cdot  \nabla v(\vec x)  \D \vec x = \int_\Omega \left  (\overline f (\vec x) + s (\vec x) \right ) v(\vec x) \D \vec x \quad \forall v(\vec x) \in \set{H}^1_0 (\Omega)   \, , 
\end{equation}
where~$v(\vec x)$ is a test function and~$\set{H}^1_0 (\Omega) \subset \set{H}^1 (\Omega) $ is the standard Sobolev space of functions satisfying the boundary condition~\eqref{eq:poissonB}.

The discretisation of the forward problem broadly follows the discretisation of the weak form of the stochastic partial equation~\eqref{eq:spdeWeakForm} for the Mat\'ern random field. We use the same finite element mesh and basis functions~$\phi_i (\vec x)$ to discretise both.  Hence, we approximate all functions in the weak form~\eqref{eq:poissonWeak} by 
\begin{subequations}	  \label{eq:discrForward}
\begin{alignat}{2} 
	u(\vec x) & \approx u^h(\vec x)  = \sum_{i=1}^{n_u} \phi_i (\vec x ) u_i \, , \quad && v(\vec x) \approx v^h(\vec x) =  \sum_{i=1}^{n_u} \phi_i (\vec x ) v_i \, ,   \label{eq:discrForwardA} \\ 
	 \overline f(\vec x) & \approx \overline{f}^h(\vec x)  = \sum_{i=1}^{n_u} \phi_i (\vec x ) \overline f_i \, ,  \quad  && s(\vec x) \approx  s^h(\vec x) =  \sum_{i=1}^{n_u} \phi_i (\vec x ) s_i \, .   \label{eq:discrForwardB}
\end{alignat}
\end{subequations}
The resulting discrete system of equations reads 
\begin{equation} \label{eq:forwardDiscrete}
	\vec A \vec u = \vec M \left ( \overline{ \vec f} +\vec s \right ) \, , 
\end{equation}
where~$\vec A$ is the system matrix,~$\vec M$ is the mass matrix and the vectors $\vec u$, $\overline{\vec f}$ and~$\vec s$ collect the respective nodal coefficients in~\eqref{eq:discrForward}. It is worth emphasising that~$\overline{\vec f}$  and~$\vec s$ contain the coefficients in~\eqref{eq:discrForwardB}, are not finite element source vectors in the usual sense, and are different from  similar vectors in our earlier work~\cite{girolami2021statistical}. They become such vectors after multiplication by the mass matrix~$\vec M$.   

As introduced in Section~\ref{sec:fractionalPDE}, the finite element discretisation of the Mat\'ern field~$s(\vec x)$ yields the multivariate Gaussian density~\eqref{eq:overallSPDESolution}, which we repeat for convenience  
 \begin{equation} \label{eq:priorForcing}
	p (\vec s) = \set N \left ( \vec 0, \,  \vec Q_{s}^{-1}  \right )   = \set N \left ( \vec 0, \,  \left (  \vec F_{s,r} ^{-\trans}  \vec F_{s,l}^{\trans} \vec Q_{\tilde{s}}  \vec F_{s,l}  \vec F_{s,r}^{-1} \right )^{-1} \right )  \, .
\end{equation}
Here, the second subscript~$s$ is introduced to distinguish between the different random fields in statFEM which are all represented using the SPDE formulation. Likewise, the parameters of the respective Mat\'ern covariance function~\eqref{eq:maternKernel} and the SPDE~\eqref{eq:fracPDE} are denoted as~\mbox{$\{ \sigma_s, \, \ell_s, \, \nu_s \}$} and~\mbox{$\{\tau_s, \, \kappa_s, \, \beta_s \}$}.  The  linear transformation of the Mat\'ern field~$\vec s$ via~\eqref{eq:forwardDiscrete} yields for the solution~$\vec u $ the multivariate Gaussian density 
 \begin{equation} \label{eq:forwardPrior}
	p (\vec u ) = \set N \left(\overline {\vec u}, \, \vec Q_u^{-1}\right) = \set N \left  (\vec A^{-1} \vec M \overline{\vec f}, \,  \left ( \vec A^\trans \vec M^{-\trans}  \vec Q_s \vec M^{-1} \vec A  \right) ^{-1} \right )  \, ;
\end{equation}
see also~\ref{sec:multivariateNormalCov}.
%

%
%
\subsection{Statistical observation model \label{sec:statFEMobvModel}}
%
As in Gaussian process regression, the observation set~\mbox{$\{ (\vec x_i, \, y_i) \}_{i=1}^{n_y}$} consists of the observations~\mbox{$y_i \in \mathbb R$} at the locations~\mbox{$\vec x_i \in \mathbb R^d$}. The observation model underlying statFEM reads
\begin{align} \label{eq:statFEMmodel}
	\vec y = \vec z + \vec e = \vec P(  \vec u +  \vec d)  + \vec e  \, .
\end{align}
That is, the observed data vector~\mbox{$\vec y \in \mathbb R^{n_y}$} collecting the observations~$\{  y_i \}_{i=1}^{n_y}$ is equal to the unknown true system response~\mbox{$\vec z \in \mathbb R^{n_y}$} plus an observation error, i.e.\ noise,~\mbox{$\vec e \in \mathbb R^{n_y}$}. In turn, the true system response~$\vec z$ consists of the finite element solution~\mbox{$\vec u \in \mathbb R^{n_u}$} of the forward problem and the mismatch error, or model inadequacy,~\mbox{$\vec d \in \mathbb R^{n_u}$}. The observation matrix~\mbox{$\vec P \in \mathbb R^{n_y \times n_u}$} is the same as the one in Gaussian process regression in~\eqref{eq:statMgp}. In the observation model~\eqref{eq:statFEMmodel} all the vectors are random and assumed to be statistically independent.

We approximate the unknown mismatch error~$\vec d$ as a Gaussian field with a Mat\'ern covariance function. The finite element discretisation of the respective SPDE yields according to~\eqref{eq:overallSPDESolution} the multivariate Gaussian density 
\begin{equation} \label{eq:distD}
	p (\vec d) = \set N \left ( \vec 0, \,  \vec Q_{d}^{-1}  \right )   = \set N \left ( \vec 0, \,  \left (  \vec F_{d,r} ^{-\trans}  \vec F_{d,l}^{\trans} \vec Q_{\tilde{d}}  \vec F_{d,l}  \vec F_{d,r}^{-1} \right )^{-1} \right )  \, .
\end{equation}
Herein and in the following, the subscript~$d$ indicates that the respective quantities belong to the mismatch error. As in Gaussian process regression, the hyperparameters~\mbox{$\{ \sigma_d, \, \ell_d, \, \nu_d \}$},  of the Mat\'ern covariance function~\eqref{eq:maternKernel} or the equivalent parameters~\mbox{$\{\tau_d, \, \kappa_d, \, \beta_d \}$}, of the SPDE~\eqref{eq:fracPDE} are learned from the observations.

Finally, the probability density of the independent and identically distributed noise vector~$\vec e$ with a standard deviation~$\sigma_e$ is given by
\begin{equation}
	p(\vec e) = \set N \left(\vec 0, \, \sigma_e^2 \vec I \right) \, .
\end{equation}
%

\subsection{Posterior densities \label{sec:statFEMPosteriorDensities}}
%
Because all the vectors in the assumed observation model~\eqref{eq:statFEMmodel} are Gaussian, the joint density~$p(\vec u, \, \vec y)$ takes the form
\begin{equation}
		p(\vec u, \, \vec y) = \set{N} \left(
		\begin{pmatrix}
			\overline{\vec u} \\   \vec P  \overline{\vec u}
		\end{pmatrix} \, ,   
		\begin{pmatrix}
			\vec Q^{-1}_u &   \vec Q^{-1}_u \vec P^\trans 
			\\ 
			\vec P  \vec Q^{-1}_u &   \vec P \vec Q^{-1}_u \vec P^\trans + \vec Q_{de}^{-1}
		\end{pmatrix}
		\right)
		= \set{N} \left(
		\begin{pmatrix}
			\overline{\vec u} \\   \vec P  \overline{\vec u}
		\end{pmatrix}  \, ,
		\begin{pmatrix}
			\vec Q_u +   \vec P^\trans \vec Q_{de} \vec P  &  -   \vec P^\trans \vec Q_{de}
			\\
			- \vec Q_{de}  \vec P &  \vec Q_{de}
		\end{pmatrix}^{-1}
		\right) \, ,
	\label{eq:statFEMPreJoint}
\end{equation}
where we have introduced the abbreviation 
\begin{equation} \label{eq:qdeStatFEM}
	\vec Q_{de}  = \left ( \vec P \vec Q^{-1}_{d} \vec P^\trans + \sigma_e^2 \vec I  \right )^{-1} \, .
\end{equation}
The sought posterior density~$p(\vec u | \vec y)$ is according to the conditioning property~\eqref{eq:gaussianDensityConditionalFormal} given by 
\begin{equation} \label{eq:statFEpostU}
	p(\vec u \vert \vec y ) = \set N \left(\overline{\vec u}_{\vert y}, \,  \vec Q_{u\vert y}^{-1} \right) \, ,
\end{equation}
where
\begin{subequations}
	\begin{align}
		\overline{\vec u}_{\vert y} &= \overline{\vec u} +\vec Q_{u\vert y}^{-1}  \vec P^\trans \vec Q_{de} \left( \vec y -    \vec P \overline{\vec u} \right)  \, ,
		\label{eq:statFEConditionalMean}
		\\
		\vec Q_{u\vert y} & =  \vec Q_u +    \vec P^\trans \vec Q_{de}  \vec P  \, .
		\label{eq:statFEConditionalCovariance} 
		\end{align}
\end{subequations}
The inverse, or the factorisation, of the precision matrix~$\vec Q_{u\vert y}$ is needed to determine the posterior mean~$\overline{\vec u}_{\vert y} $ and the components of  the covariance~\mbox{$\vec C_{u\vert y}= \vec Q_{u\vert y}^{-1}$}. The precision matrix~$\vec Q_{u\vert y}$  is sparse when both~$ \vec Q_u$ and~$\vec P^\trans \vec Q_{de}  \vec P$ are sparse. As discussed, the precision matrix~$\vec Q_{u}$ is sparse when~$\beta_s \in \mathbb N$. On the other hand, the matrix~\mbox{$\vec P^\trans \vec Q_{de}  \vec P$} is sparse when~$n_y \ll n_u$, in which case we can rewrite~\eqref{eq:qdeStatFEM} using the Sherman-Morrison-Woodbury formula as 
\begin{equation}
	\vec Q_{de}  = \left ( \vec P \vec Q^{-1}_{d}  \vec P^\trans + \sigma_e^2 \vec I  \right )^{-1} =  \frac{1}{\sigma_e^2}  \vec I - \frac{1}{\sigma_e^2} \vec P \left( \sigma_e^2  \vec Q_{d}  +  \vec P^\trans  \vec P  \right)^{-1}  \vec P^\trans    \, .
\end{equation}

In most engineering applications, there is limited data so that the assumption~$n_y \ll n_u$ appears sensible. If, however,~$n_y \approx n_u$, we can compute as a possible remedy a low-rank approximation to the sparse matrix~\mbox{$ \sigma_e^2  \vec Q_{d}  +  \vec P^\trans  \vec P $}. Alternatively, the observation model~\eqref{eq:statFEMmodel} can be replaced by \mbox{$\vec y = \vec P \vec u + \vec e$} omitting the discrepancy term~$\vec d$. In that case, we can define the auxiliary (latent) variable~\mbox{$\vec t = \vec F_{s,r}^{-1} \vec M^{-1} \vec A$} and follow an approach similar to Gaussian process regression in~Section~\ref{sec:gpPrecision}. Briefly, the suggested change of observation model mainly pertains to the modelling of epistemic versus aleatoric uncertainties. In the limit of infinite data the inferred posterior~$\vec u$ of the suggested alternative model will converge to a deterministic vector. In contrast, the posterior~$\vec z$ in the model~\eqref{eq:statFEMmodel} remains a random vector. For a statistical finite element formulation based on the alternative model see \cite{akyildiz2022statistical}.  In this paper we do not consider the case~$n_y \approx n_u$ further. 

The marginal likelihood~$p(\vec y)$ is used to learn the hyperparameters~$\{\tau_d, \, \kappa_d, \, \beta_d \}$  of the model inadequacy~$\vec d$ in the SPDE~\eqref{eq:fracPDE}. According to the marginalisation property of Gaussians, see~\ref{sec:multivariateNormalConditional}, we can read from the joint density~\eqref{eq:statFEMPreJoint}  that 
\begin{equation}
	p(\vec y)  = \set{N}\left(  \vec P \overline{\vec u}  , \,  \vec P \vec Q_u^{-1} \vec P^\trans + \vec Q_{de}^{-1} \right)  \, .
\end{equation}
Taking its logarithm gives
\begin{equation} \label{eq:statFElogMargy}
	 \log p(\vec y) =  - \frac{1}{2} \left( \vec y -  \vec P \overline{\vec u} \right)^\trans \left( \vec P \vec Q_u^{-1} \vec P^\trans +  \vec Q_{de}^{-1} \right)^{-1}  \left( \vec y -  \vec P \overline{\vec u} \right) -
	\frac{1}{2} \log \det \left( \vec P \vec Q_u^{-1} \vec P^\trans +  \vec Q_{de}^{-1} \right) - \frac{n_y}{2} \log 2\pi  \, .
\end{equation}
To avoid any dense matrix operations the terms on the right-hand side are rewritten using the Sherman-Morrison-Woodbury formula and the matrix determinant lemma as follows
\begin{subequations}
\begin{align}
	\left( \vec P \vec Q_u^{-1} \vec P^\trans + \vec Q_{de}^{-1} \right)^{-1} &= \vec Q_{de} - \vec Q_{de}  \vec P \vec Q_{u\vert y}^{-1} \vec P^\trans \vec Q_{de}  \, ,   \\
	\det \left( \vec P \vec Q_u^{-1} \vec P^\trans + \vec Q_{de}^{-1} \right) &= \frac{  \det \left(   \vec Q_u +  {\vec P}^\trans   \vec Q_{de} \vec P  \right)  \det \left( \vec Q_{de}^{-1} \right) }{  \det \left( \vec Q_u \right)  }  \, ,  \\ 
	\det \left( \vec Q_{de}^{-1} \right) &=  \left(\sigma_e^2 \right)^{n_y -  n_u } \frac{\det \left( \sigma_e^2  \vec Q_{b}  + \vec F_{d,r}^\trans \vec P^\trans \vec P \vec F_{d,r} \right)}{\det \left( \vec Q_{b}  \right)}   \, .
\end{align}
\end{subequations}
See~\ref{sec:repeatedReadingsAppendix3} for the marginal likelihood~$p(\vec Y)$ in case of repeated readings when there is more than one observation vector available.  

Finally, with the obtained posterior density~$p(\vec u \vert \vec y)$  and the hyperparameters learned by maximising the marginal  likelihood~$p(\vec y)$, or~$p(\vec Y)$, the posterior density of the true system response~$\vec z$ is given by
\begin{equation} \label{eq:posteriorTrueResponse}
	p( \vec z \vert \vec y ) = \set{N}\left( \overline{\vec z}_{\vert y} , \vec Q_{z\vert y}^{-1} \right)   =  \set{N}\left( \vec P \overline{\vec u}_{\vert y} , \vec P \vec Q_{u\vert y}^{-1} \vec P^\trans  + \vec P \vec Q_d \vec P^\trans \right)   \, .
\end{equation}
Furthermore, as an implementation note, it is numerically more stable to sequence the computation of the posterior densities and the marginal likelihood so that  only matrix vector products and matrix solve operations are required.

%
\section{Examples \label{sec:examples}}
%
We proceed to demonstrate the utility of the SPDE representation of random fields in large-scale Gaussian process regression and statistical finite element analysis on $1$D and $2$D domains and $2$-manifolds. The SPDE is discretised using either linear Lagrange or Loop subdivision basis functions~\cite{Cirak:2000aa, Cirak:2002aa}.  Unless otherwise stated, the SPDE domain~$\widehat{\Omega}$ is chosen the same as the computational domain~$\Omega$ and its boundary condition is homogeneous Neumann.  As mentioned in Section~\ref{sec:fractionalPDE}, the solution of the stochastic PDE is a Mat\'ern field only when~$\widehat{\Omega}$ is infinitely large. However, in our experience choosing a relatively compact~$\widehat \Omega$ has a negligible impact on the obtained posterior densities when the covariance length-scale is sufficiently small, i.e.~$\ell  \ll |\Omega |$. This is possibly because some finite domain effects are mitigated by the choice of hyperparameters~\cite{lindgren2011explicit}. A sensible alternative approach not pursued here is to choose Robin boundary conditions and to determine the respective weighting factor together with all other hyperparameters by maximising the marginal likelihood. See also~\cite{roininen2014whittle,daon2016mitigating,khristenko2019analysis} for discussions on the choice of the SPDE boundary conditions. In particular, the choice of Robin boundary conditions imitating an infinite domain is proposed in~\cite{daon2016mitigating}. In computational mechanics, such boundary conditions are often referred to as Dirichlet-to-Neumann boundary conditions~\cite{givoli1999recent}. In all the examples we learn the hyperparameters by maximising the marginal likelihood using the BOBYQA algorithm in the NLopt library~\cite{powell1994direct,johnson2014nlopt}. The marginal likelihood is often non-convex so that its maximisation requires some care.
%
\subsection{Convergence of the Mat\'ern covariance function \label{sec:maternConvergence}}
%
To begin with, we investigate the convergence of the approximate Mat\'ern covariance function~$c_s^h(x , \,x ')$ obtained by discretising the SPDE~\eqref{eq:fracPDE} using linear Lagrange basis functions on uniformly refined finite element meshes.  The computational domain for discretisation is chosen as~\mbox{$\widehat{\Omega} = ( - 0.2, \,  1.2 )$} while for evaluating the $L^2$ norms only the smaller domain~\mbox{$\Omega = ( 0, \,  1 )$} is considered. The integer and non-integer exponent values are~$\beta \in \{ 1,2,3,1.25,1.75 \}$, the standard deviation is~$\sigma = 1$ and the length-scale parameters is~$\ell = 0.05$. Observe that~\mbox{$\widehat{\Omega} = ( -4 \ell, \,  1 + 4 \ell )$} so that the effect of the domain size on the SPDE solution in~\mbox{$\Omega = ( 0, \,  1 )$} can be neglected.

\begin{figure}
	\centering
	\subfloat[][$\beta \in \{ 1,2,3 \} $ \label{fig:maternConvergenceA}] {
		\includegraphics[width=0.315 \textwidth]{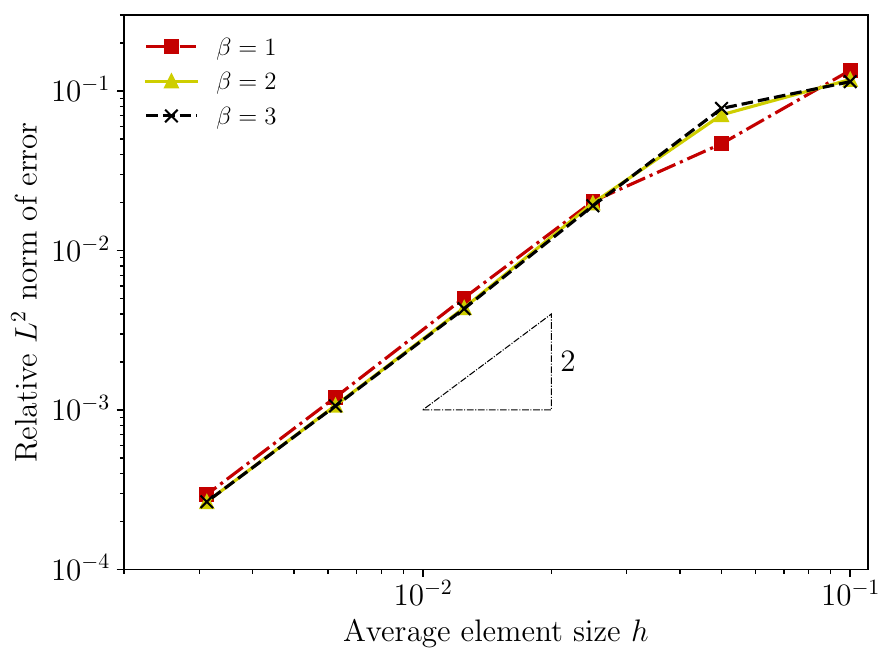} }
	\hfil
	\subfloat[][$\beta = 1.25 $ \label{fig:maternConvergenceB}] {
		\includegraphics[width=0.315 \textwidth]{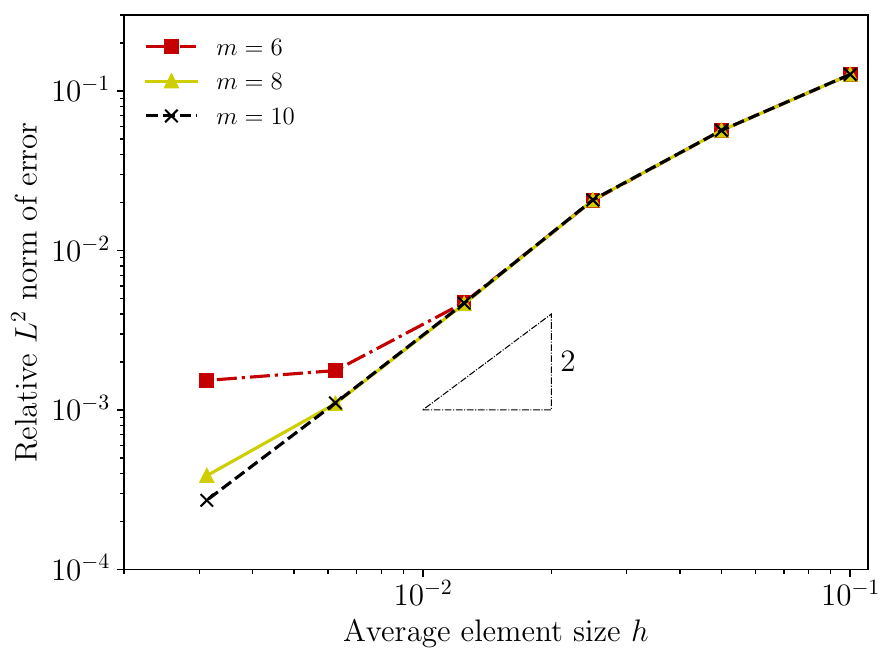} }
	\hfil
	\subfloat[][$\beta = 1.75 $ \label{fig:maternConvergenceC}] {
		\includegraphics[width=0.315 \textwidth]{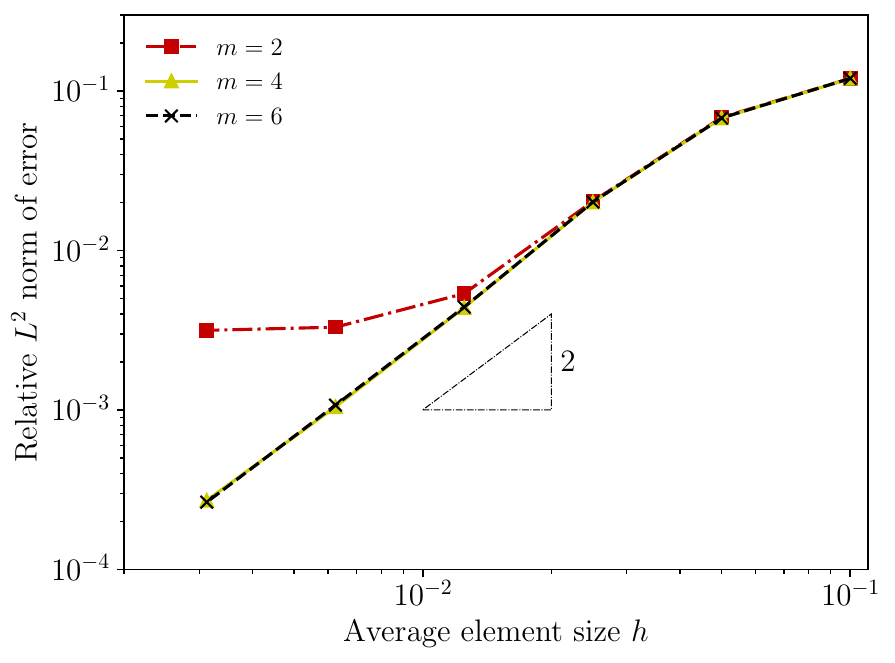} }
	\caption{Convergence of the finite element approximation of the Mat\'ern covariance function~$c_s^h( x, x')$ with~$x' = 0.5$ in one dimension. \label{fig:maternConvergence}}
\end{figure}

The exact Mat\'ern covariance function~$c_s(x , \,x ')$ is given by~\eqref{eq:maternKernel} and has according to~\eqref{eq:overallSPDESolution} the finite element approximation 
\begin{equation}
	c_s^h(x, x') = \vec \phi(x)^\trans \vec Q_s^{-1} \vec \phi(x') \, , 
\end{equation}
where~$\vec \phi(x)$ is the vector of finite element basis functions. The relative $L^2$ norm error of the approximation over the domain of interest~$\Omega$ is 
\begin{equation}
	\eta (x')= \frac{ \| c_s^h(x, x') - c_s(x, x') \|_{2} }{ \| c_s(x,x') \|_{2}  } \,.
\end{equation}
Figure~\ref{fig:maternConvergence} depicts the convergence of the relative error~$\eta (x' = 0.5)$ for different exponent values~$\beta$. The convergence rate for the integer-valued exponents~$\beta \in \{1,2,3\}$ shown in Figure~\ref{fig:maternConvergenceA} is~$2$ which is optimal for linear basis functions. The convergence of the non-integer exponents~$\beta = 1.25$ and~$\beta = 1.75$  shown in Figures~\ref{fig:maternConvergenceB} and~\ref{fig:maternConvergenceC} depend on the polynomial degree~$m$ of the polynomials in the rational interpolation~\eqref{eq:ratInterp}. As discussed in Section~\ref{sec:arbiraryExponent}, choosing a large~$m$ leads to an increasingly denser precision matrix~$\vec Q_s$ requiring~$4m$ sparse matrix multiplications for its construction. We re-emphasise that rational interpolation is only needed when~$\beta$ is non-integer. According to Figures~\ref{fig:maternConvergenceB} and~\ref{fig:maternConvergenceC} the convergence depends both on the finite element and rational approximation errors. That is, finer meshes require  higher rational approximation orders and vice versa. A practical approach for choosing the lowest possible rational approximation order~$m$ is to monitor convergence with increasing mesh refinement. 
%

%
\subsection{Gaussian process hyperparameter learning on one-dimensional domain \label{sec:regressionLine}}
%
We consider next Gaussian process hyperparameter learning on the domain~$\Omega = \widehat{\Omega} = (0, 10)$. We discretise the domain uniformly using~$101$ linear Lagrange basis functions. We assume that the data vectors used in this example are sampled from the multivariate normal density
\begin{equation}
	p(\vec y) = \set{N}\left( \vec 0, \vec P \vec Q_s^{-1} \vec P^\trans + \sigma_e^2 \vec I \right) \,,
\end{equation}
with the observation error~$\sigma_e = 0.05$. The precision matrix~$\vec Q_s$ is determined using the Mat\'ern parameters~\mbox{$\sigma = 0.15$}, $\ell = 1$, and~$\nu = 3.95$. The chosen smoothness~$\nu$ corresponds to the exponent $\beta = 2.225$. Since the exponent value is non-integer, we use rational interpolation and choose the polynomial degree as~$m = 6$. 

As illustrated earlier in Figure~\ref{fig:gpDemoPrior}, the random samples drawn from the prior~$p(\vec s) \sim \set{N}\left(\vec 0, \vec Q_s^{-1} \right)$ fluctuate smoothly about the zero mean. In the domain interior, the prescribed standard deviation $\sigma = 0.15$ is accurately approximated by the finite element interpolated standard deviation~$\sqrt{\diag\left(\vec Q_s^{-1} \right)}$. However, due to the boundary effect the standard deviation increases within a distance of approximately~$\ell = 1$ from the boundaries.

\begin{figure}
	\centering
	\hfil
	\subfloat[][$n_y = 11$, $n_o = 1$ \label{fig:gp1DLogLikelihood11}] {
		\includegraphics[height=0.315\textwidth]{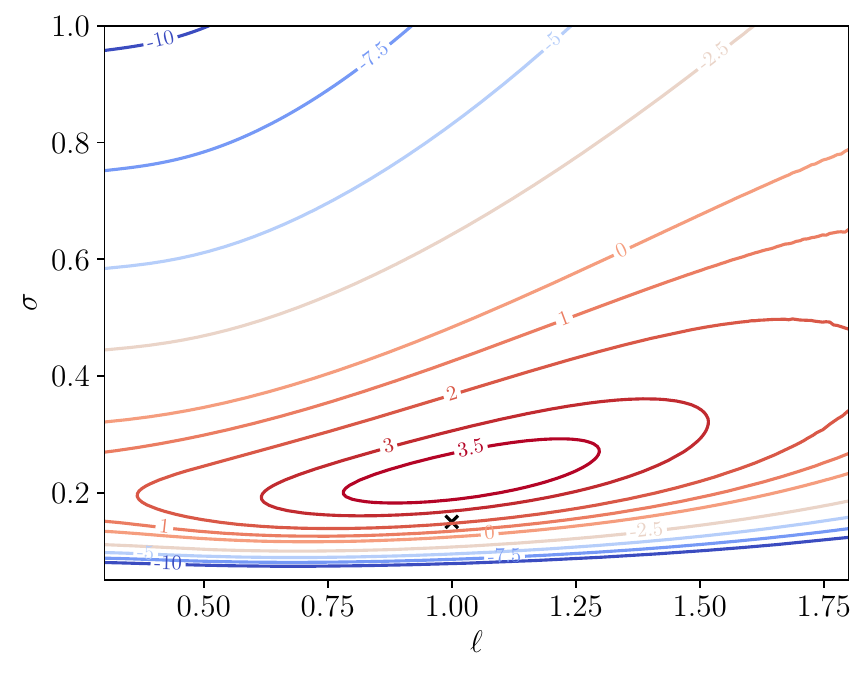} }
	\hfil
	\subfloat[][$n_y = 81$, $n_o = 1$ \label{fig:gp1DLogLikelihood81}] {
		\includegraphics[height=0.315\textwidth]{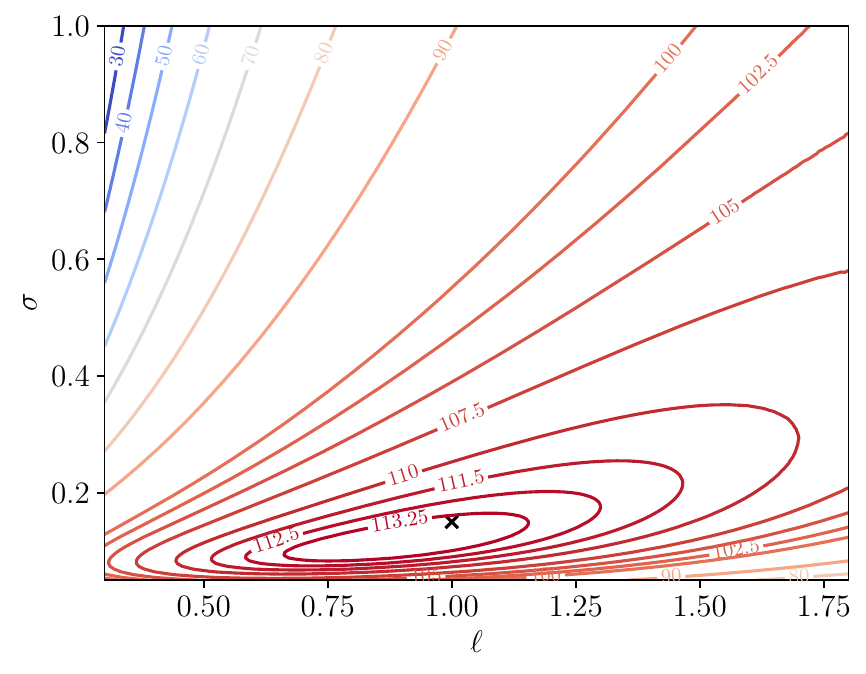} }
	\caption{Gaussian process hyperparameter learning on a one-dimensional domain. Isocontours of the log marginal likelihood $\log p(\vec y)$ computed using the sampled data vector~$\vec y$. The cross $\times$ represents the target hyperparameters~$\sigma = 0.15$ and~$\ell = 1$. \label{fig:gp1DLogLikelihood}}
\end{figure}

\begin{figure}
	\centering
	\hfil
	\subfloat[][Standard deviation~$\sigma^*$ \label{fig:gp1DHyperparamsSigma}] {
		\includegraphics[width=0.4 \textwidth]{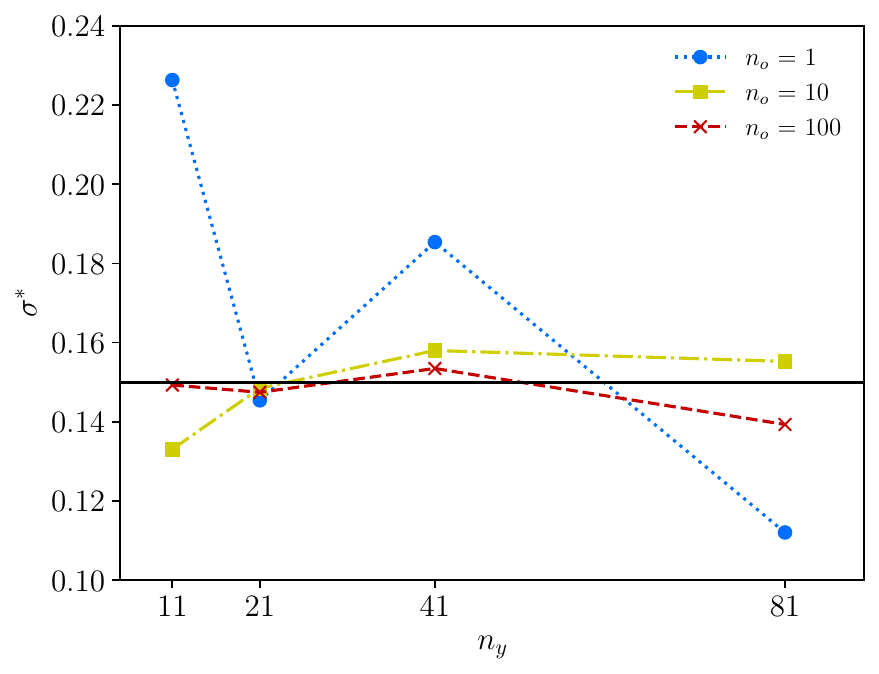} }
	\hfil
	\subfloat[][Length-scale~$\ell^*$ \label{fig:gp1DHyperparamsLength}] {
		\includegraphics[width=0.4 \textwidth]{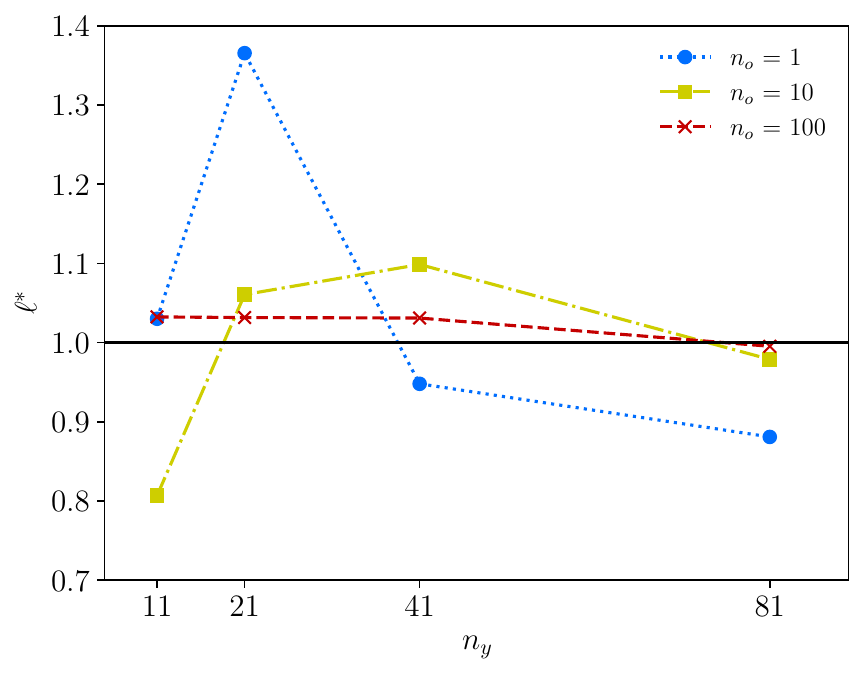} }
	\caption{Gaussian process hyperparameter learning on a one-dimensional domain. Point estimates of the hyperparameter vector~$\vec w^* = \left( \sigma^* \,\, \ell^* \right)^\trans$ obtained by maximising the log marginal likelihood. The respective black lines indicate the data generating standard deviation~$\sigma = 0.15$ and length-scale~$\ell = 1$. \label{fig:gp1DHyperparams}}
\end{figure}

We consider~$n_y \in \{11, 21, 41, 81 \} $ data points that are spaced evenly on the subdomain~$[1, 9] \subset \Omega$. For each~$n_y$, we sample $n_o = \{1, 10, 100\}$ independent readings from~$p(\vec y)$. We aim to examine if the hyperparameter vector~\mbox{$\vec w = \begin{pmatrix}	\sigma & \ell \end{pmatrix}^\trans$} can be learned effectively using maximum likelihood estimate. Note that we have excluded the smoothness $\nu$ from the hyperparameter vector~$\vec w$ to avoid non-identifiability of the Mat\'ern hyperparameters in GP regression~\cite{zhang2004inconsistent}. Figure~\ref{fig:gp1DLogLikelihood} shows the log marginal likelihood~$\log p(\vec y)$ in the vicinity of the prescribed standard deviation~$\sigma = 0.15$ and length-scale~$\ell = 1$. As visible, the local maximum of~$\log p(\vec y)$ becomes numerically closer to the prescribed hyperparameters when~$n_y$ increases. Figure~\ref{fig:gp1DHyperparams} shows that the point estimates of both~$\sigma^*$ and~$\ell^*$ converge to the respective prescribed values~$\sigma = 0.15$ and~$\ell = 1$ when both $n_y$ and~$n_o$ increase. Remarkably, by increasing~$n_o$ the hyperparameter learning yields accurate estimates even when~$n_y$ is small.

%
\subsection{Gaussian process regression on a hemispherical shell \label{sec:regressionHemisphere}}
\begin{figure}
	\centering
	\subfloat[][{Open hemispherical shell in~$\mathbb{R}^3$ \label{fig:hemisphereIntroductionA} }] {
		\includegraphics[width=0.4\textwidth]{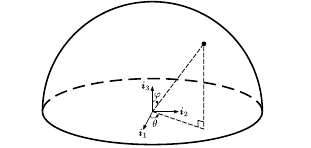} } 
	\hfill
	\subfloat[][{Spherical harmonic \label{fig:hemisphereIntroductionB}}] {
		\includegraphics[width=0.235\textwidth]{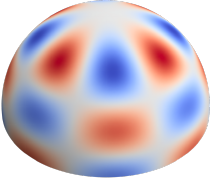} }
	\hfill
	\subfloat[][{Observation points~$n_y = 41$ \label{fig:hemisphereIntroductionC}}] {
		\includegraphics[width=0.22\textwidth]{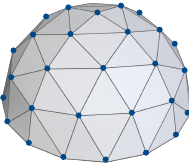} } 
	\caption{Gaussian process regression on a hemispherical shell. We discretise the open hemispherical shell shown in (a) with~$n_{\text{el}} = 1048576$ linear triangular elements. We obtain the deterministic data~$\vec y$ with $n_y = \{ 41, 145, 545, 2113, 8321 \}$ by evaluating~$s ( \theta, \, \varphi )$ in (b) at the~$5$ sets of mesh vertices. The coarse mesh and the first set of observation points with~$n_y = 41$ are shown in (c). By repeatedly quadrisecting the mesh we obtain the remaining sets of observation points. \label{fig:hemisphereIntroduction}}
\end{figure}

%
We  study next the convergence of Gaussian process regression on the hemispherical shell depicted in Figure~\ref{fig:hemisphereIntroductionA}. We choose as the function to infer the (real-valued) spherical harmonic 
\begin{equation}
	 s ( \theta, \, \varphi )=  \frac{3}{1024} \sqrt{\frac{1309}{\pi}} \cos( 4 \theta ) \sin^4( \varphi)  \left( 99 + 156 \cos(2 \varphi ) + 65 \cos( 4 \varphi ) \right) \,,
\end{equation}
where~$\theta $ and~$ \varphi$ are the polar coordinates of the surface points, see Figure~\ref{fig:hemisphereIntroductionA}. The SPDE for computing the respective prior distribution is defined on the shell and is discretised with a mesh consisting of~$n_u = 525313$ nodes and~$n_{\text{el}} = 1048576$ linear triangular elements with an average size of~$h = 5.3624 \times 10^{-3}$.  The mesh is obtained by starting from the coarse mesh depicted in Figure~\ref{fig:hemisphereIntroductionC} and successively subdividing each triangle into four triangles. The new nodes are projected onto the surface of the shell. The parameters of the SPDE are fixed throughout the convergence study and are equal to the Mat\'ern parameters~\mbox{$\sigma = 0.1$},~\mbox{$\ell = 0.4$} and~\mbox{$\nu = 1$}. The finite element discretisation of the SPDE according to Section~\ref{sec:arbiraryExponent} yields the prior for~$\vec s \in \mathbb R^{n_u}$. For the Mat\'ern parameter~\mbox{$\nu = 1$} the SPDE is non-fractional so that the prior is especially easy to compute.  

\begin{figure}
	\centering
	\subfloat[][{$n_y = 41$}] {
		\includegraphics[width=0.22\textwidth]{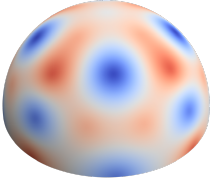} }
	\hfill
	\subfloat[][{$n_y = 145$}] {
		\includegraphics[width=0.22\textwidth]{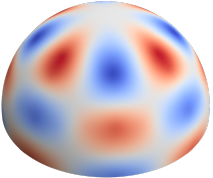} }
	\hfill
	\subfloat[][{$n_y = 545$}] {
		\includegraphics[width=0.22\textwidth]{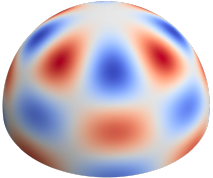} }
	\hfill
	\subfloat[][{$n_y = 8321$}] {
		\includegraphics[width=0.22\textwidth]{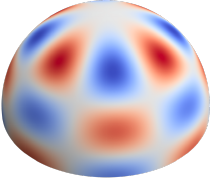} }
	\caption{Gaussian process regression on a hemispherical shell with~$n_{y} = \{ 41, 145, 545, 8321\}$ observation data. The respective posterior means~$\overline{s}_{\vert y}(\vec x)$ are obtained with~$n_{\text{el}} = 1048576$ linear triangular elements and the hyperparameters are~$\sigma = 0.1$, $\ell = 0.4$, $\nu = 1$ and $\sigma_e = 1 \times 10^{-10}$. All plotted field values range between $-0.68$ (blue) and $0.68$ (red). \label{fig:hemisphereMean}}
\end{figure}

\begin{figure}
	\centering
	\includegraphics[width=0.575 \textwidth]{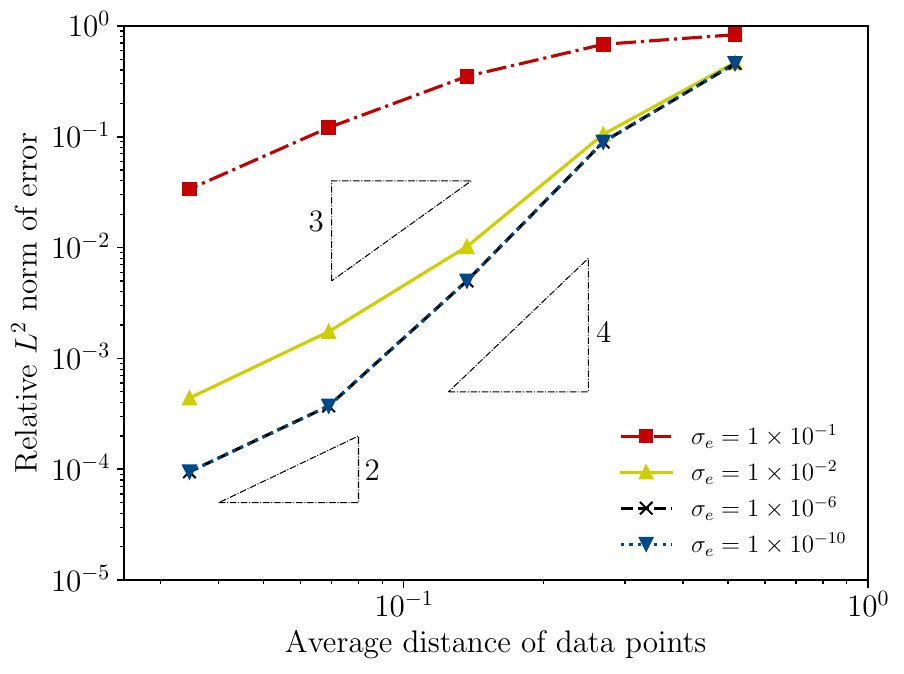}
	\caption{Gaussian process regression on a hemispherical shell. Convergence of the posterior means~$\overline{s}^h_{\vert y}(\vec x)$ with respect to the average distance of data points. In all cases, we prescribe the Mat\'ern hyperparameters as~$\sigma = 0.1$,~$\ell = 0.4$ and $\nu = 1$. For each~$n_y$, we define the average distance of data points as the average edge length of the triangulation, see Figure~\ref{fig:hemisphereIntroductionC}. \label{fig:convergenceHemisphere}}
\end{figure}

We consider five different observation vectors~$\vec y \in \mathbb R^{n_y}$, which are sampled at~$n_y \in \{ 41, 145, 545, 2113, 8321 \}$ nodes of the finite element mesh. The observation points correspond to the nodes of the original coarse mesh and the once, twice and so on subdivided meshes.  Consequently, the observation points are nearly uniformly distributed across the hemisphere. As the standard deviation for the observation noise we choose the values~\mbox{$\sigma_e \in \left \{10^{-1}, 10^{-2}, 10^{-6}, 10^{-10} \right \}$}.  The posterior distribution of~$\vec s$ is given by~\eqref{eq:sparsePrecPosterior}. The posterior mean vector~\mbox{$\overline{\vec s}_{\vert y} \in \mathbb R^{n_u}$} and precision matrix~\mbox{$\vec Q_{s\vert y} \in \mathbb R^{n_u \times n_u}$} are only defined with respect to the nodes of the finite element mesh.  Focusing in the following on the posterior mean, we obtain a piecewise linear representation by interpolating with the finite element basis functions~$\vec \phi(\vec x)$, i.e., 
\begin{equation}
	\overline{s}_{\vert y}^h (\vec x)  =  \vec \phi(\vec x)^\trans \overline{\vec s}_{\vert y} \, .
\end{equation}
As illustrated in Figure~\ref{fig:hemisphereMean} the posterior mean~$\overline{s}_{\vert  y}^h (\vec x)$ converges as expected to~$s(\vec x)$ shown in Figure~\ref{fig:hemisphereIntroductionB} with increasing~$n_y$.  Furthermore, we consider the relative~$L^2$ norm of error defined as  
\begin{equation}
	\frac{\| s (\vec x) - \overline{s}^h_{\vert y}( \vec x ) \|_{2}}{\| s(\vec x) \|_{2} } \,.
\end{equation}
Figure~\ref{fig:convergenceHemisphere} depicts the convergence of the~$L^2$ norm of error with respect to the average distance between the data points. As expected, the convergence depends on the observation error~$\sigma_e$, and a large ~$\sigma_e$ leads to a slow convergence on coarse meshes. In all cases the asymptotic convergence rate is about~$2$, which is in agreement with analytical error estimates reported in the literature~\cite{teckentrup2020convergence, karvonen2022error}. The wall-clock time for computing the posterior mean vector is approximately~$52.7 s$ on a MacBook Pro with an M1 Pro chip and 32 GB RAM using a single core and is independent of~$n_y$ and~$\sigma_e$.

We note that it is feasible to consider with the introduced sparse precision formulation up to~$n_y = n_u$ observation points without any problems.  In contrast, using the standard covariance formulation  the number of possible observation points~$n_y$ is very limited because of the need to factorise a dense covariance matrix of dimension~$n_y \times n_y$. 
%
\subsection{Statistical finite element analysis of Poisson-Dirichlet problems \label{sec:statFETwoD}}
\begin{figure}
	\centering
	\subfloat[][{Domain and partition of boundary \label{fig:plateProblemDescription}}] {
		\includegraphics[width=0.315\textwidth]{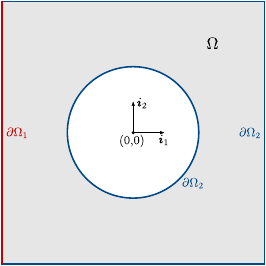} }
	\hfil
	\subfloat[][{Initial coarse mesh  \label{fig:plateProblemInitialMesh}}] {
		\includegraphics[width=0.315\textwidth]{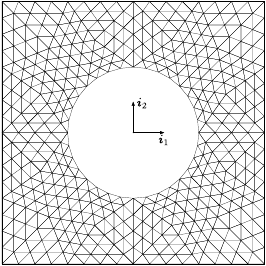} }
	\caption{Poisson-Dirichlet problem on a square plate with a hole. As shown in (a), we partition the boundary~$\partial \Omega$ into the non-overlapping sets~$\partial \Omega = \partial \Omega_1 \cup \partial \Omega_2$. We study first the case where both~$\partial \Omega_1$ and~$\partial \Omega_2$ are non-random and subsequently the case where only~$\partial \Omega_2$ is non-random. The initial coarse mesh in (b) has~$432$ vertices.  \label{fig:plateProblem}}
\end{figure}

\begin{figure}
	\centering
	\subfloat[][{$n_y = 30$}] {
		\includegraphics[width=0.22\textwidth]{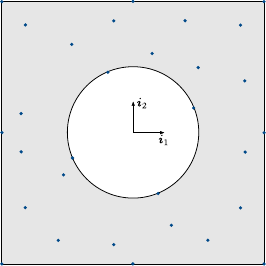} }
	\hfill
	\subfloat[][{$n_y = 60$}] {
		\includegraphics[width=0.22\textwidth]{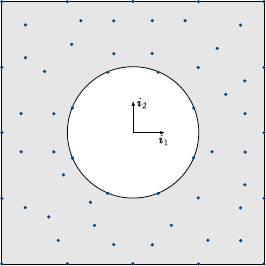} }
	\hfill
	\subfloat[][{$n_y = 120$}] {
		\includegraphics[width=0.22\textwidth]{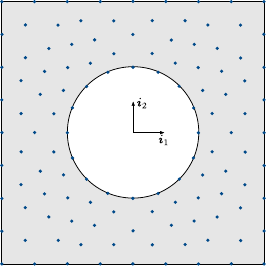} }
	\hfill
	\subfloat[][{$n_y = 432$}] {
		\includegraphics[width=0.22\textwidth]{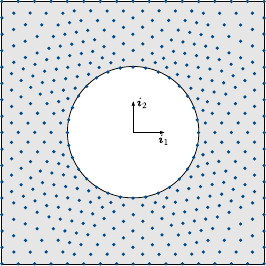} }
	\caption{Poisson-Dirichlet problem on a square plate with a hole. The dots depict the~$n_y$ data points chosen randomly from the vertices of the initial coarse mesh in Figure~\ref{fig:plateProblemInitialMesh}. \label{fig:plateDataLocation}}
\end{figure}

%
The two examples in this section concern a square domain~$\Omega$ of side length~$2$ with a centred circular hole of radius~$0.5$, see Figure~\ref{fig:plateProblemDescription}. In both examples we consider a Poisson-Dirichlet problem with a random source which we discretise with a finite element mesh consisting of~$49152$ linear triangular elements and~$n_u = 24960$ nodes. The mesh is obtained by successively subdividing the triangles in the mesh shown in Figure~\ref{fig:plateProblemInitialMesh} into four triangles. The new nodes on the hole boundary are projected onto the circle. In both examples we assume a true system response~$ z(\vec x)$ different from the solution of the Poisson-Dirichlet problem~$ u(\vec x)$ to induce a model inadequacy~$d (\vec x)$. As illustrated in Figure~\ref{fig:plateDataLocation}, we consider four sets of observation locations with~\mbox{$n_y \in \{30, 60, 120, 432\}$} points located at the nodes of the finite element mesh. At each observation point we record $n_o \in \{2, 10, 20\}$ readings.  In the second problem in Section~\ref{sec:statFETwoDBoundary}, in addition to the random source the Dirichlet boundary condition on one of the sides is also random.

%
\subsubsection{Source uncertainty \label{sec:statFETwoDSource}}
\begin{figure}[t]
	\centering
	\subfloat[][{True source:~$\overline{\chi}^h(\vec x)$ and~$\sqrt{c_\chi^h(\vec x, \vec x)}$ \label{fig:plateSourceTrueSource}}] {
		\includegraphics[width=0.22\textwidth]{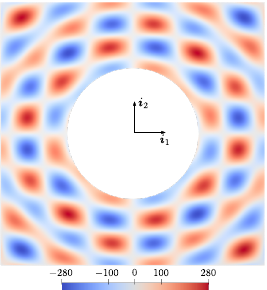} \hspace{0.75em}
		\includegraphics[width=0.22\textwidth]{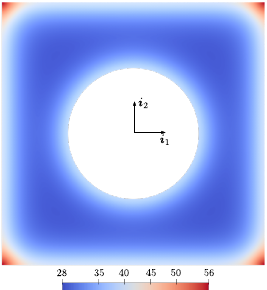}
	}
	\hfil
	\subfloat[][{True response:~$\overline{z}^h(\vec x)$ and~$\sqrt{c^h_z(\vec x, \vec x)}$ \label{fig:plateSourceTrueSolution}}] {
		\includegraphics[width=0.22\textwidth]{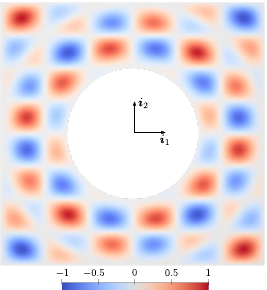} \hspace{0.75em}
		\includegraphics[width=0.22\textwidth]{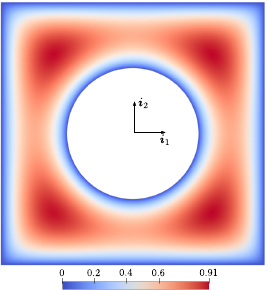}
	}
	\\
	\subfloat[][{Prior source:~$\overline{f}^h(\vec x)$ and~$\sqrt{c^h_f(\vec x, \vec x)}$ \label{fig:plateSourcePriorSource}}] {
		\includegraphics[width=0.22\textwidth]{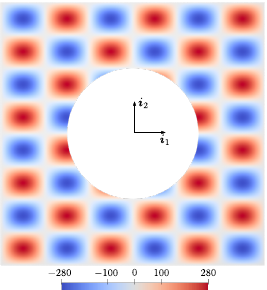} \hspace{0.75em}
		\includegraphics[width=0.22\textwidth]{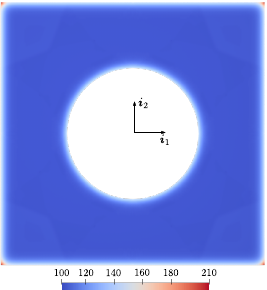}
	}
	\hfil
	\subfloat[][{Prior solution:~$\overline{u}^h(\vec x)$ and~$\sqrt{c^h_u(\vec x, \vec x)}$ \label{fig:plateSourcePriorSolution}}] {
		\includegraphics[width=0.22\textwidth]{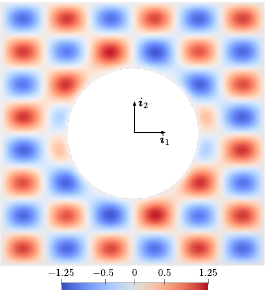} \hspace{0.75em}
		\includegraphics[width=0.22\textwidth]{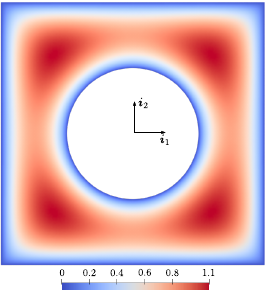}
	}
	\caption{Poisson-Dirichlet problem with source uncertainty. Mean and standard deviation of the true system response~$z^h(\vec x)$ and the finite element prior~$u^h(\vec x)$. \label{fig:plateSourceDescription}}
\end{figure}

%
\paragraph{True system response} For generating synthetic observation data we assume a true solution~$z(\vec x)$ which is the solution of the Poisson-Dirichlet problem
\begin{subequations} \label{eq:pwhExactPDE}
	\begin{alignat}{2}
		- \Delta z(\vec x) &=  \chi(\vec x)  \quad && \text{in } \Omega \,,
		\\
		z( \vec x ) &= 0  \quad && \text{on } \partial \Omega \, ,
	\end{alignat}
\end{subequations}
with the source term~$\chi(\vec x) =  \overline {\chi} (\vec x) +  \chi_s(\vec x)$ consisting of the deterministic source 
\begin{equation} \label{eq:pwhExactDeterm}
	\overline{\chi}(\vec x) = - \left (\frac{\partial^2 }{\partial x_1^2} + \frac{\partial^2}{\partial x_2^2} \right )  \left( \sin( 2 \pi \| \vec x \| ) \sin( 3 \pi x_1 ) \sin( 4 \pi x_2 ) \right) \,,
\end{equation}
and the zero-mean random source satisfying the SPDE
\begin{subequations}  \label{eq:pwhExactSPDE}
	\begin{alignat}{2}
		\left (  50   - \Delta \right ) \chi_s(\vec x) &= 280 \sqrt{2 \pi}\, g (\vec x) \quad && \text{in } \Omega \,,
		\\
		\nabla  \chi_s(\vec x) \cdot \vec n  &= 0  \quad && \text{on } \partial \Omega \,.
	\end{alignat}
\end{subequations}
The SPDE has homogeneous Neumann boundary conditions and~$\vec n$ is the outer normal to the domain~$\Omega$. The finite element discretisation of the SPDE~\eqref{eq:pwhExactSPDE} and the PDE~\eqref{eq:pwhExactPDE} consist of a mesh with~$n_z = 394752$ nodes and~$786432$ linear triangular elements yielding the multivariate Gaussian distributions 
\begin{equation}
	p(\vec \chi) = \set{N}\left( \overline{\vec \chi},  \vec Q_\chi^{-1} \right), \quad p(\vec z) = \set{N}\left( \overline{\vec z},  \vec Q_z^{-1} \right) =  \set{N}\left( \vec A^{-1} \vec M \overline{\vec \chi}, \vec A^{-1} \vec M \vec Q_\chi^{-1}  \vec M^\trans \vec A^{-\trans} \right)  \, ;
\end{equation}
cf.~\eqref{eq:forwardPrior}. These two multivariate Gaussian distributions give rise to the Gaussian processes 
\begin{subequations}
\begin{align}
	\chi^h(\vec x) & \sim \set{GP} \left (\overline{\chi}^h, \,  c_\chi^h (\vec x, \, \vec x') \right ) = \set{GP} \left (\vec \phi(\vec x)^\trans \overline{\vec \chi}, \,   \vec \phi(\vec x)^\trans  \vec Q_\chi^{-1}\vec \phi(\vec x')  \right ) \,  ,  \\
	z^h(\vec x) & \sim \set{GP} \left ( \overline{z}^h (\vec x), \,  c_z^h (\vec x, \, \vec x')  \right )  = \set{GP} \left ( \vec \phi(\vec x)^\trans \overline{\vec z}  , \,  \vec \phi(\vec x)^\trans  \vec Q_z^{-1}\vec \phi(\vec x')  \right ) \, .
\end{align}
\end{subequations}
The isocontours of the respective means and standard deviations are plotted in Figures~\ref{fig:plateSourceTrueSource} and~\ref{fig:plateSourceTrueSolution}.  In the following~$\overline{z}^h(\vec x)$ is assumed as the exact solution~$z(\vec x)$. The random source term~$\chi^h(\vec x)$ is used only for obtaining~$\overline{z}^h(\vec x)$ and is not used further. Note that the mesh for computing~$z^h(\vec x)$ is much finer than the mesh for computing the prior probability density for~$u^h(\vec x)$, i.e.\ $n_z \gg n_u$.

\paragraph{Inferred true system response}  In statFEM we aim to infer the unknown true solution~$z(\vec x)$ from a misspecified finite element model yielding the prior probability density for~$u^h(\vec x)$ and a set of observations~$\{\vec x_i, \, \vec y_i \}_{i=1}^{n_o}$. In the present example, we sample the synthetic observations~$\vec y_i$ from the multivariate Gaussian  distribution 
\begin{equation}
	\vec y_i  \sim \set{N}\left( \vec P \overline{\vec z},  \vec P \vec Q_z^{-1} \vec P^\trans + \sigma_e^2 \vec I \right) \,,
\end{equation}
with the variance observation noise chosen as~$\sigma_e^2 = 2.5 \times 10^{-5}$. Our misspecified finite element model is the discretisation of the Poisson-Dirichlet  
\begin{subequations} \label{eq:pwhPriorP}
	\begin{alignat}{2}
		- \Delta u(\vec x) &= f(\vec x)  \quad && \text{in } \Omega \,,
		\\
		u( \vec x ) &= 0  \quad && \text{on } \partial \Omega \,.
	\end{alignat}
\end{subequations}
with the source term~$f(\vec x) =  \overline{f}(\vec x)  +  s(\vec x)$ consisting of the deterministic source 
\begin{equation}
	\overline{f}(\vec x) = - 28 \pi^2 \sin( 3 \pi x_1 ) \sin( 4 \pi x_2) \,,
\end{equation}
and the zero-mean random source satisfying the SPDE 
\begin{subequations}  \label{eq:pwhExactSPDESecond}
	\begin{alignat}{2}
		\left (  800   - \Delta \right ) s(\vec x) &= 4000 \sqrt{2 \pi} \, g (\vec x) \quad && \text{in } \Omega \,,
		\\
		\nabla  s(\vec x) \cdot \vec n  &= 0  \quad && \text{on } \partial \Omega \,.
	\end{alignat}
\end{subequations}
The respective multivariate Gaussian distribution of the finite element solution~$\vec u$ has the form given in~\eqref{eq:forwardPrior}. The Gaussian processes corresponding to the discretised forcing and the finite element solution are given by
\begin{subequations}
\begin{align}
	f^h(\vec x) & \sim \set{GP} \left (\overline{f}^h (\vec x), \,  c_f^h (\vec x, \, \vec x') \right ) = \set{GP} \left (\vec \phi(\vec x)^\trans \, \overline{\vec f}, \,   \vec \phi(\vec x)^\trans  \vec Q_s^{-1}  \vec \phi(\vec x')  \right ) \,  ,  \\
	u^h(\vec x) & \sim \set{GP} \left ( \overline{u}^h (\vec x), \,  c_u^h (\vec x, \, \vec x')  \right )  = \set{GP} \left ( \vec \phi(\vec x)^\trans \vec A^{-1} \overline{\vec f} \, \overline{\vec z}  , \,  \vec \phi(\vec x)^\trans \vec A^{-1} \vec M  \vec Q_s^{-1} \vec M^\trans \vec A^{-T}\vec \phi(\vec x')  \right ) \, .
\end{align}
\end{subequations}
 Figures~\ref{fig:plateSourcePriorSource} and~\ref{fig:plateSourcePriorSolution} show the mean and standard deviation of the finite element prior source~$f^h(\vec x)$ and response~$u^h(\vec x)$, respectively. Although the prior mean~$\overline{u}^h(\vec x)$ in Figure~\ref{fig:plateSourcePriorSolution} broadly resembles the true response mean~$\overline{z}^h(\vec x)$ in Figure~\ref{fig:plateSourceTrueSolution}, they are clear differences closer to the outer boundaries of the domain. The respective standard deviations look similar even though the magnitude of the standard deviation for the true source in Figure~\ref{fig:plateSourceTrueSource} is about four times smaller than the one for the prior source in Figure~\ref{fig:plateSourcePriorSource}. 

\begin{figure}
	\centering
	\subfloat[][{$n_y = 30$, $n_o = 2$}]
	{\includegraphics[width=0.22\textwidth]{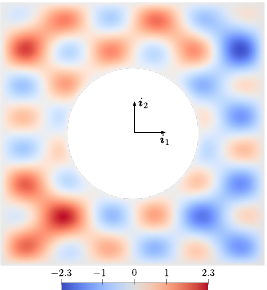}}
	\hfil
	\subfloat[][{$n_y = 60$, $n_o = 2$}]
	{\includegraphics[width=0.22\textwidth]{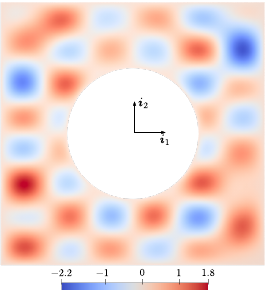}}
	\hfil
	\subfloat[][{$n_y = 120$, $n_o = 2$}]
	{\includegraphics[width=0.22\textwidth]{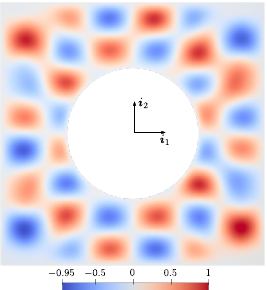}}
	\hfil
	\subfloat[][{$n_y = 432$, $n_o = 2$}]
	{\includegraphics[width=0.22\textwidth]{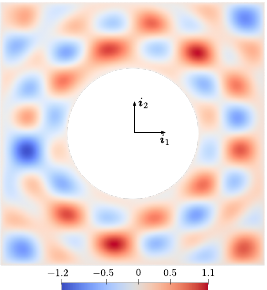}}
	\\
	\subfloat[][{$n_y = 30$, $n_o = 10$}]
	{\includegraphics[width=0.22\textwidth]{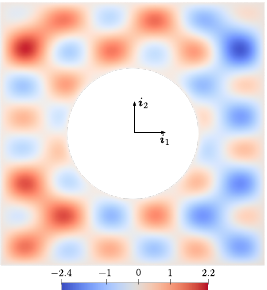}}
	\hfil
	\subfloat[][{$n_y = 60$, $n_o = 10$}]
	{\includegraphics[width=0.22\textwidth]{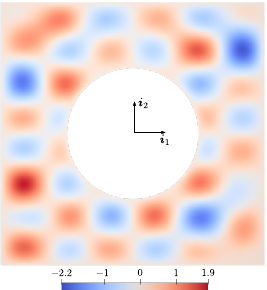}}
	\hfil
	\subfloat[][{$n_y = 120$, $n_o = 10$}]
	{\includegraphics[width=0.22\textwidth]{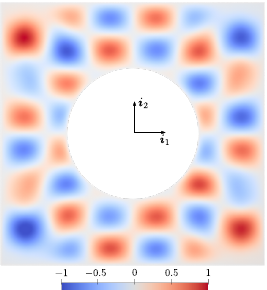}}
	\hfil
	\subfloat[][{$n_y = 432$, $n_o = 10$}]
	{\includegraphics[width=0.22\textwidth]{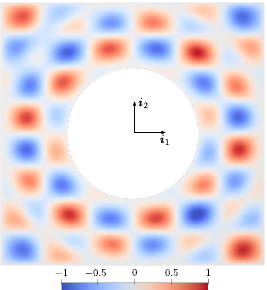}}
	\\
	\subfloat[][{$n_y = 30$, $n_o = 20$}]
	{\includegraphics[width=0.22\textwidth]{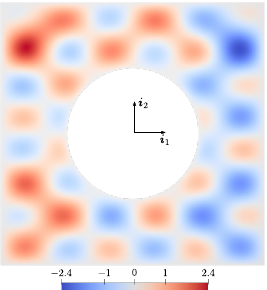}}
	\hfil
	\subfloat[][{$n_y = 60$, $n_o = 20$}]
	{\includegraphics[width=0.22\textwidth]{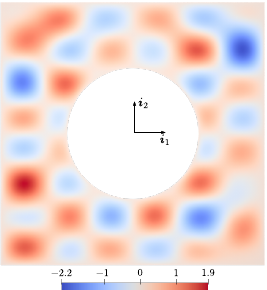}}
	\hfil
	\subfloat[][{$n_y = 120$, $n_o = 20$}]
	{\includegraphics[width=0.22\textwidth]{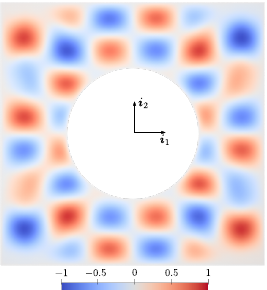}}
	\hfil
	\subfloat[][{$n_y = 432$, $n_o = 20$}]
	{\includegraphics[width=0.22\textwidth]{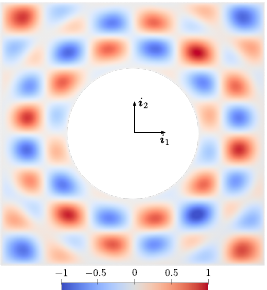}}
	\caption{Poisson-Dirichlet problem with source uncertainty. Mean of the inferred (posterior) true system response~$\overline{z}_{\vert Y}(\vec x)$. Note that all isocontours are plotted with different scales. \label{fig:plateMean}}
\end{figure}

\begin{figure}
	\centering
	\subfloat[][{$n_y = 30$, $n_o = 2$}]
	{\includegraphics[width=0.22\textwidth]{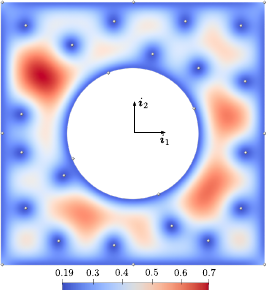}}
	\hfil
	\subfloat[][{$n_y = 60$, $n_o = 2$}]
	{\includegraphics[width=0.22\textwidth]{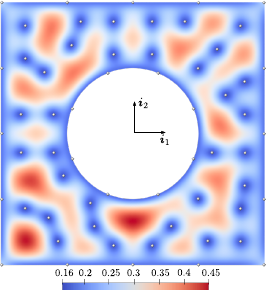}}
	\hfil
	\subfloat[][{$n_y = 120$, $n_o = 2$}]
	{\includegraphics[width=0.22\textwidth]{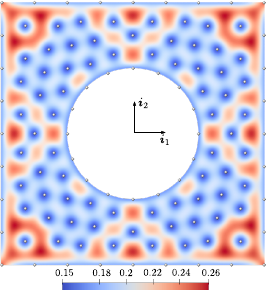}}
	\hfil
	\subfloat[][{$n_y = 432$, $n_o = 2$}]
	{\includegraphics[width=0.22\textwidth]{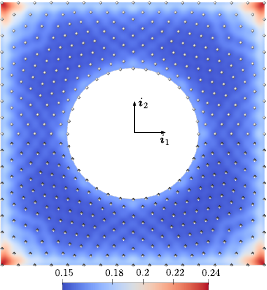}}
	\\
	\subfloat[][{$n_y = 30$, $n_o = 10$}]
	{\includegraphics[width=0.22\textwidth]{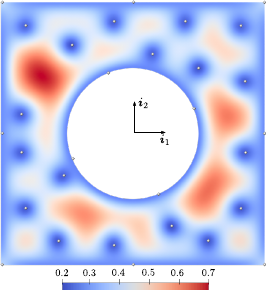}}
	\hfil
	\subfloat[][{$n_y = 60$, $n_o = 10$}]
	{\includegraphics[width=0.22\textwidth]{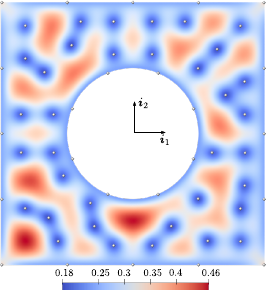}}
	\hfil
	\subfloat[][{$n_y = 120$, $n_o = 10$}]
	{\includegraphics[width=0.22\textwidth]{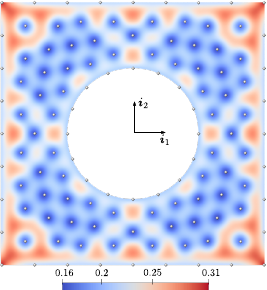}}
	\hfil
	\subfloat[][{$n_y = 432$, $n_o = 10$}]
	{\includegraphics[width=0.22\textwidth]{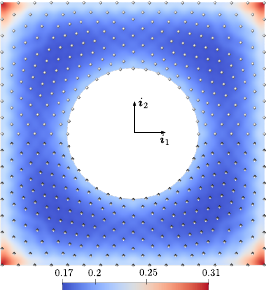}}
	\\
	\subfloat[][{$n_y = 30$, $n_o = 20$}]
	{\includegraphics[width=0.22\textwidth]{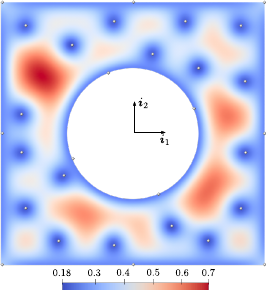}}
	\hfil
	\subfloat[][{$n_y = 60$, $n_o = 20$}]
	{\includegraphics[width=0.22\textwidth]{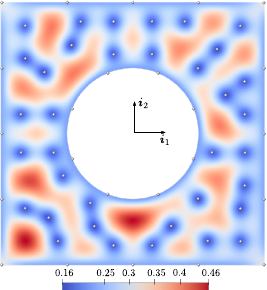}}
	\hfil
	\subfloat[][{$n_y = 120$, $n_o = 20$}]
	{\includegraphics[width=0.22\textwidth]{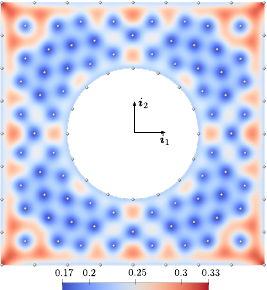}}
	\hfil
	\subfloat[][{$n_y = 432$, $n_o = 20$}]
	{\includegraphics[width=0.22\textwidth]{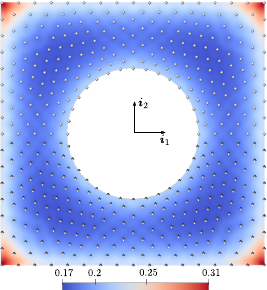}}
	\caption{Poisson-Dirichlet problem with source uncertainty. Standard deviation of the inferred (posterior) true system response~$\sqrt{c^h_{z\vert Y}(\vec x ,\vec x)}$. The grey dots denote the observation points. Note that all isocontours are plotted with different scales. \label{fig:plateStandardDeviation}}
	\centering
	\subfloat[][{True response: $\overline{z}^h(\vec x)$ and $\sqrt{c^h_z(\vec x, \vec x)}$ \label{fig:plateBoundaryTrueSolution}}] {
		\includegraphics[width=0.22\textwidth]{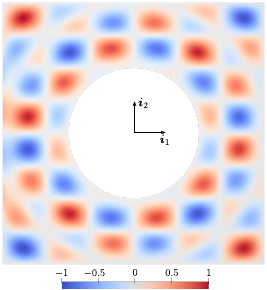} \hspace{0.75em}
		\includegraphics[width=0.22\textwidth]{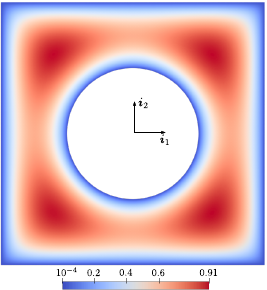}
	}
	\hfil
	\subfloat[][{Prior solution: $\overline{u}^h(\vec x)$ and $\sqrt{c^h_u(\vec x, \vec x)}$ \label{fig:plateBoundaryPriorSolution}}] {
		\includegraphics[width=0.22\textwidth]{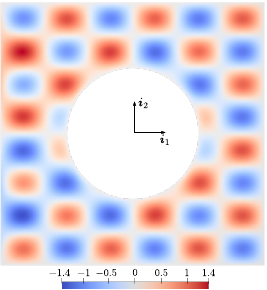} \hspace{0.75em}
		\includegraphics[width=0.22\textwidth]{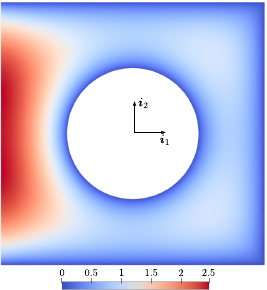}
	}
	\caption{Poisson-Dirichlet problem with source and boundary condition uncertainties. Mean and standard deviation of the true system response~$z^h(\vec x)$ and the finite element prior~$u^h(\vec x)$. \label{fig:plateBoundaryDescription}}
\end{figure}

The statistical finite element posterior distribution~$p(\vec u \vert \vec y)$ is obtained according to~\eqref{eq:statFEpostU} using the sampled data~$\vec y_i$, the finite element prior distribution~$p(\vec u)$ and the discrepancy prior distribution~$p(\vec d)$. Subsequently, the posterior true response distribution~$p(\vec z \vert \vec y)$ is determined according to~\eqref{eq:posteriorTrueResponse}. As detailed in Section~\ref{sec:statFEMobvModel}, the discrepancy prior distribution is derived from the SPDE representation of random fields and has the associated hyperparameters~\mbox{$ \{ \sigma_d, \, \ell_d, \, \nu_d \} $}. For the posterior distributions reported in the following we choose~$\nu_d=1$ and determine remaining two parameters by maximising the log marginal likelihood~$p(\vec y)$ in~\eqref{eq:statFElogMargy}, or the equivalent expression~$p(\vec Y)$ in~\eqref{eq:logMarginalLikelihoodStatFEMultipleGeneral} for repeated readings~$\vec Y \in \mathbb R^{n_y \times n_o}$. Table~\ref{table:plateSource} shows that the optimised hyperparameters~$\{ \sigma_d^*, \, \ell_d^*  \}$ are numerically consistent across~$n_o$ repeated readings for each~$n_y$. Upon obtaining the optimised hyperparameters we compute  the posterior true distribution~$p(\vec z \vert \vec Y)$ and the respective Gaussian process 
\begin{equation}
	z_{\vert Y} (\vec x) \sim \set{GP} \left (\overline{z}_{\vert Y }^h, \,  c_{z\vert Y}^h (\vec x, \, \vec x') \right ) = \set{GP} \left (\vec \phi(\vec x)^\trans \overline{\vec z}_{\vert Y}, \,   \vec \phi(\vec x)^\trans  \vec Q_{z\vert Y}^{-1}  \vec \phi(\vec x')  \right ) \, .
\end{equation}
Note that~$p(\vec z \vert \vec Y)$ has the same form as~$p(\vec z\vert \vec y)$. Its dependence on the number of repeated readings~$n_o$ is only through the log marginal likelihood. Figures~\ref{fig:plateMean} and~\ref{fig:plateStandardDeviation} show the mean and standard deviation of the inferred true system response~$z^h_{\vert Y}(\vec x)$, respectively. When both $n_y$ and~$n_o$ increase, the conditional mean~$\overline{z}^h_{\vert Y}(\vec x)$ converges to the true system response mean~$\overline{z}^h(\vec x)$. The standard deviation of the inferred true system response~$\sqrt{c^h_{z\vert Y}(\vec x ,\vec x)}$ decreases when there are more observations. Overall,~$\sqrt{c^h_{z\vert Y}(\vec x ,\vec x)}$ is smaller than~$\sqrt{c^h_{z}(\vec x ,\vec x)}$, possibly due to problems in approximating the true discrepancy with the assumed discrepancy.

\begin{table}
	\centering
	\caption{Poisson-Dirichlet problem with source uncertainty. Point estimates for the hyperparameters obtained by maximising the marginal likelihood. \label{table:plateSource}}
	\begin{tabular}{c c c c c c c c c}
		\toprule
		$n_y$ & $n_o$ &  $\sigma_d^*$ & $\ell_d^*$ & & $n_y$ & $n_o$ & $\sigma_d^*$ & $\ell_d^*$
		\\ \midrule
		\multirow{3}{*}{$30$} & $2$   & $0.1568$ & $0.1360$ & & \multirow{3}{*}{$120$} & $2$ & $0.1241$ & $0.1310$  \\ 
		& $10$  & $0.1915$ & $0.1670$ &  & & $10$ & $0.1514$ & $0.1569$ \\
		& $20$ & $0.1756$ & $0.1712$ &  & & $20$ & $0.1632$ & $0.1777$ 
		\\ \midrule
		\multirow{3}{*}{$60$} & $2$   & $0.1269$ & $0.1225$ & & \multirow{3}{*}{$432$} & $2$ & $0.1211$ & $0.2391$ \\
		& $10$  & $0.1675$ & $0.1616$ &  & & $10$ & $0.1546$ & $0.2436$\\
		& $20$ & $0.1589$ & $0.1532$ & & & $20$ & $0.1568$ & $0.2470$
		\\ \bottomrule
	\end{tabular}
\end{table}
%
\subsubsection{Source and boundary condition uncertainties \label{sec:statFETwoDBoundary}}
%
In this slightly modified example, the source terms are the same as in the previous section. The only change concerns the Dirichlet boundary conditions of the true system response and the finite element prior. To this end, we split the boundary~$\partial \Omega$ as depicted in Figure~\ref{fig:plateProblemDescription} into two non-overlapping sets 
\begin{equation}
	\partial \Omega_1 = \{ \vec x \vert x_1 = -1, -1 < x_2 < 1 \} \, , \quad  \quad \partial \Omega_2 = \partial \Omega \setminus \partial \Omega_1 \, .
\end{equation}

The true system response~$z(\vec x)$ is chosen as the solution of the Poisson-Dirichlet problem 
\begin{subequations}
	\begin{alignat}{2}
		- \Delta z(\vec x) &= \overline{\chi} (\vec x) + \chi_s(\vec x)  \quad  && \text{in } \Omega \,,  \\
		z(\vec x) &=  \frac{1}{2} \cos \left (\frac{1}{2} \pi x_2 \right )  \sin \left (3 \pi x_2 \right ) + \chi_s(\vec x)  \quad  && \text{on }  \partial \Omega_1 \, , \\
		z(\vec x)  &=  0   \quad  && \text{on }  \partial \Omega_2 \, , 
	\end{alignat}
\end{subequations}
which is with the exception of the boundary condition on~$\partial \Omega_1$ the same as~\eqref{eq:pwhExactPDE}. The resulting  true system response~$z(\vec x)$ depicted in   Figure~\ref{fig:plateBoundaryTrueSolution} is visually very similar to Figure~\ref{fig:plateSourceTrueSolution},  with slight differences in the mean~$\overline{z}^h(\vec x)$ close to the boundary~$\partial \Omega_1$. 

Similar to~\eqref{eq:pwhPriorP} the misspecified finite element prior is the solution of the Poisson-Dirichlet problem 
\begin{subequations}
	\begin{alignat}{2}
		- \Delta u(\vec x) &= \overline{f}(\vec x) + s (\vec x)  \quad && \text{in } \Omega \,,  \\ 
		u (\vec x) &= \frac{1}{2} \sin \left (\pi x_2 \right )  + s (\vec x) \quad &&  \text{on } \partial \Omega_1 \, , \\ 
		u (\vec x) &= 0 \quad  &&  \text{on } \partial \Omega_2 \, .
	\end{alignat}
\end{subequations}
To take into account that the boundary condition~$\partial \Omega_1$ is now uncertain, the computation of the prior~$p(\vec u)$ must be slightly modified. Specifically,  in computing the precision matrix~$\vec Q_u$ according to~\eqref{eq:forwardPrior}  we assume that the boundary~$\partial \Omega_1$ is a Neumann boundary by not deleting the respective rows and columns of the system matrix~$\vec A$. The so-obtained mean and standard deviation of the finite element prior are depicted  Figure~\ref{fig:plateBoundaryPriorSolution}. Notice that the standard deviation along the boundary~$\partial \Omega_1$ is non-zero as desired. Without such a modification Bayesian inversion is ill-posed considering that the mean of the true and finite element solutions~$\overline z(\vec x)$ and~$\overline u^h(\vec x)$ on~$\partial \Omega_1$ are different so that it is impossible to be certain about both values. In passing we note that the boundary~$\partial \Omega_1$ can be chosen as a Robin boundary to better reflect our confidence in the prior~$u^h(\vec x)$. 

\begin{figure}
	\centering
	\subfloat[][{$n_y = 30$, $n_o = 2$}]
	{\includegraphics[width=0.25\textwidth]{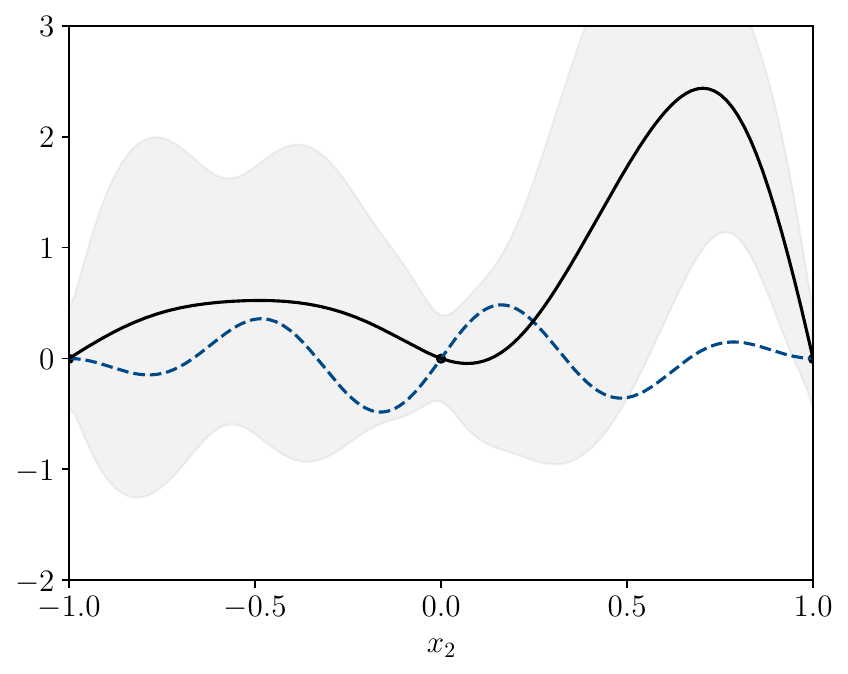}}
	\hfil
	\subfloat[][{$n_y = 60$, $n_o = 2$}]
	{\includegraphics[width=0.25\textwidth]{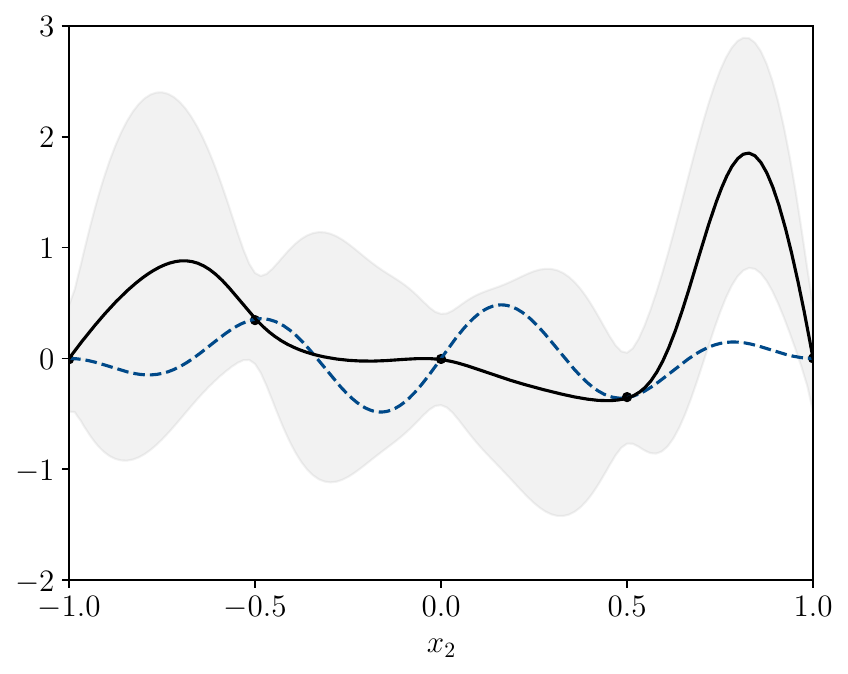}}
	\hfil
	\subfloat[][{$n_y = 120$, $n_o = 2$}]
	{\includegraphics[width=0.25\textwidth]{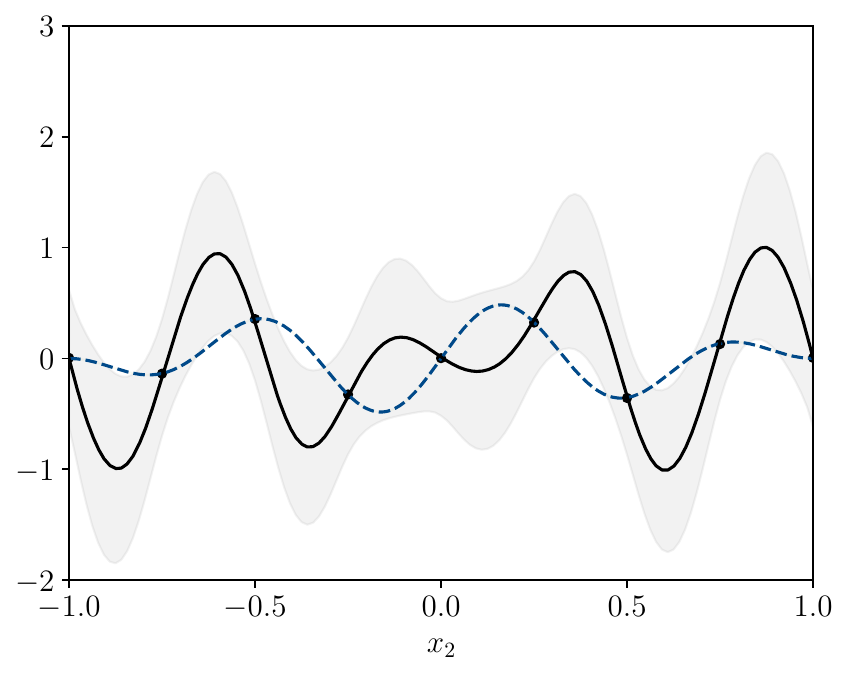}}
	\hfil
	\subfloat[][{$n_y = 432$, $n_o = 2$}]
	{\includegraphics[width=0.25\textwidth]{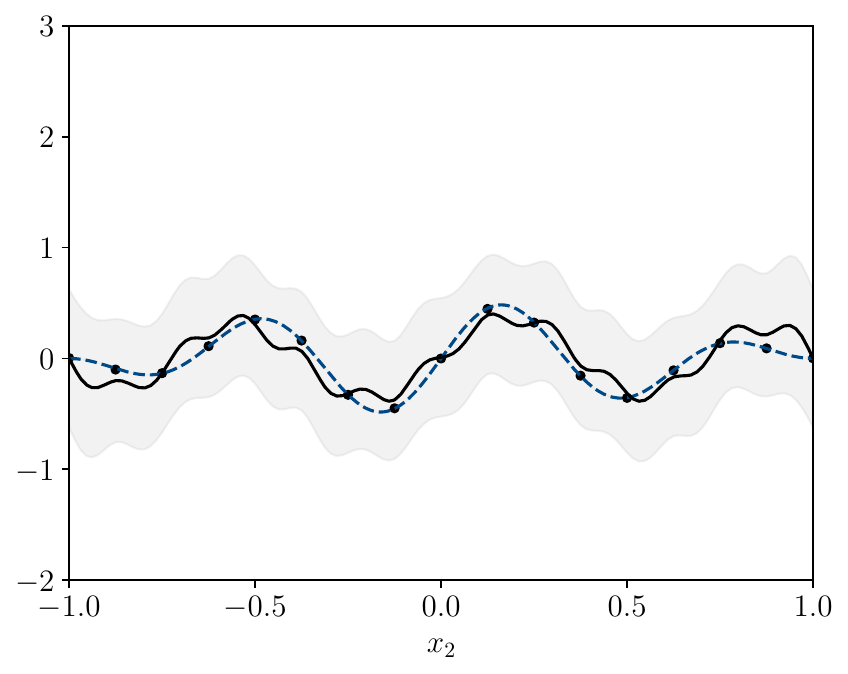}}
	\\
	\subfloat[][{$n_y = 30$, $n_o = 10$}]
	{\includegraphics[width=0.25\textwidth]{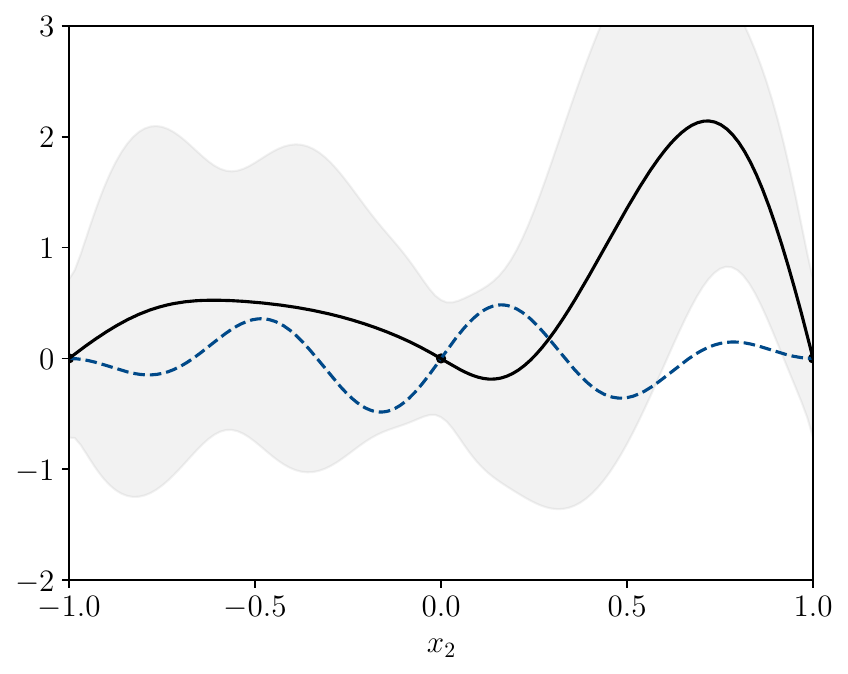}}
	\hfil
	\subfloat[][{$n_y = 60$, $n_o = 10$}]
	{\includegraphics[width=0.25\textwidth]{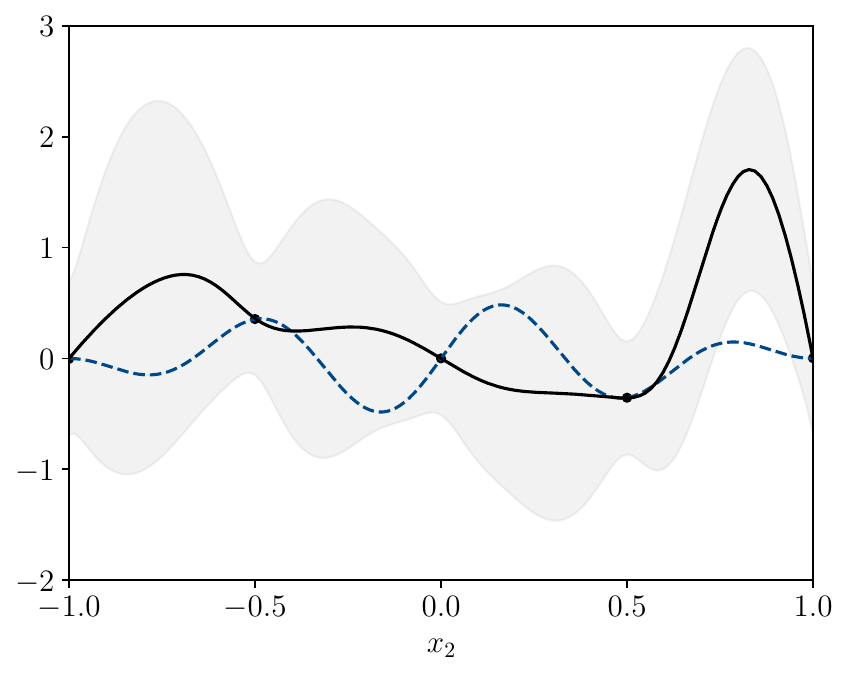}}
	\hfil
	\subfloat[][{$n_y = 120$, $n_o = 10$}]
	{\includegraphics[width=0.25\textwidth]{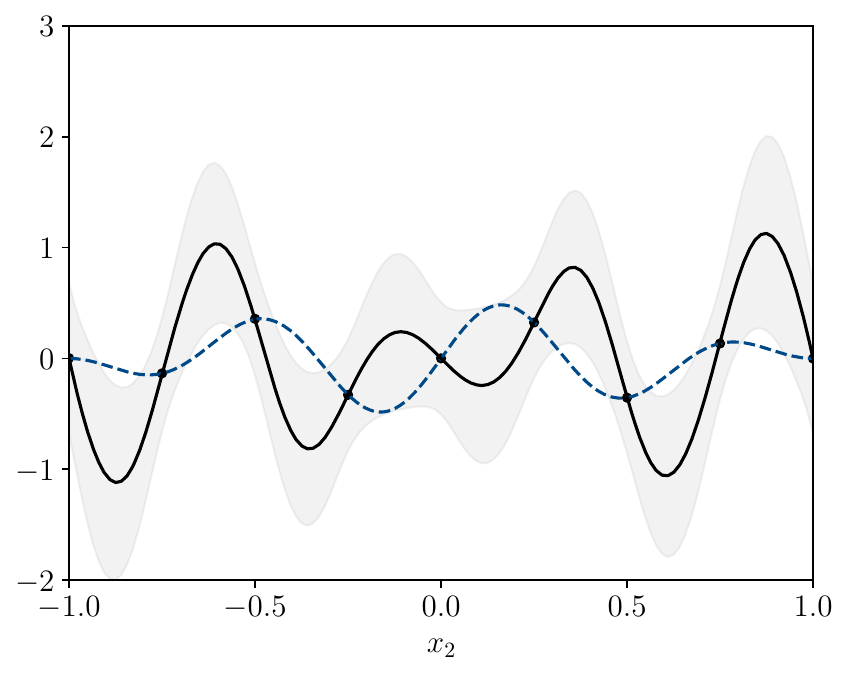}}
	\hfil
	\subfloat[][{$n_y = 432$, $n_o = 10$}]
	{\includegraphics[width=0.25\textwidth]{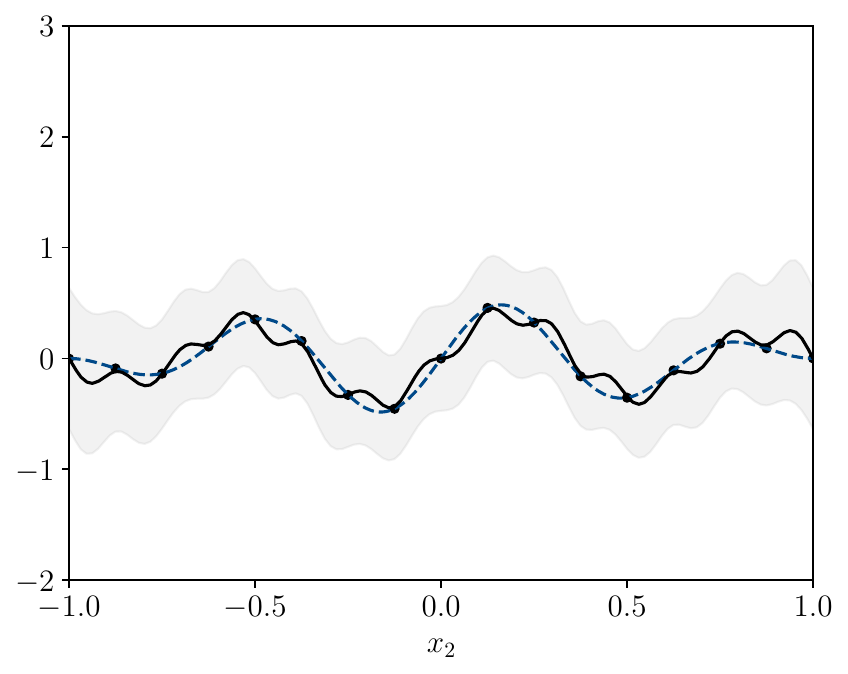}}
	\\
	\subfloat[][{$n_y = 30$, $n_o = 20$}]
	{\includegraphics[width=0.25\textwidth]{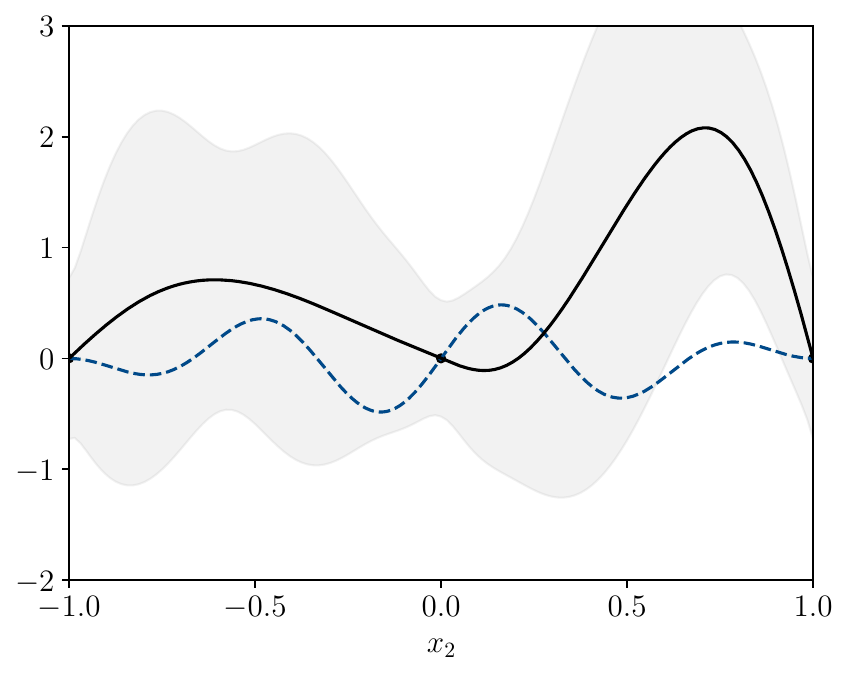}}
	\hfil
	\subfloat[][{$n_y = 60$, $n_o = 20$}]
	{\includegraphics[width=0.25\textwidth]{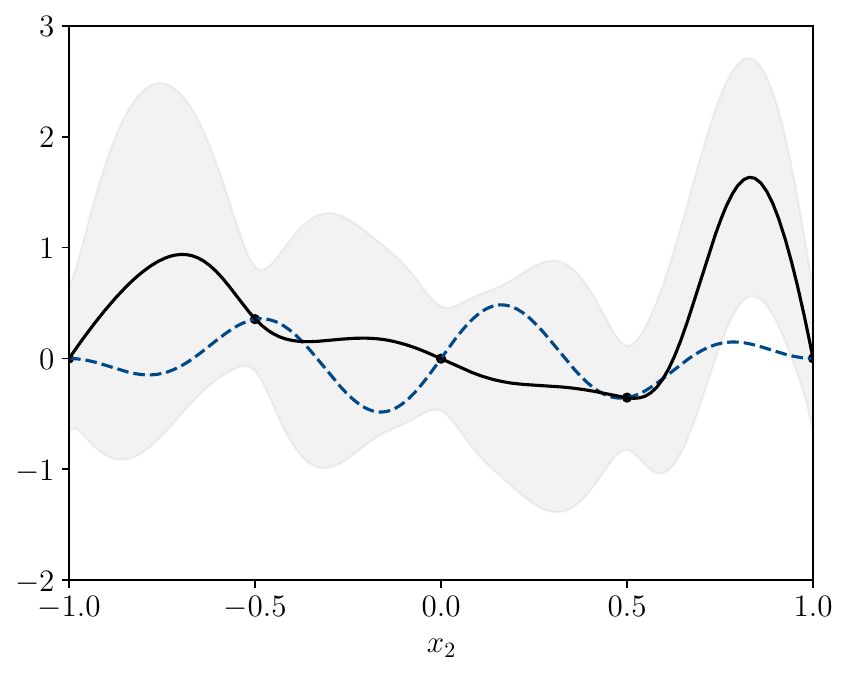}}
	\hfil
	\subfloat[][{$n_y = 120$, $n_o = 20$}]
	{\includegraphics[width=0.25\textwidth]{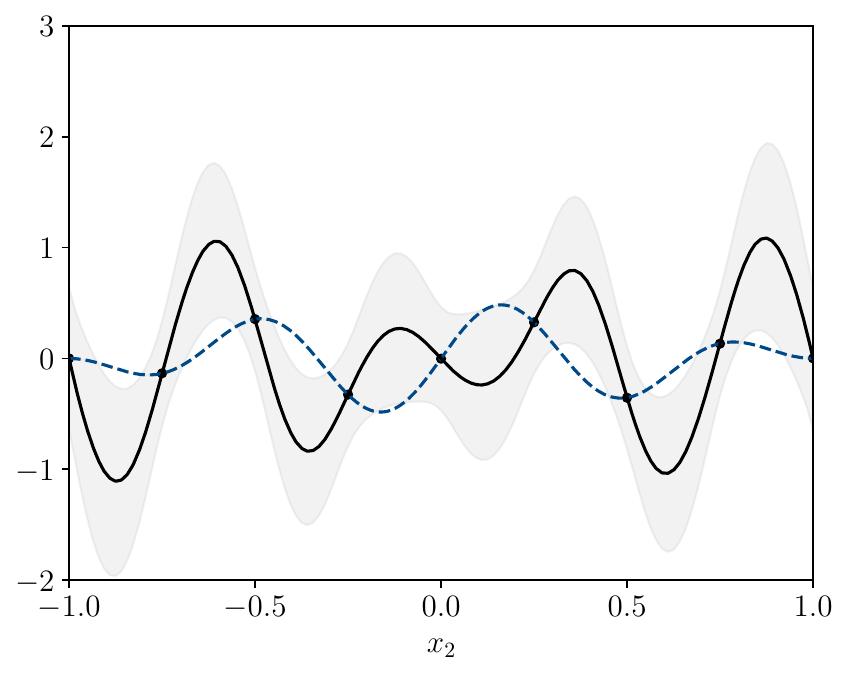}}
	\hfil
	\subfloat[][{$n_y = 432$, $n_o = 20$}]
	{\includegraphics[width=0.25\textwidth]{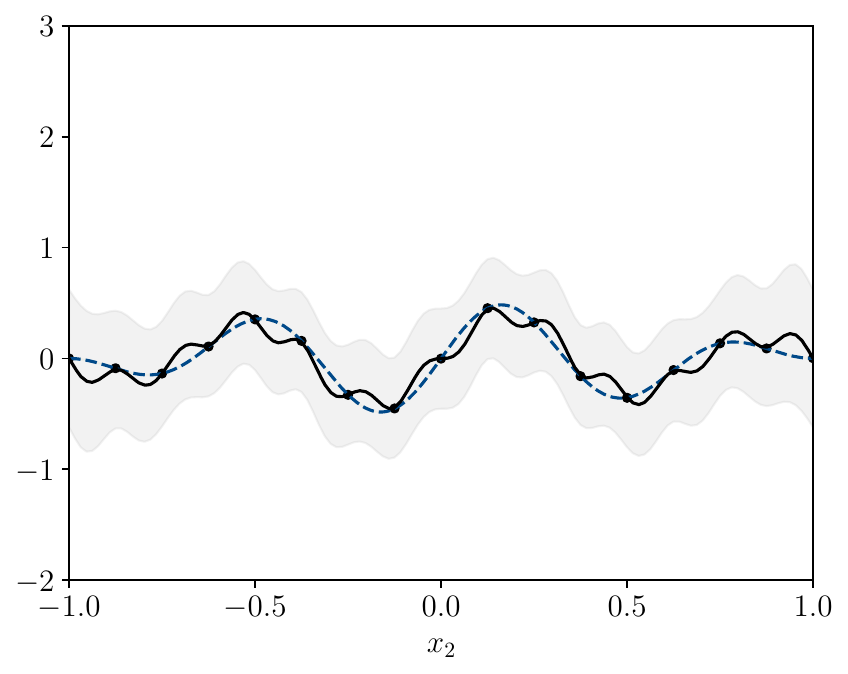}}
	\caption{Poisson-Dirichlet problem with source and boundary condition uncertainties. Convergence of inferred posterior true system response~$z_{\vert Y}^h(\vec x)$ along the boundary~$\partial \Omega_1$ (horizontal axis).  The dashed blue lines denote the mean of the true system response~$\overline{z}^h(\vec x)$, the black lines the inferred posterior mean~$\overline{z}^h _{\vert Y}(\vec x)$ and the shaded areas the~$95\%$ confidence regions corresponding to the posterior variance~$c^h_{\vert Y} (\vec x, \, \vec x')$. The black dots represent the empirical mean of the observation data~$ \sum_{i=1}^{n_o} \vec y_i / n_o$. \label{fig:plateBoundaryInfer}}
\end{figure}

We choose~$\nu_d=1$ and determine remaining two hyperparameters~$\{ \sigma_d, \, \ell_d  \}$ of the SPDE for the discrepancy by maximising the log marginal likelihood~$p(\vec y)$ in~\eqref{eq:statFElogMargy} or the equivalent expression~$p(\vec Y)$ in~\eqref{eq:logMarginalLikelihoodStatFEMultiple}.  Furthermore, Figure~\ref{fig:plateBoundaryInfer} confirms that when~$n_y$ increases, on the random boundary~$\partial \Omega_1$ the inferred mean~$\overline{z}_{\vert Y}(\vec x)$ converges to the underlying true response mean~$\overline{z}(\vec x)$ with increasing confidence. 
%
%
%
\subsection{Statistical finite element analysis of a Kirchhoff--Love thin shell}
\begin{figure} 
	\centering
	\includegraphics[width=0.275\textwidth]{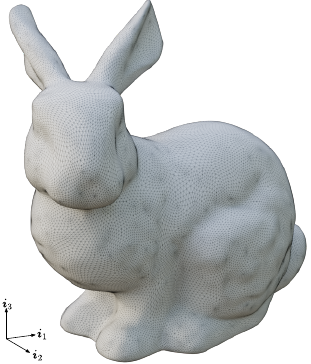}
	\hfil
	\includegraphics[width=0.275\textwidth]{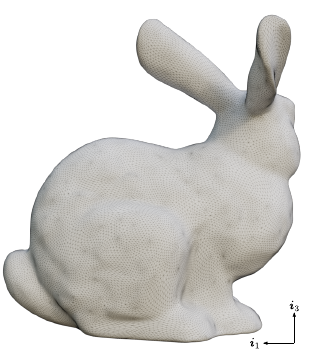}
	\caption{Statistical finite element analysis of a Kirchhoff--Love thin shell. Stanford bunny discretised with a triangular control mesh.}
	\label{fig:thinShellDescription}
\end{figure}

%
Our last example concerns the statistical finite element analysis of a Kirchhoff-Love thin shell, namely the Stanford bunny.  Different from the scalar problems discussed so far, the unknown vector-valued solution field for the shell equations consists of the three displacement components of the mid-surface. Specifically, we consider the linear elastic analysis of the Stanford bunny subjected to a uniform area load of $(0, \, 0, \, -1)^\trans$. The geometry is the same as the one introduced in Sections~\ref{sec:introductionRelated} and~\ref{sec:maternExtensions} and is parametrised with Loop subdivision surfaces. The computational control mesh shown in Figure~\ref{fig:thinShellDescription} consists of 79488 triangles and 39826 vertices. We discretise the thin shell equations with the respective Loop basis functions, see~\cite{Cirak:2000aa,Cirak:2001aa,Cirak:2002aa,Cirak:2011aa} for details. The displacements of the vertices at the bottom of the bunny are prescribed to be zero. The true shell thickness is $0.085$, the Young's modulus is~$10^5$ and the Poisson ratio is~$0.3$. In the finite element model for determining the prior density we assume that the shell thickness is misspecified and has the value~$0.1$. 

\begin{figure} 
	\centering
	\subfloat[][Prior mean \label{fig:thinShellPriorA}]
	{
		\includegraphics[width=0.275\textwidth]{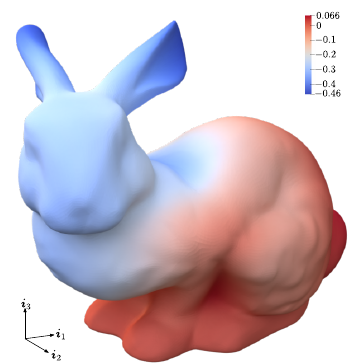}
	}
	\hfil
	\subfloat[][True mean  \label{fig:thinShellPriorB}]
	{
		\includegraphics[width=0.275\textwidth]{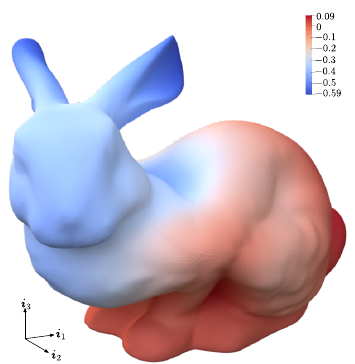}
	}
	\hfil
	\subfloat[][Prior mean in blue and true mean in transparent white \label{fig:thinShellPriorC}]
	{
		\includegraphics[width=0.275\textwidth]{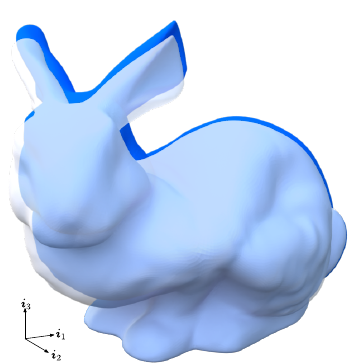}
	}
	\caption{Statistical finite element analysis of a Kirchhoff--Love thin shell. Prior and true mean deflections. In (a) and (b)  the isocontours represent the mean~$\vec i_3$-displacement of the prior and true solutions, respectively. In (c), we overlay the mean of both deflected shapes to aid visual comparison.}
	\label{fig:thinShellPriorMean}
\end{figure}

\paragraph{Finite element prior response} Similar to the Poisson-Dirichlet problems discussed so far, the discretisation of the thin-shell equation leads to  
\begin{equation} \label{eq:exampleKLfeP}
	\vec A \vec u = \vec M \left(\overline {\vec f} + \vec s \right) \, ,
\end{equation}
where~$\vec A$ is the stiffness matrix,~$\vec u$ is the displacement vector,~$\vec M$ is the mass matrix, $\overline{\vec f}$ is the deterministic external force vector and $\vec s$ is the random external force vector. We assume for the sake of discussion that the degrees of freedom are enumerated such that the displacement and force vectors have the following structure 
\begin{equation}
	\vec u = \begin{pmatrix} \vec u_x \\ \vec u_y \\ \vec u_z \end{pmatrix}  \, , \quad \overline{\vec f} = \begin{pmatrix} \overline{\vec f}_x \\ \overline{ \vec f}_y \\ \overline{ \vec f}_z \end{pmatrix}   \,, \quad \vec s = \begin{pmatrix} \vec s_x \\ \vec s_y \\ \vec s_z \end{pmatrix}   \, .
\end{equation} 
That is, the components of nodal coefficients in the three coordinate directions are ordered consecutively. Hence, given that the prescribed uniform area load is $(0, \, 0,  \, -1)^\trans$,  the external deterministic force vector is~\mbox{$\overline{\vec f} = (  \vec 0, \, \vec 0, \, - \vec 1 )^\trans$}, where~$\vec 0$ and~$\vec 1$ are all-zero and all-one vectors. We determine the probability density of the random external force vector~$\vec s$ using the SPDE representation of random fields and assuming that the components~$\vec s_x$, $\vec s_y$ and~$\vec s_z$ are uncorrelated such that 
\begin{equation} \label{eq:shellExS}
	 \vec s =  \begin{pmatrix} \vec s_x \\ \vec s_y \\ \vec s_z  \end{pmatrix} \sim  \set N \left(\vec 0, \, \vec Q_s^{-1}\right) = \set N \left  ( \begin{pmatrix} \vec 0  \\ \vec 0 \\ \vec 0 \end{pmatrix}, \, \begin{pmatrix} \vec Q_{sx}^{-1} & & \\ & \vec Q_{sy}^{-1} & \\ & & \vec Q_{sz}^{-1} \end{pmatrix} \right )  \, . 	
\end{equation}
We solve three independent SPDE problems to obtain the precision matrices~$\vec Q_{sx}$, $\vec Q_{sy}$ and $\vec Q_{sz}$ with each having a possibly different standard deviation, length-scale and smoothness parameter. In the present example, we choose the parameters for the three SPDE problems as 
\begin{equation}
	 \vec \sigma_s = ( 0.025 , \, 0.6 , \,  0.025)^\trans \, , \quad \vec \ell_s = ( 2.5 , \,  2.5 , \,  2.5 )^\trans  \, ,  \quad \vec \nu_s = ( 1 , \,  1 , \, 1 )^\trans \, .
\end{equation}
With the so-obtained precision matrices the finite element prior has the density 
\begin{equation}
	p( \vec u ) = \set N \left( \overline{\vec u}, \,  \vec Q_u^{-1}\right)  =  \set N \left( {\vec A}^{-1} \overline {\vec f} , \, \vec P    {\vec A}^{-1}  \vec M  \vec Q_{s}^{-1}  \vec M^\trans  {\vec A} ^{-1}  \vec P^\trans   \right) \, . 
\end{equation}
The mean deflection~$\overline{\vec u}$ is depicted in Figure~\ref{fig:thinShellPriorA}.

\paragraph{True system response}  For generating the used synthetic observation data~$\vec y$ we consider a true random solution vector~$\vec z$ corresponding to a shell with thickness~$0.085$ and uniform area load of $(0, \, 0, -\, 1)^\trans$.  While the area load is the same as the one for the finite element prior, the thickness is different, inducing an intrinsic model inadequacy error. The random solution vector~$\vec z$ is the solution of the discretised problem
\begin{equation}
	\widetilde {\vec A} \vec z = \vec M \left(\overline {\vec f} + \vec \chi \right) \,,
\end{equation}
where~$\widetilde {\vec A}$ is the stiffness matrix for the shell with the true thickness~$0.085$, $\overline {\vec f}$  is the deterministic external force vector as in~\eqref{eq:exampleKLfeP}. We choose a true system with a random solution  by introducing the random forcing vector $\vec \chi$ with the density 
\begin{equation}
	 p( \vec \chi)   =   \set N \left(\vec 0, \, \vec Q_\chi^{-1}\right)   \, . 	
\end{equation}
The precision matrix $\vec Q_\chi$ is block-diagonal and is obtained by solving three independent SPDE problems with the parameters \mbox{$\vec \sigma_{\chi} = (0.05, \, 0.05, \,   1)^\trans$}, \mbox{$\vec \ell_\chi = ( 1.5 , \,1.5 , \, 1.5)^\trans$} and \mbox{$\vec \nu_{\chi} = ( 1 , \,  1 , \, 1 )^\trans$}. With the so-obtained precision matrix the probability density for the true solution is given by
\begin{equation} \label{eq:shellExTrueSol}
	p(\vec z) = \set N \left(\overline {\vec z}, \, \vec Q^{-1}_z \right) = \set N \left (  \widetilde {\vec A}^{-1} \overline {\vec f} , \,   \widetilde {\vec A}^{-1}  \vec M  \vec Q_{\chi}^{-1}  \vec M^\trans  \widetilde {\vec A} ^{-\trans} \right ) \, . 
\end{equation}
The true mean deflection~$\overline{\vec z}$ is depicted in Figure~\ref{fig:thinShellPriorB}. And, the true mean deflection~$\overline{\vec z}$ and finite element prior mean deflection~$\overline{\vec u}$ are compared in Figure~\ref{fig:thinShellPriorC}. The difference between the two is due to misspecification of the shell thickness.

\begin{figure} 
	\centering
	\subfloat[][$n_y = 333$ \label{fig:thinShellData111}]
	{
		\includegraphics[width=0.275\textwidth]{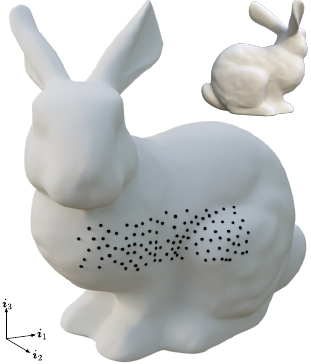}
	}
	\hfil
	\subfloat[][$n_y = 564$  \label{fig:thinShellData188}]
	{
		\includegraphics[width=0.275\textwidth]{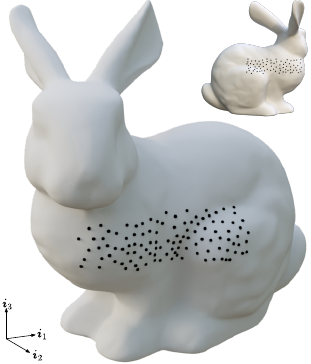}
	}
	\hfil
	\subfloat[][$n_y = 765$ \label{fig:thinShellData255}]
	{
		\includegraphics[width=0.275\textwidth]{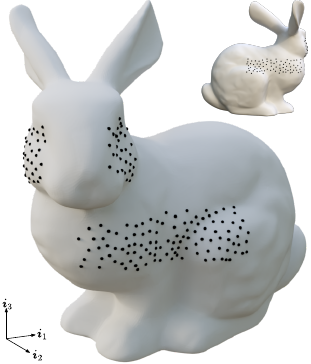}
	}
	\caption{Statistical finite element analysis of a Kirchhoff--Love thin shell. Observation locations.}
	\label{fig:thinShellDataPoints}
\end{figure}

\paragraph{Observation data.} As depicted in Figure~\ref{fig:thinShellDataPoints}, we consider three sets of observations with~\mbox{$n_y \in \{ 333, \, 564, \, 765\}$}.  There are in total~$111$,  $188$ and $255$ observation locations and at each location all the three components  of the displacement vector are observed. Taking into account the true solution~\eqref{eq:shellExTrueSol}, we sample the synthetic observations from 
\begin{equation}
	\vec y_i  \sim \set N \left( \vec P \overline{\vec z} , \, \vec P  \vec Q_z^{-1}  \vec P^\trans + \sigma_e^2 \vec I  \right) \, , 
\end{equation}
where~$\vec P$ is a suitably chosen observation matrix depending on the number and location of observations and the observation error is~$\sigma_e = 0.005$. 

\paragraph{Inferred true system response} Given that the solution field is vector-valued we choose a discrepancy which is  a vector-valued random field, modelled using the SPDE representation, such that
\begin{equation} \label{eq:shellExD}
	 \vec d =  \begin{pmatrix} \vec d_x \\ \vec d_y \\ \vec d_z  \end{pmatrix} \sim  \set N \left(\vec 0, \, \vec Q_d^{-1}\right) = \set N \left  ( \begin{pmatrix} \vec 0  \\ \vec 0 \\ \vec 0 \end{pmatrix}, \, \begin{pmatrix} \vec Q_{dx}^{-1} & & \\ & \vec Q_{dy}^{-1} & \\ & & \vec Q_{dz}^{-1} \end{pmatrix} \right )  \, . 	
\end{equation}
The components~$\vec d_x$, $\vec d_y$ and $\vec d_z$ are uncorrelated and are modelled by three independent SPDEs. For a given observation data~$\vec y$ the parameters of the SPDEs can be obtained by maximising the marginal likelihood. We found by empirical grid search that the parameters  \mbox{$\vec \sigma_{d} = ( 0.2, \, 0.3 , \,  0.2 )^\trans$},  \mbox{$\vec \ell_d = ( 3.5 \,,  3.5 , \, 3.5 )^\trans$} and \mbox{$\vec \nu_d = ( 1 , \,   1 , \,  1 )^\trans$} yield an optimal marginal likelihood for the case with~$n_y=333$ observations. For simplicity, we use the same parameters for the considered three different sets of observations. 

\begin{figure} 
	\centering
	\subfloat[][$n_y = 333$  \label{fig:thinShellPosterior111}]
	{
		\includegraphics[width=0.275\textwidth]{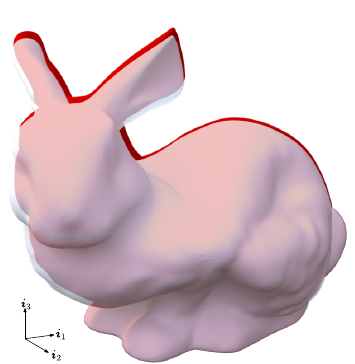}
	}
	\hfil
	\subfloat[][$n_y = 564$   \label{fig:thinShellPosterior188}]
	{
		\includegraphics[width=0.275\textwidth]{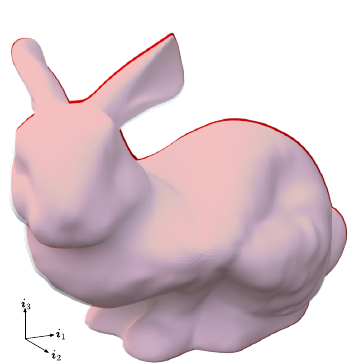}
	}
	\hfil
	\subfloat[][$n_y = 765$  \label{fig:thinShellPosterior255}]
	{
		\includegraphics[width=0.275\textwidth]{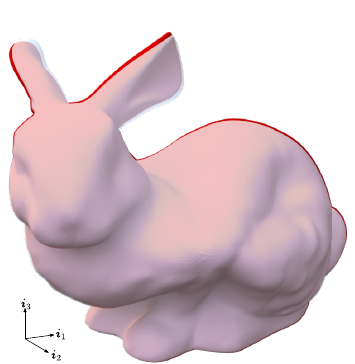}
	}
	\caption{Statistical finite element analysis of a Kirchhoff--Love thin shell. Inferred posterior (red) and true (transparent white) mean deflections.}
	\label{fig:thinShellPosteriorMean}
\end{figure}

\begin{figure}
	\centering
	\includegraphics[width=0.275\textwidth]{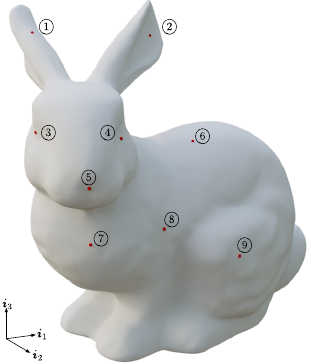}
	\hfil
	\includegraphics[width=0.275\textwidth]{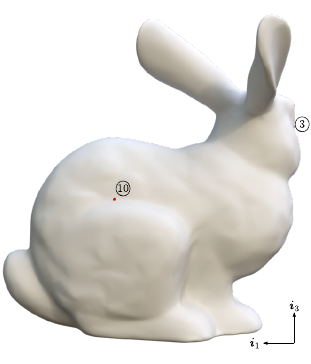}
	\caption{Statistical finite element analysis of a Kirchhoff--Love thin shell. Locations and labels of selected evaluation points.
		\label{fig:thinShellTestPoints}}
\end{figure}

The true and statFEM inferred deflected shapes are compared in Figure~\ref{fig:thinShellPosteriorMean}.  Notice the close agreement between the inferred and true deflected shapes with an increasing number of observation points. The  agreement between the true and inferred shapes is to be contrasted with the relatively poor agreement between the true and finite element deflected shapes in Figure~\ref{fig:thinShellPriorC}. As evident a relatively low number of observation points restricted to the body of the bunny yield a shape that is significantly closer to the true shape. In particular,  the close agreement between the true and inferred shapes close to the head of the bunny in Figure~\ref{fig:thinShellPosteriorMean}  is noteworthy, given that all the sampling points are on the body of the bunny.
Furthermore, we compare the true, finite element and the inferred displacements at the ten selected points depicted in Figure~\ref{fig:thinShellTestPoints}. As given in Table~\ref{table:thinShellDisplacements}, at all the ten points the true displacements are much closer to the inferred displacements than the finite element displacements. Moreover, we report that the empirical mean difference between the true and inferred deflected shapes are~$0.0411$,~$0.0167$ and~$0.0145$ for~$n_y = 333$,~$n_y = 564$ and~$n_y = 765$, respectively. The inferred displacements converge towards the true displacements with an increase of the number of observation points. Finally, the wall-clock times for computing the posterior mean vector for $n_y=333$, $n_y=564$, and $n_y = 765$ are $135.3 s$, $170.5 s$ and $203. 5 s$, respectively, on a MacBook Pro with an M1 Pro chip and 32 GB RAM using a single core. The increase in the computing time is due to the decrease in the sparsity of the precision matrix, as discussed in Section~\ref{sec:statFEMPosteriorDensities}.

\begin{table}
	\centering
	\caption{Statistical finite element analysis of a Kirchhoff--Love thin shell. Distance between the unknown true mean~$\overline{\vec z}$ and the inferred true posterior mean~$\overline{\vec z}_{\vert y}$ and distance between~$\overline{\vec z}$ and the finite element prior mean~$\overline{\vec u}$ at the evaluation points in Figure~\ref{fig:thinShellTestPoints}.  \label{table:thinShellDisplacements}}
	\begin{tabular}{c c c c c c c c c}
		\toprule
		\multirow{2}{*}{Label} &  \multirow{2}{*}{$\text{dist}\left(\overline{\vec u}, \overline{\vec z} \right)$ } &  \multicolumn{3}{c}{$\text{dist}\left(\overline{\vec z}_{\vert y}, \overline{\vec z}\right)$} \\ \cmidrule{3-5}
		& &  $n_y =111 $ & $n_y = 188 $ & $n_y = 255 $ 
		\\ \midrule
		\large \textcircled{\footnotesize $1$} & $0.1172$ & $0.0639$ &  $0.0256$ &  $0.0406$
		\\
		\large \textcircled{\footnotesize $2$} & $0.1406$ & $0.0781$  &  $0.0305$ & $0.0396$
		\\
		\large \textcircled{\footnotesize $3$} &  $0.1305$ &  $0.0738$ & $0.0224$ & $0.0127$
		\\
		\large \textcircled{\footnotesize $4$} & $0.1195$ & $0.0689$  & $0.0281$  & $0.0230$
		\\
		\large \textcircled{\footnotesize $5$} & $0.1347$ & $0.0778$  &  $0.0282$ & $0.0148$
		\\
		\large \textcircled{\footnotesize $6$} & $0.1152$ & $0.0685$ & $0.0310$ & $0.0250$
		\\
		\large \textcircled{\footnotesize $7$} & $0.0993$ & $0.0571$ & $0.0201$  & $0.0106$
		\\
		\large \textcircled{\footnotesize $8$} & $0.1671$ & $0.1044$ & $0.0534$ & $0.0448$
		\\
		\large \textcircled{\footnotesize $9$} &  $0.0836$  &  $0.0531$ &  $0.0295$  & $0.0255$
		\\
		\large \textcircled{\footnotesize $10$} & $0.0944$ & $0.0630$ & $0.0363$ & $0.0301$
		\\ \bottomrule
	\end{tabular}
\end{table}

%
\section{Conclusions \label{sec:conclusions}}
%
We introduced an approach for large-scale Gaussian process regression and statistical finite element analysis on Euclidean domains and shells discretised by finite elements. According to the SPDE representation of Gaussian random fields, the precision matrix of a Mat\'ern field is equal to the finite element system matrix of the respective SPDE. The order of the SPDE depends on the smoothness of the random field and can be fractional. When its order is finite and non-fractional, the system matrix of the SPDE is sparse even though its inverse, the covariance matrix, is dense. The sparse precision matrix formulation makes it possible to efficiently solve large-scale Gaussian process regression and statistical finite element analysis problems using only sparse matrix operations. This is in stark contrast to the prevalent covariance-based GP regression and statFEM approaches requiring onerous dense matrix operations. Hence, the SPDE approach is suitable for large-scale problems with hundreds of thousands of unknowns, whereas common covariance-based techniques only apply to relatively low-dimensional problems. A further key advantage of the SPDE representation over current covariance-based approaches is its ease of extensibility to anisotropic and non-stationary random fields.

In applications, it may be necessary to consider fractional order SPDEs depending on the smoothness of the physical random field analysed. The numerical solution of fractional partial differential equations has been intensely studied recently, and various approximation techniques are available. In this paper, we used the rational series expansion of the fractional operator. Owing to the nonlocality of the fractional operator, the resulting precision matrices are typically dense. It is, however, possible to reformulate GP regression so that it is still sufficient to perform only sparse matrix operations. This reformulation is accomplished by introducing an auxiliary random field. The statFEM observation model assumed in this paper does not allow this reformulation, so that we considered for statFEM only non-fractional SPDEs. Nevertheless, it is possible to postulate an alternative observation model without a discrepancy term by assuming only epistemic uncertainties in terms of the choice of the model, which requires only sparse matrix operations.

In closing, we note several promising research directions and extensions of the proposed approach. In this paper, we did not provide any estimates for computational complexity and memory usage for evaluating the posterior mean and covariance in GP regression and statFEM. It appears straightforward to derive such estimates based on known standard estimates for sparse matrix operations, see, e.g.~\cite{ duff2017direct}. Furthermore, the flexibility of SPDE representation makes it appealing to other governing equations beyond the linear elliptic and shell equations considered. Especially, for engineering structures consisting of a combination of beams, trusses, solids and shells, see e.g.~\cite{febrianto2022digital}, the SPDE formulation seems to be the only principled approach for representing random fields. To this end, recent research in machine learning and statistics on the SPDE representation of random fields on graphs~\cite{ borovitskiy2021matern, nikitin2022non, bolin2023regularity} is noteworthy. Furthermore, in time-dependent problems, a spatio-temporal random field can be represented using a time-dependent SPDE, taking into account the correlations in space and time~\cite{nikitin2022non}. This is particularly appealing in connection with the extension of statFEM to time-dependent problems~\cite{duffin2021statistical} and, more broadly, Bayesian filtering and Kalman filtering for finite elements. Beyond GP regression and statFEM, the SPDE formulation is equally promising for representing and sampling generalised Mat\'ern random fields in large-scale stochastic forward and Bayesian inverse problems and their variants, see e.g. \cite{ranftl2022stochastic, nitzler2022generalized, rixner2021probabilistic, vadeboncoeur2023fully} amongst the extensive literature. Finally, the SPDE representation of random fields provides an interesting link between finite elements and the kernel methods for PDEs~\cite{chen2021solving, owhadi2015}, which use a chosen covariance function for discretisation. This covariance function may indeed come from the finite element discretisation of an SPDE.

%
\section*{Acknowledgements \label{sec:acknowledgements}}
%
This work was supported by Wave $1$ of The UKRI Strategic Priorities Fund under the EPSRC Grant EP/T001569/1, particularly the “Digital twins for complex engineering systems" theme within that grant, and The Alan Turing Institute.

\appendix
%
\section{Rational interpolation \label{sec:rationalInterp}}
%
The barycentric rational interpolation of the power function~$ r(x) = x^{\gamma} $, with $x, \, \gamma \in \mathbb R$, using its values at the~$m + 1$ distinct interpolation points with the coordinates, that is~$\{ x_i, \, r(x_i) \}_{i=1}^{m+1} $, reads 
\begin{equation} \label{eq:barycentRat}
	r( x)  \approx \frac{1}{ \sum_{i=1}^{m + 1} \dfrac{w_i}{ x - x_i}    }
			\sum_{i=1}^{m+1} \dfrac{w_i}{ x - x_i} r(x_i)
 \, .
\end{equation} 
Different choices for the weights~$w_i \neq 0$ yield different interpolation schemes passing through the~$m + 1$ points. For instance, the following weights yield the Lagrange interpolation
\begin{equation}
	w_i = \frac{1}{\prod_{k=1, k\neq i}^{m+1} (x_i - x_k)} \, ; 
\end{equation}
see~\cite{berrut2004barycentric} for a discussion on barycentric Lagrange interpolation. The flexibility in choosing the coordinates~$x_i$ and weights~$w_i$ can be used to interpolate~$r(x)$ at~$2m + 1$ rather than the original~$m + 1$ points to reduce the interpolation errors~\cite{knockaert2008simple}. In the BRASIL algorithm proposed by Hofreither~\cite{hofreither2021algorithm} both the interpolation points and weights are iteratively adjusted to minimise the interpolation errors. The respective open-source software implementation provides in addition to the weights~$w_i$ and the coordinates~$x_i$ also the roots of the numerator and denominator~$c_i$ and~$d_j$, respectively, so that the factorised form of~\eqref{eq:barycentRat} is given by 
\begin{equation}
	r( x)  \approx \frac{a_{m+1} \prod_{i=1}^{m}  (x- c_i   ) }{b_{m+1} \prod_{j=1}^{m} ( x- d_j  ) } \, .
\end{equation}
As an illustrative example, Figure~\ref{fig:rationalInterpolation} depicts the interpolation of $r(x) = \sqrt{x}$ using two rational interpolants with $m=2$ and $m=4$.
%
%
%
%
\section{Gaussian random vectors \label{sec:multivariateNormal}}
%
We summarise in the following the relevant properties of multivariate Gaussian random variables used throughout this paper. The corresponding proofs can be found, for instance, in~\cite[Ch.\ 2.3]{bishop2006pattern} and~\cite[Ch.\ 2.5]{murphy2012machine}.   A Gaussian random vector~$\vec u \in \mathbb R^{n_u}$ has the probability density 
\begin{equation}
	p(\vec u ) = \set N \left(\overline{\vec u}, \, \vec C_u \right) = \frac{1}{\sqrt{ (2 \pi)^{n_u} \det \left( \vec C_u \right) }} \left (  - \frac{1}{2} \left (\vec u - \overline{\vec u} \right )^\trans \vec C_u^{-1} \left (\vec u - \overline{\vec u} \right ) \right ) \, .
\end{equation}
With the help of the expectation operator~$\expect$ the respective mean vector and covariance matrix are defined as 
\begin{subequations}
\begin{align}
	\overline {\vec u} &= \expect \left [\vec u \right ]  \,, \\ 
	\vec C_u &= \cov( \vec u, \vec u ) = \expect \left[ \left( \vec u - \overline{\vec u} \right) \left( \vec u - \overline{\vec u} \right)^\trans \right]  \, . 
\end{align}
\end{subequations}
Furthermore, the precision matrix is defined as~$\vec Q_u= \vec C_u^{-1} $.
%

\subsection{Linear transformations \label{sec:multivariateNormalCov}}
%
A random vector~$\vec B \vec u + \vec b$ given as the linear transformation of the random vector~\mbox{$\vec u \sim \set N \left(\overline{\vec u}, \, \vec C_u \right) $} via the deterministic matrix~$\vec B$ and the deterministic vector~$\vec b$ has the Gaussian density
\begin{equation} \label{eq:linearTransformApp}
	\vec B \vec u + \vec b \sim  \set{N} \left(\vec B \, \overline{\vec u} + \vec b,  \vec B \vec C_u  \vec B^\trans \right) \, .
\end{equation}
This result is easy to verify considering the linearity of the expectation operator, i.e.,  
\begin{subequations}  \label{eq:appMappedMeanCov}
\begin{align} 
	 \expect \left [\vec B \vec u + \vec b \right ] &=  \vec B \expect \left [\vec u \right ] + \vec b = \vec B \, \overline{\vec u} + \vec b\, , \\ 
	 \cov( \vec B \vec u + \vec b  ,  \vec B \vec u + \vec b  ) &= \expect \left[ \left(  \vec B \vec u  - \vec B \, \overline{\vec u} \right) \left( \vec B \vec u - \vec B \, \overline{\vec u} \right)^\trans \right ] = \vec B  \expect \left[ \left(  \vec u  - \overline{\vec u} \right) \left(  \vec u -  \overline{\vec u} \right)^\trans \right] \vec B^\trans =  \vec B \vec C_u  \vec B^\trans \, .
\end{align}
\end{subequations}

Similarly, with the deterministic matrices~$\vec B$ and~$\vec L$ and the deterministic vector~$\vec b$ the linear combination $\vec B \vec u + \vec L \vec v + \vec b$ of two independent Gaussian random vectors~\mbox{$\vec u \sim  \set{N} \left(  \overline{\vec u}, \vec C_u \right)$}  and~\mbox{$\vec v \sim  \set{N} \left(  \overline{\vec v}, \vec C_v \right)$} has the Gaussian density
\begin{equation}
	\vec B \vec u + \vec L \vec v + \vec b \sim  \set{N} \left(\vec B \, \overline{\vec u} + \vec L  \overline{\vec v}  + \vec b,  \vec B \vec C_u  \vec B^\trans +  \vec L \vec C_v  \vec L^\trans \right) \, .
\end{equation}

%
\subsection{Jointly Gaussian random vectors  \label{sec:multivariateNormalConditional}}
%
The two jointly Gaussian random vectors~$\vec u \in \mathbb R^{n_u}$  and~$\vec v \in \mathbb R^{n_v}$ have the probability density
\begin{equation}
	\begin{pmatrix} \vec u \\ 
		\vec  v 
	\end{pmatrix} 
	\sim p(\vec u, \, \vec v) =  \set{N} \left( 
	\begin{pmatrix}
		\overline{\vec u} \\
		\overline{\vec v} 
	\end{pmatrix} \,,
	\begin{pmatrix}
		\vec C_{uu} & \vec C_{uv} \\
		\vec C_{vu} & \vec C_{vv}
	\end{pmatrix}
	\right)
	= 
	 \set{N} \left( 
	\begin{pmatrix}
		\overline{\vec u} \\
		\overline{\vec v} 
	\end{pmatrix} \,,
	\begin{pmatrix}
		\vec Q_{uu} & \vec Q_{uv} \\
		\vec Q_{vu} & \vec Q_{vv}
	\end{pmatrix}^{-1} 
	\right)\, . 
	\label{eq:GaussianDensityJoint}
\end{equation}
The covariance and precision matrices are decomposed such that~$\vec C_{uv}, \, \vec Q_{uv} \in \mathbb R^{n_u \times n_v}$,  and so on. The marginal distribution of~$\vec u$ is defined as 
\begin{equation}
	p (\vec u) = \int p(\vec u, \, \vec v)  \D \vec v  \,,
	\label{eq:GaussianDensityMarginalisation}
\end{equation}
and is given by~$p(\vec u) =  \set N\left(\overline{\vec u}, \, \vec C_{uu} \right) $,~$p(\vec v) =  \set N \left(\overline{\vec v}, \, \vec C_{vv} \right)$.

The conditional density of the random vector~$\vec u$ when the vector~$\vec v$ is given, i.e.\ is observed, is a  Gaussian given by  
\begin{equation} \label{eq:gaussianDensityConditionalFormal}
	p( \vec u \vert \vec v ) = \set{N}\left( \overline{\vec u}_{ \vert v}, \vec C_{u \vert v} \right) \, , 
\end{equation}
with the mean vector and covariance matrix 
\begin{subequations}\label{eq:gaussianDensityConditional}
	\begin{align}
		\overline{\vec u}_{\vert v} &= \overline{\vec u}  + \vec C_{uv} \vec C_{vv}^{-1} \left( \vec v - \overline{\vec v} \right) \,,
		\label{eq:gaussianDensityConditionalMean}
		\\
		\vec C_{ u \vert v} &= \vec C_{uu} - \vec C_{uv} \vec C_{vv}^{-1} \vec C_{vu} \, .
		\label{eq:gaussianDensityConditionalCovariance}
	\end{align}
\end{subequations}
Both expressions can be expressed in terms of precision matrices such that 
\begin{subequations}\label{eq:gaussianDensityConditionalQ}
	\begin{align}
		\overline{\vec u}_{\vert v} &= \overline{\vec u} - \vec Q_{uu}^{-1} \vec Q_{uv} \left( \vec v - \overline{\vec v}  \right) \,,
		\label{eq:gaussianDensityConditionalMeanQ}
		\\
		\vec C_{u \vert v} &= \vec Q_{uu}^{-1} \,.
		\label{eq:gaussianDensityConditionalCovarianceQ}
	\end{align}
\end{subequations}
%
%
\subsection{Sampling \label{sec:multivariateNormalSampling}}
%
We first express the random vector~$\vec u \sim \set{N} \left( \overline{\vec u}, \vec C_u \right) $ as the linear transformation of~$\vec t \sim \set{N}\left( \vec 0, \vec I \right)$ using the Cholesky decomposition~$\vec C_u = \vec L \vec L^\trans$. That is,  
\begin{equation}
	\vec u = \overline{\vec u} + \vec L \vec t \, .
\end{equation}
According to~\eqref{eq:appMappedMeanCov} the mean and covariance of the mapped vector are~$\expect(\vec u) = \overline{\vec u} $  and~$\cov( \vec u, \vec u )= \vec L \vec L^\trans = \vec C_u$ as desired. Hence, we can generate samples of~$\vec u$ using the Cholesky factor~$\vec L$ and sampling~$\vec t$ with standard algorithms available in most software libraries. 

To sample using the Cholesky decomposition of the precision matrix note that 
\begin{equation}
	\vec Q_u =  \vec C_u^{-1} =   \vec L^{-\trans} \vec L^{-1} \, .
\end{equation}
Accordingly, the affine transformation 
\begin{equation}
	\vec u = \overline{\vec u} + \vec L^{-\trans} \vec t \, .
\end{equation}
yields a random vector with the mean and covariance~$\expect(\vec u) = \overline{\vec u} $  and~$\cov( \vec u, \vec u )= \vec L^{-\trans} \vec L^{-1} = \vec C_u^{-1} = \vec Q_u$. 

To obtain multiple samples it is sufficient to factorise the covariance or precision matrix only once.  See also Rue and Held~\cite{rue2005gaussian} for sampling using precision matrices of manipulated forms.
\section{Repeated readings \label{sec:repeatedReadingsAppendix}}
%
In the following we derive the posterior density and log marginal likelihood for repeated readings, i.e.\ $n_o > 1$. In case of repeated readings there are at the~$\{ \vec x_i \}_{i=1}^{n_y}$  locations~$\{\vec y_j \}_{j=1}^{n_o}$  observations available, which are collected in the matrix~$\vec Y \in \mathbb R^{n_y \times n_o}$. We first consider the Gaussian process regression in Section \ref{sec:gpRegression}, and then the statistical finite element analysis in Section~\ref{sec:statFEM}. We adopt an inductive approach by first considering~$n_o = 2$ and then generalising the obtained formulas to~$n_o > 2$. Crucially, in both formulations, we show that the computational complexity for evaluating the posterior density and log marginal likelihood has the same complexity when~$n_o > 1$.
\subsection{Gaussian process regression \label{sec:repeatedReadingsAppendix2}}
\subsubsection*{Posterior density}
The joint probability density for $n_o = 2$ is given by
\begin{equation}
	\begin{pmatrix}
		\vec t \\ \vec y_1 \\ \vec y_2
	\end{pmatrix} \sim \set{N} \left(
	\begin{pmatrix}
		\vec 0 \\ \vec 0 \\ \vec 0 
	\end{pmatrix} \,,
	\begin{pmatrix}
		\vec Q_t^{-1} &  \vec Q_t^{-1} \vec F_r^\trans \vec P^\trans  &  \vec Q_t^{-1} \vec F_r^\trans \vec P^\trans
		\\ 
		\vec P \vec F_r \vec Q_t^{-1} & \vec P \vec F_r \vec Q_t^{-1}  \vec F_r^\trans \vec P^\trans + \sigma_e^2 \vec I & \vec P \vec F_r \vec Q_t^{-1}  \vec F_r^\trans \vec P^\trans
		\\
		\vec P \vec F_r \vec Q_t^{-1} & \vec P \vec F_r \vec Q_t^{-1}  \vec F_r^\trans \vec P^\trans & \vec P \vec F_r \vec Q_t^{-1}  \vec F_r^\trans \vec P^\trans + \sigma_e^2 \vec I
	\end{pmatrix}
	\right) \,.
\end{equation}
The covariance matrix can be partitioned according to
\begin{equation}
	\text{cov}(\vec t, \left(
	\vec y_1, \vec y_2
	\right) )
	=
	\left( \begin{array}{c|c}
	\vec C_{aa} & \vec C_{ab} \\ \hline
	\vec C_{ba} &  \vec C_{bb}
	\end{array}\right)
	=
	\left( \begin{array}{c|c}
		\vec Q_t^{-1}  &\begin{pmatrix} 
			\vec Q_t^{-1} \vec F_r^\trans \vec P^\trans &  \vec Q_t^{-1} \vec F_r^\trans \vec P^\trans
		\end{pmatrix}
		\\   \hline  \\ [-1em]
		\begin{pmatrix}
			\vec P \vec F_r \vec Q_t^{-1}  \\  \vec P \vec F_r \vec Q_t^{-1}
		\end{pmatrix} & \begin{pmatrix}
			\vec P \vec F_r \vec Q_t^{-1}  \vec F_r^\trans \vec P^\trans + \sigma_e^2 \vec I & \vec P \vec F_r \vec Q_t^{-1}  \vec F_r^\trans \vec P^\trans
			\\ \vec P \vec F_r \vec Q_t^{-1}  \vec F_r^\trans \vec P^\trans & \vec P \vec F_r \vec Q_t^{-1}  \vec F_r^\trans \vec P^\trans + \sigma_e^2 \vec I
		\end{pmatrix}
	\end{array}\right)\,.
	\label{eq:jointDensityQPartitionBB}
\end{equation}
The corresponding partitioned precision matrix reads
\begin{equation}
	\left( \begin{array}{c|c}
		\vec Q_{aa}  & \vec Q_{ab} \\ \hline
		\vec Q_{ba}  & \vec Q_{bb}
	\end{array}\right)
	=
	\left( \begin{array}{c|c}
	\vec C_{aa} & \vec C_{ab} \\ \hline
	\vec C_{ba} &  \vec C_{bb}
	\end{array}\right)^{-1} \,,
\end{equation}
where the relevant components are
\begin{subequations}
	\begin{align}
		\vec Q_{aa} &=  \vec C^{-1}_{aa} +  \vec C^{-1}_{aa}   \vec C_{ab} \left( \vec C_{bb} - \vec C_{ba} \vec C^{-1}_{aa} \vec C_{ab} \right)^{-1}  \vec C_{ba} \vec C^{-1}_{aa} 
		= \vec Q_t +  2 \sigma_e^{-2} \vec F_r^\trans \vec P^\trans \vec P \vec F_r  \,,
		\\
		\vec Q_{ab} & = - \vec C_{aa}^{-1} \vec C_{ab} \left( \vec C_{bb} - \vec C_{ba} \vec C^{-1}_{aa} \vec C_{ab} \right)^{-1} 
		= -  \begin{pmatrix}
			\vec F_r^\trans \vec P^\trans & \vec F_r^\trans \vec P^\trans  
		\end{pmatrix} \begin{pmatrix}
			\sigma_e^{-2} \vec I & \vec 0 
			\\
			\vec 0 &  \sigma_e^{-2} \vec I
		\end{pmatrix} \,.
	\end{align}
\end{subequations}
Hence, according to~\eqref{eq:gaussianDensityConditionalQ}, the conditional mean and covariance are given by
\begin{subequations}
	\begin{align}
		\overline{\vec t}_{\vert y_1, y_2} &= - \vec Q^{-1}_{aa} \vec Q_{ab}  
		\begin{pmatrix}
			\vec y_1 
			\\
			\vec y_2
		\end{pmatrix} 
		=  \left( \sigma_e^{2} \vec Q_t +  {n_o}  \vec F_r^\trans \vec P^\trans  \vec P \vec F_r \right)^{-1}  \vec F_r^\trans \vec P^\trans  \sum_{i=1}^{{2}}  \vec y_i \,,
		\\
		\vec Q^{-1}_{t\vert y_1, y_2} & =  \vec Q_{aa}^{-1} = \sigma_e^2 \left( \sigma_e^{2} \vec Q_t + {2}  \vec F_r^\trans \vec P^\trans  \vec P \vec F_r \right)^{-1}
		\label{eq:gpQCovarianceTGivenY}
		\,.
	\end{align}
\end{subequations}
Using induction we generalise the conditional mean and covariance for arbitrary~$n_o$ to
\begin{subequations}
	\begin{align}
		\overline{\vec t}_{\vert Y} &=
		\left( \sigma_e^{2} \vec Q_t +  {n_o}  \vec F_r^\trans \vec P^\trans  \vec P \vec F_r \right)^{-1}  \vec F_r^\trans \vec P^\trans  \sum_{i=1}^{{n_o}}  \vec y_i \,,
		\\
		\vec Q^{-1}_{t\vert Y} & = \sigma_e^2 \left( \sigma_e^{2} \vec Q_t + {n_o}  \vec F_r^\trans \vec P^\trans  \vec P \vec F_r \right)^{-1}
		\label{eq:gpQCovarianceTGivenYGeneral}
		\,.
	\end{align}
\end{subequations}
Finally, we recover the desired conditional density~$p(\vec s \vert \vec Y)$ with the mean and covariance
\begin{equation}
		\overline{\vec s}_{\vert Y} = \vec F_r \overline{\vec t}_{\vert Y} \,, \qquad 
		\vec Q^{-1}_{s\vert Y} = \vec F_r \vec Q^{-1}_{t\vert Y} \vec F_r^\trans \,.
\end{equation}
\subsubsection*{Log marginal likelihood}
The marginal likelihood for~$n_o = 2$ is given by
\begin{equation}
	p(\vec y_1, \vec y_2 ) = \mathcal{N} \left( 
	\begin{pmatrix}
		\vec 0
		\\
		\vec 0
	\end{pmatrix},
	\begin{pmatrix}
		\vec P \vec F_r \vec Q^{-1}_t\vec F_r^\trans \vec P^\trans + \sigma_e^2 \vec I & \vec P \vec F_r \vec Q^{-1}_t\vec F_r^\trans \vec P^\trans
		\\ \vec P \vec F_r \vec Q^{-1}_t\vec F_r^\trans \vec P^\trans & \vec P \vec F_r \vec Q^{-1}_t\vec F_r^\trans \vec P^\trans + \sigma_e^2 \vec I \end{pmatrix}
	\right) \,,
\end{equation}
and has the logarithm
\begin{equation}
	\log p(\vec y_1, \vec y_2) = - \frac{1}{2} 
	\begin{pmatrix}
		\vec y_1 
		\\
		\vec y_2
	\end{pmatrix}^\trans
	\vec C_{bb}^{-1}
	\begin{pmatrix}
		\vec y_1 
		\\
		\vec y_2
	\end{pmatrix} -  \frac{1}{2} \log \det \left( \vec C_{bb} \right)  -  {n_y} \log 2\pi \,,
\end{equation}
where~$\vec C_{bb}$, as defined in~\eqref{eq:jointDensityQPartitionBB}, can be decomposed as
\begin{equation}
	\vec C_{bb} = \begin{pmatrix}
		\vec P \vec F_r
		\\  \vec P \vec F_r  
	\end{pmatrix} \vec Q^{-1}_t
	\begin{pmatrix}
		\vec F_r^\trans \vec P^\trans
		& \vec F_r^\trans \vec P^\trans  
	\end{pmatrix} + 
	\begin{pmatrix}
		\sigma_e^2 \vec I & \vec 0 \\
		\vec 0   & \sigma_e^2 \vec I
	\end{pmatrix} \,.
\end{equation}
We use the Sherman--Morrison--Woodbury formula to evaluate the quadratic term
\begin{equation}
		-\frac{1}{2} \begin{pmatrix}
			\vec y_1 
			\\
			\vec y_2 
		\end{pmatrix}^\trans
		\vec C_{bb}^{-1}
		\begin{pmatrix}
			\vec y_1 
			\\
			\vec y_2
		\end{pmatrix} 
		= -\frac{1}{2} \sum_{j=1}^{{2}}  \vec y_j  ^\trans \left( \sigma_e^{-2}   \vec y_j   -   \sigma_e^{-4} \vec P \vec F_r \vec Q^{-1}_{t\vert y_1, y_2} \vec F_r^\trans \vec P^\trans   \sum\limits_{i=1}^{{2}} \vec y_i   \right) \,,
		\label{eq:gpLogMarginalLikelihoodQuadraticTerm}
\end{equation}
where~$\vec Q^{-1}_{t \vert y_1, y_2}$ is given by~\eqref{eq:gpQCovarianceTGivenY}. Similarly, we use the matrix determinant lemma to evaluate the determinant
\begin{equation}
		\det \left( \vec C_{bb} \right) =   \det \begin{pmatrix}
			\sigma_e^2 \vec I & \vec 0
			\\
			\vec 0  & \sigma_e^2 \vec I
		\end{pmatrix}  \det \left( {\vec Q_t + {2} \sigma_e^{-2} \vec F_r^\trans \vec P^\trans \vec P \vec F_r  }  \right)    \left( \det \left( \vec Q_t \right)   \right)^{-1}  \,.
	\label{eq:gpLogMarginalLikelihoodDetTerm}
\end{equation}
%
%
Using induction we generalise the log marginal likelihood for arbitrary $n_o$ to
\begin{equation}
	\log p(\vec Y) = M_1 - \frac{1}{2} \log M_2 -  \frac{{n_o} n_y}{2} \log 2\pi \,.
\end{equation}
By generalising~\eqref{eq:gpLogMarginalLikelihoodQuadraticTerm} we obtain for the quadratic term
\begin{equation}
	M_1
	= - \frac{1}{2} \sum_{j=1}^{{n_o}}  \vec y_j  ^\trans \left( \sigma_e^{-2}   \vec y_j   -   \sigma_e^{-4} \vec P \vec F_r \vec Q^{-1}_{t\vert Y} \vec F_r^\trans \vec P^\trans   \sum\limits_{i=1}^{{n_o}} \vec y_i   \right) \,,
\end{equation}
where~$\vec Q_{t \vert Y}$ is given by~\eqref{eq:gpQCovarianceTGivenYGeneral}. Similarly, by generalising~\eqref{eq:gpLogMarginalLikelihoodDetTerm} we obtain for the determinant
\begin{equation}
M_2 =  \left(\sigma_e^{2} \right)^{{n_o} n_y - n_u}   \det \left(  \sigma_e^{2} \vec Q_t + {n_o}  \vec F_r^\trans \vec P^\trans \vec P \vec F_r \right)    \left( \det \left( \vec Q_t \right)   \right)^{-1}  \,.
\end{equation}
\subsection{Statistical finite element method \label{sec:repeatedReadingsAppendix3}}
%
\subsubsection*{Posterior density}
The joint probability density for~$n_o = 2$ is given by
\begin{equation}
	\begin{pmatrix}
		\vec u \\ \vec y_1 \\ \vec y_2
	\end{pmatrix} \sim \set{N} \left(
	\begin{pmatrix}
		\overline{\vec u} \\  \vec P  \overline{\vec u} \\  \vec P  \overline{\vec u}
	\end{pmatrix} \,,
	\begin{pmatrix}
		\vec Q^{-1}_u &  \vec Q^{-1}_u \vec P^\trans &   \vec Q^{-1}_u \vec P^\trans
		\\ 
		\vec P  \vec Q^{-1}_u &  \vec P \vec Q^{-1}_u \vec P^\trans + \vec Q_{de}^{-1} &  \vec P \vec Q_u^{-1} \vec P^\trans
		\\
		\vec P  \vec Q^{-1}_u &  \vec P \vec Q_u^{-1} \vec P^\trans &  \vec P \vec Q^{-1}_u \vec P^\trans + \vec Q_{de}^{-1}
	\end{pmatrix}
	\right) \,.
\end{equation}
The covariance matrix can be partitioned according to
\begin{equation}
	\text{cov}(\vec u, \left(
	\vec y_1, \vec y_2
	\right) ) = 
	\left( \begin{array}{c|c}
		\vec C_{aa} &  \vec C_{ab}
		\\ \hline
		\vec C_{ba} &  \vec C_{bb}    
	\end{array}\right)
	=
	\left( \begin{array}{c|c}
		\vec Q^{-1}_u  & \begin{pmatrix}
			\vec Q^{-1}_u \vec P^\trans  &   \vec Q^{-1}_u \vec P^\trans 
		\end{pmatrix}
		\\ \hline \\ [-1em]
		\begin{pmatrix}
			\vec P  \vec Q^{-1}_u  \\   \vec P  \vec Q^{-1}_u
		\end{pmatrix}  &    \begin{pmatrix}
			\vec P \vec Q^{-1}_u \vec P^\trans + \vec Q_{de}^{-1}  &   \vec P \vec Q_u^{-1} \vec P^\trans
			\\  \vec P \vec Q_u^{-1} \vec P^\trans &  \vec P \vec Q^{-1}_u \vec P^\trans + \vec Q_{de}^{-1}
		\end{pmatrix}
	\end{array}\right)
	\,.
	\label{eq:jointDensityStatFEMQPartitionBB}
\end{equation}
The corresponding partitioned precision matrix reads
\begin{equation}
	\left( \begin{array}{c|c}
		\vec Q_{aa}  & \vec Q_{ab} \\ \hline
		\vec Q_{ba}  & \vec Q_{bb}
	\end{array}\right)
	=
	\left( \begin{array}{c|c}
		\vec C_{aa} & \vec C_{ab} \\ \hline
		\vec C_{ba} &  \vec C_{bb}
	\end{array}\right)^{-1} \,,
\end{equation}
where the relevant components are
\begin{subequations}
	\begin{align}
		\vec Q_{aa} &=  \vec C^{-1}_{aa} +  \vec C^{-1}_{aa}   \vec C_{ab} \left( \vec C_{bb} - \vec C_{ba} \vec C^{-1}_{aa} \vec C_{ab} \right)^{-1}  \vec C_{ba} \vec C^{-1}_{aa} 
		= \vec Q_u +  {2} \vec P^\trans \vec Q_{de} \vec P  \,,
		\\
		\vec Q_{ab} & = - \vec C_{aa}^{-1} \vec C_{ab} \left( \vec C_{bb} - \vec C_{ba} \vec C^{-1}_{aa} \vec C_{ab} \right)^{-1} 
		= -   \begin{pmatrix}
			\vec P^\trans & \vec P^\trans  
		\end{pmatrix} \begin{pmatrix}
			\vec Q_{de} & \vec 0 
			\\
			\vec 0 & \vec Q_{de} 
		\end{pmatrix} \,.
	\end{align}
\end{subequations}
Hence, according to~\eqref{eq:gaussianDensityConditionalCovarianceQ}, the conditional mean and covariance, respectively, are given by
\begin{subequations}
	\begin{align}
		\overline{\vec u}_{\vert y_1, y_2} 
		&= \overline{\vec u} + \left( \vec Q_u +   {2}  \vec P^\trans \vec Q_{de} \vec P \right)^{-1}   \vec P^\trans \vec Q_{de} \left( \sum_{i=1}^{{2}}  \vec y_i -  {2} \vec P \overline{\vec u} \right) \,,
		\\
		\vec Q^{-1}_{u\vert y_1, y_2} &  = \left( \vec Q_u +  {2}  \vec P^\trans \vec Q_{de}  \vec P \right)^{-1} 
		\,.
		\label{eq:statFECovarianceGivenY}
	\end{align}
\end{subequations}
Using induction we generalise the conditional mean and covariance for arbitrary~$n_o$ to
\begin{subequations}
	\begin{align}
		\overline{\vec u}_{\vert Y} 
		&= \overline{\vec u} + \left( \vec Q_u +   {n_o}  \vec P^\trans \vec Q_{de} \vec P \right)^{-1}   \vec P^\trans \vec Q_{de} \left( \sum_{i=1}^{{n_o}}  \vec y_i -  {n_o} \vec P \overline{\vec u} \right) \,,
		\\
		\vec Q^{-1}_{u\vert Y} &  = \left( \vec Q_u +  {n_o}  \vec P^\trans \vec Q_{de}  \vec P \right)^{-1} 
		\,.
		\label{eq:statFECovarianceGivenYGeneral}
	\end{align}
\end{subequations}
\subsubsection*{Log marginal likelihood}
The marginal likelihood for~$n_o = 2$ is given by
\begin{equation}
	\begin{split}
		p(\vec y_1, \vec y_2 ) & = \mathcal{N} \left( 
		\begin{pmatrix}
			\vec P \overline{\vec u}
			\\
			\vec P \overline{\vec u}
		\end{pmatrix},
		\begin{pmatrix}
			\vec P \vec Q^{-1}_u \vec P^\trans + \vec Q_{de}^{-1} &  \vec P \vec Q^{-1}_u \vec P^\trans
			\\  \vec P \vec Q^{-1}_u \vec P^\trans &  \vec P \vec Q^{-1}_u \vec P^\trans + \vec Q_{de}^{-1} \end{pmatrix}
		\right) \,,
	\end{split} 
\end{equation}
and has the logarithm
\begin{equation}
	\log p(\vec y_1, \vec y_2) = - \frac{1}{2} 
	\begin{pmatrix}
		\vec y_1 -  \vec P \overline{\vec u}
		\\
		\vec y_2 -  \vec P \overline{\vec u}
	\end{pmatrix}^\trans
	\vec C_{bb}^{-1}
	\begin{pmatrix}
		\vec y_1 -  \vec P \overline{\vec u}
		\\
		\vec y_2 -  \vec P \overline{\vec u}
	\end{pmatrix} - \frac{1}{2} \log \det \left( \vec C_{bb} \right) - n_y \log 2\pi \,,
\label{eq:logMarginalLikelihoodStatFEMultiple}
\end{equation}
where~$\vec C_{bb}$, as defined in~\eqref{eq:jointDensityStatFEMQPartitionBB}, can be decomposed as
\begin{equation}
	\vec C_{bb} = \begin{pmatrix}
		\vec P 
		\\  \vec P  
	\end{pmatrix} \vec Q^{-1}_u 
	\begin{pmatrix}
		\vec P^\trans
		& \vec P^\trans  
	\end{pmatrix} + 
	\begin{pmatrix}
		\vec Q_{de}^{-1} & \vec 0 \\
		\vec 0   & \vec Q_{de}^{-1}
	\end{pmatrix}  \,.
\end{equation}
We use the Sherman--Morrison--Woodbury formula to obtain the matrix inverse
\begin{equation}
	\vec C_{bb}^{-1} = 
	\begin{pmatrix}
		\vec Q_{de} & \vec 0
		\\
		\vec 0 & \vec Q_{de}
	\end{pmatrix} - \begin{pmatrix}
		\vec Q_{de} & \vec 0
		\\
		\vec 0 & \vec Q_{de}
	\end{pmatrix}		
	\begin{pmatrix}
		\vec P \vec Q^{-1}_{u\vert y_1, y_2} \vec P^\trans & \vec P \vec Q^{-1}_{u\vert y_1, y_2} \vec P^\trans
		\\
		\vec P \vec Q^{-1}_{u\vert y_1, y_2} \vec P^\trans & \vec P \vec Q^{-1}_{u\vert y_1, y_2} \vec P^\trans
	\end{pmatrix}
	\begin{pmatrix}
		\vec Q_{de} & \vec 0
		\\
		\vec 0 & \vec Q_{de}
	\end{pmatrix} \,.
\end{equation}
Hence, the quadratic term can be evaluated as
\begin{equation}
		-\frac{1}{2} \begin{pmatrix}
			\vec y_1 -  \vec P \overline{\vec u}
			\\
			\vec y_2 -  \vec P \overline{\vec u}
		\end{pmatrix}^\trans
		\vec C_{bb}^{-1}
		\begin{pmatrix}
			\vec y_1 -  \vec P \overline{\vec u}
			\\
			\vec y_2 -  \vec P \overline{\vec u}
		\end{pmatrix} 
		= -\frac{1}{2}\sum_{j=1}^{{2}} \left( \vec y_j -  \vec P \overline{\vec u} \right)^\trans \left( \vec Q_{de} \left(  \vec y_j -  \vec P \overline{\vec u} \right) -   \vec Q_{de} \vec P \vec Q^{-1}_{u\vert y_1, y_2} \vec P^\trans \vec Q_{de} \left( \sum\limits_{i=1}^{{2}} \vec y_i - {2}  \vec P \overline{\vec u} \right) \right) \,,
		\label{eq:statFELogMarginalLikelihoodQuadraticTerm}
\end{equation}
where~$\vec Q^{-1}_{u \vert y_1, y_2}$ is given by~\eqref{eq:statFECovarianceGivenY}. Similarly, we use the matrix determinant lemma to evaluate the determinant
\begin{equation}
	\det \left( \vec C_{bb} \right) = \det \left( \vec Q_u +  {2} {\vec P}^\trans   \vec Q_{de} \vec P\right) \left( \det \left(  {\vec Q_{de}^{-1}}  \right) \right)^{2} \left( \det \left( \vec Q_u  \right) \right)^{-1} \,.
	\label{eq:statFELogMarginalLikelihoodDetTerm}
\end{equation}
%
%
Using induction we generalise the log marginal likelihood for arbitrary $n_o$ to
\begin{equation}
	\log p(\vec Y) = M_1 - \frac{1}{2} \log M_2 -  \frac{{n_o} n_y}{2} \log 2\pi \,.
	\label{eq:logMarginalLikelihoodStatFEMultipleGeneral}
\end{equation}
By generalising~\eqref{eq:statFELogMarginalLikelihoodQuadraticTerm} we obtain for the quadratic term
\begin{equation}
	M_1 = -\frac{1}{2}\sum_{j=1}^{{n_o}} \left( \vec y_j -  \vec P \overline{\vec u} \right)^\trans \left( \vec Q_{de} \left(  \vec y_j -  \vec P \overline{\vec u} \right) -   \vec Q_{de} \vec P \vec Q^{-1}_{u\vert Y} \vec P^\trans \vec Q_{de} \left( \sum\limits_{i=1}^{{n_o}} \vec y_i - {n_o}  \vec P \overline{\vec u} \right) \right) \,,
\end{equation}
where~$\vec Q^{-1}_{u \vert Y}$ is given by~\eqref{eq:statFECovarianceGivenYGeneral}. Similarly, by generalising~\eqref{eq:statFELogMarginalLikelihoodDetTerm} we obtain for the determinant
\begin{equation}
	M_2 = \det \left( \vec Q_u +  {n_o} {\vec P}^\trans   \vec Q_{de} \vec P\right) \left( \det \left(  {\vec Q_{de}^{-1}}  \right) \right)^{n_o} \left( \det \left( \vec Q_u  \right) \right)^{-1} \,.
\end{equation}

\bibliographystyle{elsarticle-num-names}
\bibliography{statFEMsparse}

\end{document}